\documentclass[reqno,11pt]{amsart}
\usepackage{amssymb}
\usepackage{graphicx}
\usepackage{amsmath}
\usepackage{mathbbol}
\usepackage[margin=1.20in]{geometry}

\usepackage{ulem}

\usepackage[usenames, dvipsnames]{color}
\usepackage{hyperref}
\usepackage{verbatim}
\usepackage{mathrsfs}
\usepackage{bm}
\usepackage{color}
\usepackage{esint}
\allowdisplaybreaks[4]
\usepackage{dutchcal}

\numberwithin{equation}{section}

\newtheorem{theorem}{Theorem}[section]

\newtheorem{lemma}[theorem]{Lemma}
\newtheorem{prop}[theorem]{Proposition}

\theoremstyle{definition}
\newtheorem{remark}[theorem]{Remark}

\theoremstyle{definition}

\theoremstyle{definition}

\makeatletter
\def\dashint{\operatorname%
{\,\,\text{\bf-}\kern-.98em\DOTSI\intop\ilimits@\!\!}}
\makeatother

\def\\det{\text{\det}}

\def\.5{\frac{1}{2}}

\newcommand{\RN}[1]{%
\textup{\uppercase\expandafter{\romannumeral#1}}%
}

\renewcommand{\epsilon}{\varepsilon}

\newcounter{marnote}

\begin{document}

\title[Higher derivative estimates]{Higher derivative estimates for Stokes equations with closely spaced rigid inclusions in three dimensions}

\author[H.J. Dong]{Hongjie Dong}
\address[H.J. Dong]{Division of Applied Mathematics, Brown University, 182 
George Street, Providence, RI 02912, USA}
\email{Hongjie\_Dong@brown.edu}
\thanks{H. Dong was partially supported by the NSF under agreement DMS-2350129.}

\author[H.G. Li]{Haigang Li}
\address[H.G. Li]{School of Mathematical Sciences, Beijing Normal University, 
Laboratory of Mathematics and Complex Systems, Ministry of Education, Beijing 
100875, China.}
\email{hgli@bnu.edu.cn}
\thanks{H. Li was partially Supported by Beijing Natural Science Foundation (No. 
1242006), National Natural Science Foundation of China (No. 12471191), and the 
Fundamental Research Funds for the Central Universities (No. 2233200015).}

\author[H.J. Teng]{Huaijun Teng}
\address[H.J. Teng]{School of Mathematical Sciences, Beijing Normal University, 
Laboratory of Mathematics and Complex Systems, Ministry of Education, Beijing 
100875, China.}
\email{hjt2021@mail.bnu.edu.cn}

\author[P.H. Zhang]{Peihao Zhang}
\address[P.H. Zhang]{School of Mathematics and Statistics, Henan University, Kaifeng 475004, China.}
\email{phzhang@henu.edu.cn}


\date{\today} 



\begin{abstract}
In this paper, we establish higher-order derivative estimates for the Stokes equations in a three-dimensional domain containing two closely spaced rigid inclusions. We construct a sequence of auxiliary functions via an inductive process to isolate the leading singular terms of higher-order derivatives within the narrow region between the inclusions. For a class of convex inclusions of general shapes, the construction of three-dimensional auxiliary functions --- unlike the two-dimensional case --- relies on the decay properties of solutions to a class of two-dimensional partial differential equations with singular coefficients. Taking advantage of this, we obtain pointwise upper bounds of derivatives up to the seventh order for general inclusions. Under additional symmetry conditions, we derive optimal estimates for derivatives of arbitrary order. Consequently, we obtain precise blow-up rates for the Cauchy stress and its higher-order derivatives in the narrow region between the inclusions.
\end{abstract}

\maketitle

\section{Introduction}
In this paper, we investigate higher-order derivative estimates for the solution to the Stokes equations in three dimensions, in the presence of two adjacent rigid particles suspended in a fluid. Let $D \subset \mathbb{R}^3$ be a smooth, bounded domain, and let $D_1$ and $D_2$ be two subdomains of $D$ whose boundaries are of class $C^{m+1,\alpha}$, where $m \geq 1$ and $0 < \alpha < 1$. We assume that $D_1$ and $D_2$ are separated by a distance $\varepsilon$ and are located strictly inside $D$, away from the outer boundary $\partial D$. More precisely, we assume
$$
\overline{D}_1, \overline{D}_2 \subset D, \quad \varepsilon := \mathrm{dist}(D_1, D_2) > 0, \quad \mathrm{dist}(D_1 \cup D_2, \partial D) > \kappa_0 > 0,
$$
where $\kappa_0$ is a constant independent of $\varepsilon$. We also assume that the $C^{m+1,\alpha}$-norms of $\partial D_1$ and $\partial D_2$ are uniformly bounded independently of $\varepsilon$.

We consider the following Dirichlet problem for the Stokes system with two adjacent rigid inclusions $D_1$ and $D_2$:
\begin{align}\label{maineqs}
\begin{cases}
\mu\Delta{\bf u}=\nabla p,\quad\nabla\cdot{\bf u}=0\quad&\hbox{in}\ 
\Omega:=D\setminus\overline{D_1\cup D_2 },\\
{\bf u}|_{+}={\bf u}|_{-},&\hbox{on}\ \partial{D}_{i},\,i=1,2,\\
e({\bf u})=0,&\hbox{in}~~D_{i},\,i=1,2,\\ 
\int_{\partial{D}_{i}}\boldsymbol{\psi}_{\alpha}\cdot\frac{\partial {\bf u}}{\partial \nu}\Big|_{+}
\,ds-\int_{\partial{D}_{i}}p\boldsymbol{\psi}
_{\alpha}\cdot\nu\ ds =0,&i=1,2,\,\alpha=1,2,\dots,6,\\ 
{\bf u}=\boldsymbol{\varphi},&\hbox{on}\ \partial{D},
\end{cases}
\end{align}
where $\mu > 0$ is the viscosity, $e(\mathbf{u}) := \tfrac{1}{2}(\nabla \mathbf{u} + \nabla \mathbf{u}^T)$ is the strain tensor, and $\frac{\partial \mathbf{u}}{\partial \nu}\big|_{+} := \mu (\nabla \mathbf{u} + \nabla \mathbf{u}^T) \nu$, with $\nu$ denoting the unit outward normal to $D_i$ for $i = 1, 2$.
The functions $\{ \boldsymbol{\psi}_\alpha \}$ form a basis of the six-dimensional space of rigid displacements
$\Psi := \left\{ \boldsymbol{\psi} \in C^1(\mathbb{R}^3; \mathbb{R}^3) \;\middle|\; e(\boldsymbol{\psi}) = 0 \right\}$.
Specifically, the basis functions of $\Psi$ are given by
$$
\boldsymbol{\psi}_\alpha = \mathbf{e}_\alpha,\quad \alpha = 1, 2, 3, \quad \text{where } \{\mathbf{e}_\alpha,\alpha=1,2,3\} \text{ is the standard basis of } \mathbb{R}^3,$$
$$\boldsymbol{\psi}_4 = (x_2, -x_1, 0)^T,\quad \boldsymbol{\psi}_5 = (x_3, 0, -x_1)^T,\quad \boldsymbol{\psi}_6 = (0, x_3, -x_2)^T.~
$$
Throughout the paper, the subscripts $\pm$ denote the limits from the exterior and interior of the inclusions, respectively.
Since $D$ is a bounded domain, the divergence-free condition $\nabla \cdot \mathbf{u} = 0$, together with the divergence theorem, implies that the prescribed boundary velocity $\boldsymbol{\varphi}$ must satisfy the compatibility condition
$\int_{\partial D} \boldsymbol{\varphi} \cdot \nu = 0,$
to ensure the existence and uniqueness of a solution. See, for instance, \cite{Lady}. For a detailed treatment of the steady Stokes equations, we refer the interested reader to the monograph by Galdi \cite{GaldiBook}.

Investigating field enhancement in the narrow region between two adjacent inclusions is of fundamental importance in both fluid mechanics and materials science. The mathematical modeling of rigid inclusions moving in a fluid is one of the most classical and long-standing problems in fluid dynamics, tracing back to the foundational work of Stokes \cite{Stokes}, Kirchhoff \cite{Kirchhoff}, and Jeffery \cite{Jeffery}. Complex fluids, particularly those containing suspended particles or inclusions, often exhibit highly intricate flow behaviors. Therefore, understanding the interaction between inclusions and the surrounding fluid is crucial for practical applications involving particulate flows \cite{DDj}.

In experimental regimes with very low Reynolds numbers, the fluid motion is well approximated by the Stokes equations. In two dimensions, Ammari, Kang, Kim, and Yu \cite{AKKY} studied stress concentration in Stokes flow with two identical disk-shaped inclusions. Since the stress is represented by the gradient of the velocity field, they showed --- using bipolar coordinates --- that the gradient blows up at the rate of $\varepsilon^{-1/2}$ as the distance $\varepsilon$ between the inclusions vanishes. They also derived an asymptotic formula that explicitly characterizes the singular behavior of the solution.

However, as noted in \cite{AKKY} and emphasized by Kang in his ICM 2022 presentation \cite{K}, extending these results to general inclusions and higher dimensions remains a significant challenge. In this direction, Li and Xu recently developed a pointwise iteration method to obtain upper bounds for the gradient in the narrow gaps between inclusions with strictly convex boundaries \cite{LX1,LX2}. They demonstrated that the optimal blow-up rate is $\varepsilon^{-1/2}$ in two dimensions and $|\varepsilon \ln \varepsilon|^{-1}$ in three dimensions for general inclusions. Furthermore, they obtained pointwise estimates for second-order derivatives. The asymptotic formula of the concentrated stress for general inclusions has been further explored by Li, Zhang, and Zhang in \cite{LZZ}.

To better understand the stress enhancement in the narrow region, it is crucial --- both from engineering and computational mathematics perspectives \cite{LLZ} --- to investigate higher-order derivative estimates for solutions to the Stokes equations. In the two-dimensional Stokes flow setting, we previously extended the energy iteration technique developed in \cite{LX1,LX2} and introduced a novel method for constructing a sequence of auxiliary functions \cite{DLTZ1}. This enabled us to derive the optimal blow-up rate $\varepsilon^{-m/2}$ for the $m$-th order derivatives of the solution under certain symmetry conditions on the domain. For the case of general inclusions, we also established pointwise upper bounds for higher-order derivatives.

In the two-dimensional setting, the auxiliary function was chosen as a special solution to a first-order ordinary differential equation tailored to satisfy the necessary structural conditions. However, this approach cannot be directly extended to higher dimensions, as finding a suitable solution to a first-order partial differential equation in $(d-1)$ variables with the required properties is no longer feasible. Therefore, new analytical tools must be developed.

As a continuation of \cite{DLTZ1}, the present paper aims to construct a sequence of appropriate auxiliary functions for the three-dimensional case. These constructions are intricately connected to the decay properties of solutions to a class of two-dimensional inhomogeneous second-order partial differential equations with singular coefficients. We precisely quantify the singular behavior of higher-order derivatives in the Stokes system with closely spaced rigid inclusions in three dimensions.

This work is also partially motivated by the study of physical field enhancement between two hard inclusions in high-contrast materials, due to its close connection with linear elasticity. When two inclusions are closely spaced within a material, physical fields --- such as stress --- can become arbitrarily large in the narrow gap between them. The mathematical investigation of this phenomenon dates back to the seminal work of Babuška \cite{Bab}. It is now well established that such field blow-up occurs in both electrostatics and elastostatics. A quantitative understanding of this behavior is essential, as it plays a critical role in applications such as light confinement in electrostatics and failure analysis in elastostatics.

Over the past two decades, substantial progress has been made by numerous researchers in characterizing field enhancement. For studies related to the conductivity problem. See, for example, \cite{BV,AKL,AKL3,BLY,BT1,KangLY,BT2,LiYang,Ben,DLY1,LLL}. For investigations on linear elasticity with hard inclusions, we refer the reader to \cite{BLL,BLL2,KY,Li2021,LX24,DLTZ}. For results concerning the Stokes system or elliptic equations with general discontinuous coefficients and their applications, we refer the reader to \cite{LIL07, DK18, KW18, CK19, KMW21,LV,LN} and the references therein.

Before presenting our main results, we first describe the domain setup and introduce the necessary notation. For convenience, we translate and rotate the coordinate system so that the points
$P_1 = (0', \varepsilon/2) \in \partial D_1$ and $P_2 = (0', -\varepsilon/2) \in \partial D_2$
satisfy $\mathrm{dist}(P_1, P_2) = \mathrm{dist}(\partial D_1, \partial D_2) = \varepsilon$. There exists a constant $R > 0$, independent of $\varepsilon$, such that the portions of $\partial D_1$ and $\partial D_2$ near the origin can be represented as graphs of functions:
$$
x_3 = \frac{\varepsilon}{2} + h_1(x') \quad \text{and} \quad x_3 = -\frac{\varepsilon}{2} - h_2(x') \quad \text{for } |x'| \leq 2R,
$$
where we denote points in $\mathbb{R}^3$ by $x = (x', x_3)$, with $x' \in \mathbb{R}^2$.
 Here \( h_1, h_2 
\in C^{m+1, \alpha}(B_{2R}(0'))\), \(\alpha\in(0,1)\) are {\it radial} functions and satisfy the following conditions, for $ 2 
\leq k \leq m+1$, $i = 1, 2,$
\begin{equation}\label{h1h2prop}
\begin{aligned}
&-\frac{\varepsilon}{2} - h_2(x') < \frac{\varepsilon}{2} + h_1(x')\quad\text{for } |x'| \leq 2R,\\
&h_1(0') = h_2(0') = 0, \quad \nabla_{x'}h_1(0') = \nabla_{x'}h_2(0') = 0, \\
&h_1(x')+h_2(x')=\sum_{j=1}^{m}a_j|x'|^{j+1}+O(|x'|^{m+1+\alpha})\quad \text{for } |x'| \leq 2R,\\
&|h_i(x')| \leq C|x'|^2, \quad |\nabla_{x'}h_i(x')| \leq C|x'|, \quad |\nabla_{x'}^kh_i(x')| \leq C,
\end{aligned}
\end{equation}
where $a_1> 0$ and $a_j$, $j=2,\dots,m$, are constants. For \( 0 \leq r \leq 2R 
\), we define the neck region between \( D_1 \) and \( D_2 \) as
\[\Omega_r := \left\{ (x', x_3) \in \Omega : -\frac{\varepsilon}{2} - h_2(x') < 
x_3 < \frac{\varepsilon}{2} + h_1(x'), \, |x'| < r \right\}.\]
The top and bottom boundaries of the neck region are denoted by
$\Gamma^+_r := \big\{ (x', x_3) \in \Omega : x_3 = \frac{\varepsilon}{2} + h_1(x'), \, |x'| < r \big\}$ and
$\Gamma^-_r := \big\{ (x', x_3) \in \Omega : x_3 = -\frac{\varepsilon}{2} - h_2(x'), \, |x'| < r \big\}$. Denote the vertical distance at \(x'\) in the narrow region between \(D_1\) and \(D_2\) by
\[\delta(x') := \varepsilon + h_1(x') + h_2(x'), \quad |x'| \leq 2R.\]
It follows from \eqref{h1h2prop} that 
$\frac{1}{C}(\varepsilon+|x'|^2)\le \delta(x') \le C (\varepsilon+|x'|^2)$.
For $|z'|\leq R$ and $s<R$, set
\begin{equation}\label{omegasr}
\Omega_{s}(z'):=\left\{(x',x_{3})\big| -\frac{\varepsilon}{2}-h_{2}(x')<x_{3}
<\frac{\varepsilon}{2}+h_{1}(x'),\,|x'-z'|<s \right\}.
\end{equation}
Throughout this paper, we use $C$ to denote a generic constant, whose value may change from line to line. Unless otherwise specified, $C$ depends only on $\kappa_0$, $\kappa$, $R$, the boundary data $\boldsymbol{\varphi}$, and the upper bounds of the $C^{m+1,\alpha}$-norms of $\partial D_1$ and $\partial D_2$, but is independent of the small parameter $\varepsilon$.

\subsection{Upper bounds of the higher derivatives}

According to the standard theory for Stokes equations, we have the following bound for $m\ge 0$: 
$$\|\nabla^{m+1} {\bf u}\|_{L^{\infty}(\Omega\setminus\Omega_{R})} + \|\nabla^{m} (p-p(z))\|_{L^{\infty}(\Omega\setminus\Omega_{R})} \leq C,$$ 
for some fixed point $z\in\Omega\setminus\Omega_{R}$. See, for instance, \cite{GaldiBook}. Our primary focus is to obtain higher-order derivative estimates for problem \eqref{maineqs} in the narrow region $\Omega_{R}$.
 
To this end, following the idea in \cite{DLTZ1} for the two-dimensional problem and instead of solving an ordinary differential equation, we construct a sequence of divergence-free auxiliary functions which are, however, directly connected to the following inhomogeneous elliptic partial differential equation in two dimensions:
\begin{equation}\label{zyfc/}
\mathcal{L}\mathcal{U}:=\Delta_{x'} \mathcal{U}(x')+\frac{3}{\varepsilon+\mathcal{h}(x')}\nabla_{x'}\mathcal{h}(x')\cdot\nabla_{x'}\mathcal{U}(x')=\mathcal{F}(x')\quad \text{in}~B_{2R}(0')
\end{equation}
with given $\mathcal{F}(x')$, where $\mathcal{h}(x')$ is a positive function satisfying the third line in \eqref{h1h2prop}. In our setting, we will later take $\mathcal{h}(x')=h_{1}(x')+h_{2}(x')$.

If $\mathcal{F}(x')$ is radial, i.e., $\mathcal{F}(x')=f(|x'|)$, then a special solution to \eqref{zyfc/} can be obtained directly by ordinary differential equation theory --- similar to the approach for the two-dimensional Stokes problem in \cite{DLTZ1}. However, when $\mathcal{F}(x')$ is not radial --- even for example, $\mathcal{F}(x')=f(|x'|)x_{1}$ --- an explicit form of $\mathcal{U}(x')$ is generally not be easily obtained. Instead to find a special solution to \eqref{zyfc/}, we establish the following estimates.

\begin{prop}\label{U_esti}
Suppose that $\mathcal{F}$ is smooth. Let $\mathcal{U}\in C^{m+1}(B_{2R}(0'))$ be a solution to 
\eqref{zyfc/}. If 
\begin{equation}\label{festi}
|\nabla_{x'}^l\mathcal{F}(x')|\leq C(\varepsilon+\mathcal{h}(x'))^{\gamma-\frac l2}\quad\mbox{for}~~ 0\le l\le m,\quad x'\in B_{2R}(0')
\end{equation}
for some constants $C$, then for $-3\le \gamma<-1$, we have
\begin{align}\label{highestofu}
|\nabla_{x'}^l \mathcal{U}(x')|\leq C(\varepsilon+\mathcal{h}(x'))^{\gamma-\frac l2+1}\quad\mbox{for}~~ 0\le l\le m+1,
\quad x'\in B_{2R}(0').
\end{align}
\end{prop}

Here, we conjecture that \eqref{highestofu} remains valid even when $\gamma \ge -1$ and in all dimensions. With the aid of Proposition \ref{U_esti}, we are able to derive higher-order derivative estimates for the solution to \eqref{maineqs} --- up to the seventh order --- for inclusions of various shapes, by constructing appropriate auxiliary functions.

\begin{theorem}\label{main thm1}
Let \(D_{1}\), \(D_{2}\), \(D\), and \(\Omega\) be defined as above. Let \({\bf u} 
\in H^{1}(D; \mathbb{R}^3) \cap C^{m+1}(\bar{\Omega}; \mathbb{R}^3)\) and 
\(p \in L^2(D) \cap C^{m}(\bar{\Omega})\) be the solution to 
\eqref{maineqs}, with \(\boldsymbol{\varphi} \in C^{m+1,\alpha}(\partial D; 
\mathbb{R}^3)\), for an integer \(m \geq 0\) and \(0 < \alpha < 1\). Then, 
for sufficiently small \(0<\varepsilon < 1/2\) and $x=(x',x_3)\in\Omega_{R}
$, we have the following estimates: 
\begin{equation}\label{upper_esti1}
\quad\text{(i)}~
|\nabla {\bf u}(x)| \leq \frac{C(1+|\ln\varepsilon||x'|)}{|\ln\varepsilon|\delta(x')}, \quad
|p(x) - p(z', 0)| \leq \frac{C}{|\ln\varepsilon|\delta(x')^{3/2}}+C
\quad \text{for}~x\in \Omega_{R},
\end{equation}
for some fixed point \((z', 0) \in \Omega_{R}\) with \(|z'| = R/2\); and

(ii) for \(1\le m \le 6\),
\begin{equation}\label{upper_esti2}
|\nabla^{m+1}{\bf u}(x)| + |\nabla^m p(x)|\leq 
\frac{C}{|\ln\varepsilon|}\delta(x')^{-\frac{m+3}{2}}+C{\delta(x')^{-\frac{m+1}{2}}} \quad \text{for}~x\in \Omega_{R}.
\end{equation}
\end{theorem}

The core of the proof for Theorem \ref{main thm1} involves constructing a sequence of auxiliary functions to establish pointwise higher-order derivative estimates for solutions $({\bf u}_{i}^{\alpha},p_{i}
^{\alpha})$ to the following boundary value problem with narrow region 
\begin{equation}\label{u,peq1}
\begin{cases}
\mu\Delta{\bf u}_{i}^\alpha=\nabla p_{i}^{\alpha},~~~\nabla\cdot 
{\bf u}_{i}^{\alpha}=0&\mathrm{in}~\Omega,\\
{\bf u}_{i}^{\alpha}={\boldsymbol\psi}_{\alpha}&\mathrm{on}~\partial{D}
_{i},\\
{\bf u}_{i}^{\alpha}=0&\mathrm{on}~\partial{D_{j}}\cup\partial{D},~j\neq 
i,
\end{cases}~~\text{for}~\alpha=1,2,\dots,6,~i=1,2.
\end{equation}

We use $({\bf u}_{1}^{1},p_{1}^{1})$ as an example to illustrate our method. As mentioned previously, we construct auxiliary functions satisfying \eqref{zyfc/} to suppress the growth of the inhomogeneous term near the origin, thereby isolating the singularities in higher-order derivatives of the solution in the narrow region. 

First, we choose a divergence-free function ${\bf v}_{1}^{1}$ such that ${\bf v}_{1}^{1}
={\bf u}_{1}^{1}$ on $\partial\Omega$ and a suitable $\bar{p}_1$. The difference
${\bf w}:={\bf u}_{1}^{1}-{\bf v}_{1}^{1}$ and $q:=p_{1}^{1}-\bar{p}_{1}$ then
satisfy the zero boundary value problem 
\begin{equation}\label{c1w11}
\begin{cases}
-\mu\Delta{\bf w}+\nabla q={\bf f}_{1}^{1},\quad\nabla\cdot {\bf w}=0&\mathrm{in}\ \Omega,\\
{\bf w}=0\quad&\mathrm{on}\ \partial\Omega,
\end{cases}
\end{equation}
where ${\bf f}_{1}^{1}:=\mu\Delta{\bf v}_{1}^{1}-\nabla\bar{p}_{1}$. By Proposition \ref{keyprop2}, if $|{\bf f}_{1}^{1}|\leq C\delta(x')^{-1}$ in the narrow region $\Omega_{2R}$, it follows that
$|\nabla {\bf u}_1^1|\leq |\nabla{\bf v}_{1}^{1}|+C.$ To achieve this, recalling that ${\boldsymbol\psi}_{1}=(1,0,0)^{T}$, we first define the auxiliary function
\begin{equation*}
\tilde{\bf v}_{1}^{1}=\boldsymbol\psi_{1}\Big(k(x)+\frac{1}{2}\Big),\quad\mbox{where}~k(x)=\frac{x_3-\frac12(h_1-h_2)(x')}{\delta(x')}\quad\hbox{in}\ \Omega_{2R},
\end{equation*}
for $h_{1}(x')$ and $h_{2}(x')$ satisfying \eqref{h1h2prop}. It is easy to check that $k(x)+\frac{1}{2}=1$ on $\Gamma_{2R}^+$ and $k(x)+\frac{1}{2}=0$ 
on $\Gamma_{2R}^-$, and so $\tilde{\bf v}_{1}^{1}
={\bf u}_{1}^{1}$ on $\Gamma^{\pm}_{2R}$, but it fails to satisfy $\nabla\cdot \tilde{\bf v}_{1}^{1}=0$, since $\nabla\cdot \tilde{\bf v}_{1}^{1}=\partial_{x_1}k(x)\neq0$. We therefore need to introduce a correction term
$$\hat{\bf v}_{1}^{1}=\Big(F_{1}(x'),G_{1}(x'),H_{1}(x)\Big)^{T}\big(k(x)^2-\frac14\big)\quad\hbox{in}\ \Omega_{2R},$$
vanishing on $\Gamma^{+}_{2R}$ (using $k(x)^2-\frac14=0\,~\text{on}\,\Gamma_{2R}^+
\cup\Gamma_{2R}^-$) and ensuring $\nabla\cdot (\tilde{\bf v}_{1}^{1}+\hat{\bf v}_{1}^{1})=0$. A sufficient condition is that $F_{1}(x')$, $G_{1}(x')$, and $H_{1}(x)$ satisfy
\begin{equation}\label{c12by1}
\partial_{x_{1}}F_{1}(x')+\partial_{x_{2}}G_{1}(x')+\frac{\partial_{x_{1}}\delta(x')}{\delta(x')}F_{1}(x')+\frac{\partial_{x_{2}}\delta(x')}{\delta(x')}G_{1}(x')=-\frac{3\partial_{x_{1}}(h_1-h_2)(x')}{\delta(x')}
\end{equation}
and
$$H_{1}(x)=\partial_{x_{1}}(h_1-h_2)(x')k(x)+\frac12\partial_{x_{1}}(h_1+h_2)(x')-
\delta(x')\big(F_{1}(x)\partial_{x_1}k(x)-G_{1}(x)\partial_{x_2}k(x)\big).
$$
See \eqref{ffz2} and \eqref{2by1} below. This construction results in $\partial_{x_3x_3}(\hat{\bf v}_{1}^{1})^{(1)}$ and $\partial_{x_3x_3}(\hat{\bf v}_{1}^{1})^{(2)}$ being the leading terms in $\mu\Delta(\tilde{\bf v}_{1}^{1}+\hat{\bf v}_{1}^{1})$, exhibiting a higher order $\delta(x')^{-2}$ singularity. To reduce this singularity, we introduce 
$\hat{p}_{1}(x)=\hat{p}_{1}(x')$ to cancel these terms, namely,
\begin{equation}\label{c1tj1}
\partial_{x_{1}}\hat{p}_{1}(x')=\mu\partial_{x_3x_3}({\bf v}_1^1)^{(1)}=\frac{2\mu F(x')}{\delta(x')^2},\quad
\partial_{x_{2}}\hat{p}_{1}(x')=\mu\partial_{x_3x_3}({\bf v}_1^1)^{(2)}=\frac{2\mu G(x')}{\delta(x')^2}.
\end{equation}
Thus, substituting \eqref{c1tj1} into \eqref{c12by1} yields
$$
\Delta_{x'} \hat{p}_{1}(x')+\frac{3}{\delta(x')}
\nabla_{x'} \delta(x')\cdot\nabla_{x'}\hat{p}_{1}(x')=-\frac{6\mu\partial_{x_{1}}(h_1-h_2)(x')}
{\delta(x')^3}:=\mathcal{F}_1(x').
$$
By Proposition \ref{U_esti}, $\hat{p}_{1}(x')$ satisfies \eqref{zyfc/} with $\mathcal{h}=h_{1}+h_{2}$ and $\gamma=-\frac52$, ensuring the existence of auxiliary functions ${\bf v}_1^1$ and $\hat{p}_1^1$ with desired estimates. Then, the leading term in $\mu\Delta{\bf v}_{1}^{1}-\nabla\hat{p}_{1}$ only appears in the third term, $\mu\partial_{x_{3}x_3}({\bf v}_1^1)^{(3)}$, which is a polynomial in $x_3$. Thus, by choosing a polynomial $\tilde{p}_1$ in $x_3$, we eliminate this term, guaranteeing $|{\bf f}_{1}^{1}|\leq C\delta(x')^{-1}$, and thus by Proposition \ref{keyprop2}, $|\nabla {\bf u}_1^1|\leq |\nabla{\bf v}_{1}^{1}|+C.$ 

To estimate the second-order derivatives $\nabla^2{\bf u}_1^1$, we further construct auxiliary functions $({\bf v}_1^2, \bar{p}_2)$ to reduce the upper bounds of the inhomogeneous term in \eqref{c1w11}. We divide ${\bf v}_1^2$ into two parts, ${\bf v}_1^2=\tilde{\bf v}_1^2+\hat{\bf v}_1^2$. The first part, $\tilde{\bf v}_1^2$, cancels the leading term of ${\bf f}_1^1$, while the second part, $\hat{\bf v}_1^2$, ensures $\nabla\cdot{\bf v}_1^2=0$. Both parts vanish on $\Gamma_{2R}^{\pm}$. From the first step, we know that the leading terms in ${\bf f}_1^1$ are $({\bf f}_{1}^{1})^{(1)}$ and $({\bf f}_{1}^{1})^{(2)}$. Therefore, we choose $\tilde{\bf v}_1^2$ such that
$$\partial_{x_3x_3}(\tilde{\bf v}_1^2)^{(i)}=({\bf f}_{1}^{1})^{(i)},~~i=1,2.$$
Since $({\bf f}_{1}^{1})^{(i)}$ is expressible as a polynomial in $x_3$, $(\tilde{\bf v}_1^2)^{(i)}$ can be solved explicitly. See \eqref{u11-v12-1,2} and \eqref{F12i22iform} below. If we simply set $(\tilde{\bf v}_1^2)^{(3)}$ = 0 such that $\nabla\cdot \tilde{\bf v}_1^2=R(x',x_3)$, for some explicit function $R(x',x_3)$, then it is challenging to find a suitable $\hat{\bf v}_1^2$ satisfying $\nabla\cdot\hat{\bf v}_1^2=-R(x',x_3)$ in this case. Here, a key idea is to first design $(\tilde{\bf v}_1^2)^{(3)}$ such that $\nabla\cdot \tilde{\bf v}_1^2=R(x')$, independent of $x_3$. Then we construct $\hat{\bf v}_1^2$ in the form
$$\hat{\bf v}_{1}^{2}=(F_{2}(x'),G_{2}(x'),H_{2}(x))^{T}\big(k(x)^2-\frac14\big)\quad\hbox{in}\ \Omega_{2R},$$
with $\nabla\cdot\hat{\bf v}_{1}^{2}=-R(x')$. Here, we take
\begin{equation*}
H_{2}(x)=-2\delta(x')R(x')k(x)-\delta(x')\partial_{x_1}k(x)F_{2}(x')-\delta(x')\partial_{x_2}k(x)G_{2}(x')
\end{equation*}
and
\begin{equation*}
\partial_{x_{1}}F_{2}(x')+\partial_{x_{2}}G_{2}(x')+
\frac{\partial_{x_{1}}\delta(x')}{\delta(x')}F_{2}(x')+
\frac{\partial_{x_{2}}\delta(x')}{\delta(x')}G_{2}(x')=6R(x').
\end{equation*}
See \eqref{tilF32eqs} and \eqref{TilFb2andR} below. 

Analogous to the first step, the functions $F_{2}(x')$ and $G_{2}(x')$ introduce new singularities in ${\bf f}_1^2$. To cancel these leading singularities, we choose $\hat{p}_2$ such that
\begin{equation*}
\partial_{x_1}\hat{p}_{2}=\frac{2\mu F_{2}(x')}{\delta(x')^2},\quad
\partial_{x_2}\hat{p}_{2}=\frac{2\mu G_{2}(x')}{\delta(x')^2}.
\end{equation*}
Then
\begin{equation*}
\Delta_{x'}\hat{p}_{2}+\frac{3}{\delta(x')}\nabla_{x'}\delta(x')\cdot
\nabla_{x'}\hat{p}_{2}=\frac{12\mu R(x')}{\delta(x')^2}:=\mathcal{F}_2(x').
\end{equation*}
One can see that $\hat{p}_{2}(x')$ satisfies \eqref{zyfc/} with $\mathcal{h}=h_{1}+h_{1}$ and $\gamma=-\frac32$. Hence, by virtue of Proposition \ref{U_esti}, there exist auxiliary functions ${\bf v}_1^2$ and $\hat{p}_2$ with desired decays. We choose $\tilde{p}_2$ as a polynomial in $x_3$ to eliminate the singularity in the third term; see \eqref{tilp12} below. Therefore, $|{\bf f}_{1}^{2}|=|\mu\Delta({\bf v}_{1}^{1}+{\bf v}_{1}^{2})-\nabla(\bar{p}_{1}+\bar{p}_{2})|\leq C$, which implies $|\nabla^2 {\bf u}_1^1|\leq |\nabla^2({\bf v}_{1}^{1}+{\bf v}_{1}^{2})|+C$ by Proposition \ref{keyprop2}.

However, we cannot further proceed this procedure for the third and higher order derivatives, as it remains unclear whether Proposition \ref{U_esti} holds for $\gamma\ge -1$.
Therefore, it is no longer possible to find a solution to \eqref{zyfc/} to construct an auxiliary divergence-free function ${\bf v}_1^3$. Instead, we introduce a non-divergence-free ${\bf v}_1^3$ and pressure ${\bar p}_3$ (see \eqref{v13-12}, \eqref{v13-3form}, and \eqref{barp3form}) such that
\begin{equation}\label{c1m=61}
\Big|\nabla^{m+1}({\bf v}_1^1+{\bf v}_1^2+{\bf v}_1^3)\Big|+\Big|\nabla^{m}({\bar p}_1+{\bar p}_2+{\bar p}_3)\Big|\le\,C\delta(x')^{-\frac{m+3}{2}},\quad~ m\geq1,
\end{equation}
and 
\begin{equation}\label{c1m=62}
\Big\|\nabla^{m+1}({\bf u}_1^1 -\sum_{i=1}^{3}{\bf v}_1^i)\Big\|_{L^{\infty}(\Omega_{\delta(x')/
2}(x'))} +\Big\|\nabla^m(p_1^1-\sum_{i=1}^{3}{\bar p}_i)\Big\|_{L^{\infty}(\Omega_{\delta(x')/2}
(x'))} \le \,C\delta(x')^{\frac{3}
{2}-m}.
\end{equation}
By comparing the decay orders on the right-hand sides \eqref{c1m=61} and \eqref{c1m=62}, we observe that when $m\geq 7$, the orders of $\delta(x')$, $m-\frac32>\frac{m+3}{2}$. This means the estimates can no longer capture the leading behavior of $\nabla^{8}{\bf u}_1^1$ and $\nabla^{7}p_1^1$. Therefore, we can only establish the derivative estimates up to seventh-order for the general inclusion case. 

For derivatives of order $8$ or higher, we obtain a comparatively larger upper bound:
\begin{equation*}
|\nabla^{m+1}{\bf u}_1^1(x)| + |\nabla^m p_1^1(x)|\le C\delta(x')^{\frac32-m}
,\quad m\geq 7, \quad \text{for}~x\in \Omega_{R}.
\end{equation*}
An open question is whether \eqref{highestofu} holds for $\gamma\ge1$. If confirmed, iterative application of our method would extend \eqref{upper_esti2} to all $m\ge1$. The detailed proofs are carried out component-wise for the decomposition, and the construction process is performed separately for each component $({\bf u}_{i}^{\alpha},p_{i}^{\alpha})$ in Sections \ref{sec_estu11,12}, \ref{sec4}, and \ref{sec5}, and Theorem \ref{main thm1} is then proved in Section \ref{sec6}.

However, by incorporating symmetry assumptions on both the domain and boundary data, we can directly obtain the optimal estimates. Specifically, we assume that $D_1 \cup D_2$ and $D$ are symmetric about the coordinate axes,  the boundary data satisfy $\boldsymbol{\varphi}(-x) = -\boldsymbol{\varphi}(x)$, and $h_1(x') = h_2(x') = \frac{1}{2} h(x')$ in $\Omega_{2R}$. This setting clearly includes the most important case in practical applications, where $D_1$ and $D_2$ are two balls of equal radius.

\begin{theorem}\label{main thm2}
Let \(D_{1}\), \(D_{2}\), \(D\), and \(\Omega\) be defined as above. Let \({\bf u} 
\in H^{1}(D; \mathbb{R}^3) \cap C^{m+1}(\bar{\Omega}; \mathbb{R}^3)\) and \(p 
\in L^2(D) \cap C^{m}(\bar{\Omega})\) be the solution to \eqref{maineqs}, 
where \(\boldsymbol{\varphi} \in C^{m+1,\alpha}(\partial D; \mathbb{R}^3)\)
for an integer \(m \geq 0\) and \(0 < \alpha < 1\), and \(\boldsymbol{\varphi}(-
x)= -\boldsymbol{\varphi}(x)\). Then, for sufficiently small \(0 <\varepsilon < 
1/2\) and $x = (x',x_3)\in\Omega_{R}$, the following estimates hold:
\begin{equation*}
\text{(i)}~~~~
|\nabla {\bf u}(x)| \leq\frac{C}{|\ln\varepsilon|
\delta(x')}+C, \quad
|p(x) - p(z', 
0)| \leq 
\frac{C\varepsilon}{\delta(x')^2}+C
~\text{for}~(z', 0) \in \Omega_{R}~\text{with}~ |z'| = R/2;
\end{equation*}
and 
\begin{equation*}
\text{(ii)}~~~~|\nabla^{m+1}{\bf u}(x)| + |\nabla^m p(x)| \leq 
C\delta(x')^{-\frac{m+2}{2}}
 +C\quad\mbox{for}~m\ge 1.\qquad\qquad
\end{equation*}
As a consequence, for the Cauchy stress tensor \(\sigma[{\bf u}, p] = 2\mu 
e({\bf u}) - p\mathbb{I}\), where \(\mathbb{I}\) is the identity matrix, we have
\begin{equation*}
\big|\nabla^m \sigma[{\bf u}, p - p(z', 0)]\big| \leq 
C\delta(x')^{-\frac{m+2}{2}}
 +C. 
\end{equation*}
\end{theorem}

\begin{remark}
We would like to point out that our method works well in higher dimensions $d\geq4$. The construction of the auxiliary functions is associated with the decaying solutions of $(d-1)$-dimensional partial differential equations \eqref{zyfc/} and the method we utilize to prove Proposition \ref{U_esti} is valid for any $d\ge3$.
\end{remark}

\subsection{Lower bounds of the higher derivatives}
We will establish a lower bound for $|\nabla^{m+1} {\bf u}(x)|$ that exhibits the same blow-up rate as stated in Theorem \ref{main thm2}, thereby confirming the optimality of the derived rates. This result is obtained under the assumption that $D_1 \cup D_2$ and $D$ are symmetric with respect to the $x_i$-axis for $i = 1, 2, 3$, and that both $h_1(x')$ and $h_2(x')$ are quadratic functions.
Specifically, we set \( h_1(x') = h_2(x')=\frac{1}{2}| x'|^2 \) for \( |x'| \leq 2R \). We define a linear continuous functional of 
${\boldsymbol\varphi}$:
\begin{equation*}
\tilde b_{j}^{*\alpha}[{\boldsymbol\varphi}]:=\int_{\partial D_j^0}
{\boldsymbol\psi}_\alpha\cdot\sigma[{\bf u}^*,p^*]\nu,\quad\alpha=1,2,\dots,6,~j=1,2,
\end{equation*}
where $({\bf u}^*,p^*)$ verifies
\begin{align*}
\begin{cases}
\mu\Delta{\bf u}^{*}=\nabla p^{*},\quad\nabla\cdot {\bf u}^{*}
=0\quad&\hbox{in}\ \Omega^{0},\\
{\bf u}^{*}=\sum_{\alpha=1}^{6}C_{*}^{\alpha}{\boldsymbol\psi}_{\alpha}
&\hbox{on}\ \partial D_{1}^{0}\cup\partial D_{2}^{0},\\
\int_{\partial{D}_{1}^{0}}{\boldsymbol\psi}_\alpha\cdot\sigma[{\bf u}
^*,p^{*}]\nu+\int_{\partial{D}_{2}^{0}}{\boldsymbol\psi}
_\alpha\cdot\sigma[{\bf 
u}^*,p^{*}]\nu=0,&\alpha=1,2,\dots,6,\\
{\bf u}^{*}={\boldsymbol\varphi}&\hbox{on}\ \partial{D},
 \end{cases}
\end{align*}
where $D_{1}^{0}:=\{x\in\mathbb R^{3}~\big|~ x+P_{1}\in D_{1}\},\quad~ D_{2}
^{0}:=\{x\in\mathbb R^{3}~\big|~ x+P_{2}\in D_{2}\}$ and $
\Omega^{0}:=D\setminus\overline{D_{1}^{0}\cup D_{2}^{0}}$. 

\begin{theorem}\label{main thm3}
Let $D_1,D_2\subset D$ be defined as above and let ${\bf u}\in H^{1}(D;\mathbb 
R^{3})\cap C^{m}(\bar{\Omega};\mathbb R^{3})$ be the solution to \eqref{maineqs}. 
Suppose that there exists a $\boldsymbol{\varphi}$ such that $
\boldsymbol{\varphi}(-x)=-\boldsymbol{\varphi}(x)$ and $\tilde b_{1}^{*1}
[\boldsymbol{\varphi}],$ $\tilde b_{1}^{*3}[\boldsymbol{\varphi}]\neq 0$. Then for 
$m\ge 1$ and sufficiently small $0<\varepsilon<1/2$, there exists a small 
constant $r>0$ which may depend on $m$, such that
\begin{align}\label{lower_esti}
	|\partial_{x_{3}}{\bf u}^{(1)}(r\sqrt{\varepsilon},0,0)|
	\ge \frac{C|\tilde{b}_{1}^{*1}[\boldsymbol{\varphi}]|}{|\ln\varepsilon|}
	\varepsilon^{-1},\quad |\nabla_{x'}^{m-1}\partial_{x_3}^{2}({\bf u})^{(1)}
	(r\sqrt{\varepsilon},0,0)|\ge C |\tilde{b}_{1}^{*3}[\boldsymbol{\varphi}]|
	\varepsilon^{-\frac{m+2}{2}}.
\end{align}
Consequently,
$$\big|\nabla_{x'}^{m-1}\partial_{x_3}\sigma[{\bf u},p]\big|
(r\sqrt{\varepsilon},0,0)\ge C|\tilde{b}_{1}^{*3}[\boldsymbol{\varphi}]|
\varepsilon^{-\frac{m+2}{2}}.$$
\end{theorem}
\begin{remark}
The upper bound of $\sigma[{\bf u},p-p(z',0)]$, as established in Theorem \ref{main thm2}, is of order $\frac{1}{\varepsilon}$. This blow-up rate is considered optimal, as it satisfies 
$$\big|\sigma[{\bf u},p-p(z',0)]\big|(0',x_3)\ge \frac{1}{C\varepsilon},$$
as shown in \cite[Theorem 4.6]{LX2}.

The optimality of the gradient estimates in \eqref{upper_esti1} is confirmed by the lower-bound \eqref{lower_esti} with \(m=0\). The optimality of the higher-order derivative bounds in Theorem \ref{main thm2} is established through Theorem \ref{main thm3}. 
However, for arbitrary convex functions \(h_{1}(x')\) and \(h_{2}(x')\), the optimality of the high-order derivative estimates in Theorem \ref{main thm1} remains unclear. It is speculated that the asymmetry \(h_{1}(x') \neq h_{2}(x')\) could lead to additional singularities, as observed in the two-dimensional case \cite{DLTZ1}.
\end{remark}

The rest of this paper is organized as follows. In Section \ref{sec_pri}, we first present a general $W^{m,\infty}$ estimate for Stokes equations with partially zero boundary conditions in a narrow region. The proof, given in the Appendix, is based on the $W^{m,p}$ estimates and a modified version of the energy method. Specifically, for an auxiliary function pair \(({\bf v},\bar{p})\) satisfying the same boundary conditions as \(({\bf u},p)\) on \(\Gamma^{\pm}_{2R}\), with \(\nabla \cdot {\bf v} = g(x)\) in \(\Omega_{2R}\) for an explicit function $g(x)$, we reduce the local estimates of the differences \( |\nabla^{m+1}({\bf u} - {\bf v})|\) and \(|\nabla^m(p - \bar{p})| \) to controlling the right-hand side terms \({\bf f} := \mu \Delta {\bf v}(x) - \nabla \bar{p}(x)\) and $g(x)$. Proposition \ref{U_esti} is proven in Subsection \ref{subsec22}. In Section \ref{sec_estu11,12}, we consider the Stokes system with prescribed Dirichlet boundary data \(\boldsymbol{\psi}_{\alpha}\) on \(\partial{D}_{1}\) and \(\partial{D}_{2}\), for \(\alpha = 1, 2\) (see \eqref{u,peq1} below). We demonstrate that if $h_1(x')=h_2(x')$, then the non-homogeneous term in equation \eqref{zyfc/} becomes radial. We then establish $m$-order derivative estimates for any $m$ by constructing a sequence of divergence-free auxiliary functions \({\bf v}_\alpha^m\) and corresponding \(\bar{p}_m\) in \(\Omega_{2R}\). This construction ensures that the quantity 
\[
\Big|\sum_{l=1}^{m+1} (\mu \Delta {\bf v}_{\alpha}^{l} - \nabla \bar{p}_{l})(x)\Big|
\] 
satisfies the required pointwise upper bound in \(\Omega_{2R}\), as stated in Proposition \ref{keyprop1}. However, in cases where $h_1(x')\neq h_2(x')$, we are only able to construct three pairs of auxiliary functions \(({\bf v}_\alpha^l, \bar{p}_l)\), \(l=1,2,3\), by utilizing Proposition \ref{U_esti}. This allows us to accurately capture the leading terms of the high-order derivative of equation \eqref{maineqs} up to the seventh order. For derivatives of order eight or higher, we are only able to provide relatively coarse upper bounds. In Sections \ref{sec4} and \ref{sec5}, we investigate the estimates of high-order derivatives for specific Dirichlet boundary data \(\boldsymbol{\psi}_{\alpha}\) on \(\partial{D}_{1}\) and \(\partial{D}_{2}\), where $\alpha=3,4,5,6$, respectively. Finally, in Section \ref{sec6}, we prove Theorems \ref{main thm1}, \ref{main thm2}, and \ref{main thm3}, based on the estimates established in Sections \ref{sec_estu11,12}, \ref{sec4}, and \ref{sec5}.

\section{Preliminary results}\label{sec_pri}
In this section, we present high-order derivatives for a class of homogeneous Stokes equations and prove Proposition \ref{U_esti}. Let \(D\) be a domain in \(\mathbb{R}^d\), \(d\ge 2\), with two closely disjoint subdomains \(D_1\) and \(D_2\), and denote \(\Omega = D \setminus\overline{D_1 \cup D_2}\). We consider the following Stokes equations:
\begin{equation}\label{ueqs}
\begin{cases}
-\mu\Delta{\bf u}+\nabla p=0\quad&\mathrm{in}\,\,\Omega,\\
\nabla\cdot {\bf u}=0\quad&\mathrm{in}\,\,\Omega,\\
{\bf u}=\boldsymbol{\varphi}\quad&\mathrm{on}\,\,\partial\Omega,
\end{cases}
\end{equation}
with given smooth boundary data \(\boldsymbol{\varphi}\) satisfying \( 
\int_{\partial \Omega} \boldsymbol{\varphi}\cdot\nu = 0 \). The existence 
and uniqueness of the solution to \eqref{ueqs} are well known. 
By standard regularity theory, we have
\[ \|{\bf u}\|_{C^{m+1}(\Omega \setminus \Omega_R)} + \|p\|_{C^m(\Omega 
\setminus \Omega_R)} \leq C.
\]
We only need to estimate the high-order derivatives of \eqref{ueqs} in 
the narrow region $\Omega_{2R}$. 

Our approach involves constructing a series of different auxiliary functions to capture the singular behavior of the derivatives of any order of the solution to \eqref{ueqs}. Let \(({ \bf v}, \bar{p})\) be a known function pair that satisfies \({\bf v} = {\bf u} = {\bf g}\) on \(\Gamma^+_{2R} \cup \Gamma^-_{2R}\) and \(\nabla \cdot {\bf v} = g(x)\) in \(\Omega_{2R}\), with
\[\|{\bf v}\|_{C^{m+1}(\Omega \setminus \Omega_R)} + \|\bar{p}\|_{C^m(\Omega \setminus \Omega_R)} \leq C.\]
Define the difference between \(({ \bf u}, p)\) and \(({ \bf v}, \bar{p})\) as
\[{\bf w} := {\bf u} - {\bf v} \quad \text{and} \quad q := p - \bar{p}.\]
Then, \(({ \bf w}, q)\) satisfies the following boundary value problem for the nonhomogeneous Stokes equations in the narrow region:
\begin{equation}\label{weqs1}
\begin{cases}
-\mu\Delta{\bf w}+\nabla q={\bf f}\quad&\mathrm{in}\,\,\Omega_{2R},\\
\nabla\cdot {\bf w}=-g\quad&\mathrm{in}\,\,\Omega_{2R},\\
{\bf w}=0\quad&\mathrm{on}\,\,\Gamma^+_{2R}\cup\Gamma^-_{2R},
\end{cases}
\end{equation}
where \({\bf f} := \mu \Delta {\bf v} - \nabla \bar{p}\) and \(g=\nabla \cdot {\bf v}\). We have the following higher derivative estimates for \eqref{weqs1}.

\subsection{Estimates for Stokes equations with non-zero divergence in the narrow region}

\begin{prop}\label{keyprop1}
Let ${\bf w}\in H^{1}(\Omega_{\frac{3}{2}R})\cap L^{\infty}(\Omega_{2R}
\setminus\Omega_{R})$ and $p\in L^{\infty}(\Omega_{2R}\setminus\Omega_{R})$ 
be the solution to \eqref{weqs1}. For any $m\ge1$, if
\begin{equation*}
|\nabla^k{\bf f}(x)|\leq 
C\delta(x')^{l-k}, \quad\forall\,0\le k\leq m,~l\ge -\frac{3}{2}\quad\mbox{for}
~x=(x',x_d)\in\Omega_{2R},
\end{equation*}
and
\begin{equation*}
|\nabla^kg(x)|\leq 
C\delta(x')^{\alpha-k}, \quad\forall\,0\le k\leq m+1,~\alpha\ge -\frac{1}{2}\quad\mbox{for}
~x=(x',x_d)\in\Omega_{2R},
\end{equation*}
then it holds that
\begin{equation*}
\|\nabla{\bf w}\|_{L^{\infty}(\Omega_{\delta(z')/2}(z'))}\leq 
C\delta(z')^{l+1}+C\delta(z')^{\alpha},\quad z=(z',z_d) \in \Omega_{R},
\end{equation*}
and for $z =(z',z_d) \in \Omega_{R}$,
\begin{equation*}
\|\nabla^{m+1}{\bf w}\|_{L^{\infty}(\Omega_{\delta(z')/2}(z'))}+\|\nabla^{m}q\|_{L^{\infty}(\Omega_{\delta(z')/2}(z'))}\le 
 C\delta(z')^{l+1-m}+C\delta(z')^{\alpha-m},
\end{equation*}
where $\Omega_{s}(z')$ is defined in \eqref{omegasr}.
\end{prop}

The proof can be found in the Appendix. It can be deduced from Proposition \ref{keyprop1} that in order to obtain higher derivative estimates of \({\bf u}\) and \(p\) in $\Omega_{R}$, it is sufficient to construct suitable auxiliary functions \({\bf v}\) and \(\bar{p}\) in the region \(\Omega_{2R}\), such that \({\bf f} \) and $g$ and their derivatives have appropriate upper bounds. This is the central idea of the paper and the main challenge also lies here. In the construction process of auxiliary divergence-free functions, Proposition \ref{U_esti} plays an important role.
 
\subsection{Proof of Proposition \ref{U_esti}}\label{subsec22}

In this subsection, we give the proof of Proposition \ref{U_esti}. To this end, the following lemma is needed.
\begin{lemma}\label{delta-center}
For any $\alpha\in\mathbb{R},$ there exists a universal constant $C_{\alpha}$, 
such that 
\begin{equation}\label{dela alpha esti}
\|(\varepsilon+\mathcal{h})^{\alpha}\|_{L^{\infty}\big(B_{\frac{1}{C_{\alpha}}
\sqrt{\varepsilon+\mathcal{h}(x'_{0})}}(x'_{0})\big)} \le C_{\alpha}(\varepsilon+\mathcal{h}(x'_{0}))^{\alpha}, 
\quad x'_{0}\in B_{2R}(0').
\end{equation}
In particular, there exists a constant $\mathcal{C}:=C_{-1/2}$ such that 
\begin{equation}\label{partdelta/del esti}
\Big\|\frac{\nabla_{x'}\mathcal{h}}{\varepsilon+\mathcal{h}}\Big\|_{L^{\infty}\big(B_{\frac{1}
{\mathcal{C}} \sqrt{\varepsilon+\mathcal{h}(x'_{0})}}(x'_{0})\big)} \le \mathcal{C}
(\varepsilon+\mathcal{h}(x'_{0}))^{-\frac12}, \quad i=1,2,\quad x'_{0}\in B_{2R}(0').
\end{equation}
\end{lemma}

\begin{proof}[Proof of Lemma \ref{delta-center}]
It follows from the third line in \eqref{h1h2prop} that $\frac{1}{C}|x'|^2\leq\mathcal{h}(x')\leq C|x'|^2$, thus we only need to prove \eqref{dela alpha esti} for $\mathcal{h}(x')=|x'|^2$. For any $x'\in B_{\sqrt{\varepsilon+\mathcal{h}(x'_{0})}}(x'_{0})$, if $\alpha\ge 0$, by a direct calculation, we have
\begin{equation*}
(\varepsilon+\mathcal{h}(x'))^\alpha\leq C\Big(\varepsilon+\Big(|x'_{0}|+\sqrt{\varepsilon+|x_0'|
^2}\Big)^2 \Big)^{\alpha}\le C(\varepsilon+\mathcal{h}(x_0'))^\alpha.
\end{equation*}
Thus, \eqref{dela alpha esti} holds for $\alpha\ge 0$. If $\alpha<0$, we divide the proof into two cases.

Case 1: if $|x'_{0}|\le \sqrt{\varepsilon}$, for any $x'\in 
B_{\sqrt{\varepsilon+\mathcal{h}(x'_{0})}}(x'_{0})$, we have
\begin{equation*}
(\varepsilon+\mathcal{h}(x'))^\alpha\leq \varepsilon^\alpha\leq 2^{-\alpha}(\varepsilon+|x_0'|
^2)^{\alpha}\leq C (\varepsilon+\mathcal{h}(x_0'))^{\alpha}.
\end{equation*}

Case 2: if $|x'_{0}|> \sqrt{\varepsilon}$, there exists a small enough 
$r<\frac14$ such that $\big||x_0'|-r\sqrt{\varepsilon+\mathcal{h}(x_0')}\big|>0$, for any 
$x'\in B_{r\sqrt{\varepsilon+\mathcal{h}(x'_{0})}}(x'_{0})$, we have
\begin{align*}
(\varepsilon+\mathcal{h}(x'))^{\alpha}&\,\le\Big(\varepsilon+|x_0'|^2+r^2(\varepsilon+\mathcal{h}(x_0'))-2r|x_0'|
\sqrt{\varepsilon+\mathcal{h}(x_0')}\Big)^\alpha\\
&\,\le C \Big(\varepsilon+(1-2\sqrt{2}r)|x_0'|
^2\Big)^\alpha\le C(\varepsilon+\mathcal{h}(x_0'))^{\alpha}.
\end{align*}
We thus complete the proof.
\end{proof}

The following well-known iteration lemma is needed for our proof, which is from e.g. \cite[Chapter V, Lemma 3.1]{M1}.
\begin{lemma}\label{itLemma}
 Let $f(t)$ be a nonnegative bounded function defined in $[\tau_{0},\tau_{1}]$, $\tau_{0}\geq 0$. Suppose that for $\tau_{0}\leq t<s\leq\tau_{1}$, we have
 $$f(t)\leq [A(s-t)^{-a}+B]+\theta f(s),$$
 where $A,B, a,\theta$ are nonnegative constants with $0\leq \theta<1$. Then for all$\tau_{0}\leq \rho<R\le \tau_{1}$ there holds
 $$f(\rho)\leq\,c\,[A(R-\rho)^{-a}+B],$$
 where $c$ is a constant depending on $a$ and $\theta$.
\end{lemma}

With Lemmas \ref{delta-center} and \ref{itLemma} at hand, we now prove Proposition \ref{U_esti}.
\begin{proof}[Proof of Proposition \ref{U_esti}]Let $\mathcal{U}_{1}:=\mathcal{U}_1(r)$ and 
$\mathcal{U}_{2}:=\mathcal{U}_2(r)$, $r=|x'|$, be the solution to the following equation
\begin{equation}\label{u1u2eqs}
\mathcal{L}\mathcal{U}_{1}=C(\varepsilon+\mathcal{h}(x'))^{\gamma},\quad \mathcal{L}\mathcal{U}_{2}=-C(\varepsilon+\mathcal{h}(x'))^{\gamma} \quad\text{in}~B_{2R}(0').
\end{equation}
The boundary condition will be determined later. By \eqref{zyfc/}, we rewrite \eqref{u1u2eqs} as 
\begin{equation}\label{rad form}
\mathcal{U}_1''+\Big(\frac{1}{r}+\frac{3\delta(r)'}{\delta(r)}\Big)\mathcal{U}_1'=C(\varepsilon+\mathcal{h}(r))^{\gamma}, 
\,\, \mathcal{U}_2''+\Big(\frac{1}{r}+\frac{3\delta(r)'}{\delta(r)}\Big)\mathcal{U}_2'=-
C(\varepsilon+\mathcal{h}(r))^{\gamma} ~\text{in}~B_{2R}(0').
\end{equation}
We first prove \eqref{highestofu} for $l=0$. By 
\eqref{rad form}, we choose $U_1$ and $U_2$ in the following form,
\begin{equation*}
\mathcal{U}_{1}(x')=
-C\int_{r}^{2R}\frac{1}{s(\varepsilon+\mathcal{h}(s))^3}\int_{0}^{s}(\varepsilon+\mathcal{h}(t))^{\gamma+3}t \, dt \ ds=-\mathcal{U}_2(x').
\end{equation*}
Adirect calculation gives
\begin{equation*}
\Big|\int_{r}^{2R}\frac{1}{s(\varepsilon+\mathcal{h}(s))^3}\int_{0}^{s}(\varepsilon+\mathcal{h}(t))^{\gamma+3}t \, dt \ ds\Big|\le \Big|\int_{r}^{2R}(\varepsilon+\mathcal{h}(s))^\gamma s \ ds\Big|\le C( \varepsilon+\mathcal{h}(x'))^{\gamma+1},
\end{equation*}
which implies
\begin{equation}\label{u1u2aes}
 |\mathcal{U}_1|,|\mathcal{U}_2|\leq C(\varepsilon+\mathcal{h}(x'))^{\gamma+1}\quad \text{in}~B_{2R}(0'). 
\end{equation}

Let $\mathcal{U}$ be the solution to \eqref{zyfc/} such that $\mathcal{U}(x')=0$, $|x'|=2R$. Notice 
that $\mathcal{U}_1(2R)=\mathcal{U}_2(2R)=0$, by comparison principle and \eqref{u1u2aes}, we have
\begin{equation}\label{uesti}
 |\mathcal{U}(x')|\le \max\{|\mathcal{U}_{1}|,|\mathcal{U}_2|\}\le C (\varepsilon+\mathcal{h}(x'))^{\gamma+1}\quad\text{for}~x\in 
 B_{2R}(0').
\end{equation}

For $l=1$, we rewrite \eqref{zyfc/} as 
\begin{equation}\label{lapu=F}
\Delta \mathcal{U} =-\frac{3}{\varepsilon+\mathcal{h}(x')}\nabla_{x'}\mathcal{h}(x')\nabla_{x'}\mathcal{U}(x')+\mathcal{F}(x'):=F(x') \quad\text{in}~B_{2R}(0').
\end{equation}
For any $y'=(y_{1},y_{2})\in B_{2R}(0')$, we denote the cube centered at $y'$ with radius $d$ as follows $$Q_{d}(y')=\{(x_1,x_2) : |x_1-y_1|<d,|x_2-y_2|<d \}.$$ 
By gradient estimate for the Possion equation (see, for instance, \cite[(3.15)]{GT}), we have 
\begin{align}\label{pointesti grad u}
|\nabla_{x'}\mathcal{U}(y')|\le \frac{4}{d}\sup_{\partial Q_{d}(y')}|\mathcal{U}|+ d\sup_{Q_{d}
(y')}|F|, \quad y'\in B_{2R}(0').
\end{align}
For any fixed point$x'_{0}\in B_{2R}(0'),\,0<t<s\le\frac{1}{2(C_{a}+C_{a+1}+
\mathcal{C})}<1$, assume $|\nabla_{x'}\mathcal{U}|$ attains its maximum $
\sup_{Q_{t\sqrt{\varepsilon+\mathcal{h}(x'_{0})} }(x'_{0})} |\nabla_{x'}\mathcal{U}|$ at $z'_{0}
\in \overline{Q_{t\sqrt{\varepsilon+\mathcal{h}(x'_{0})}}(x'_0)}\subset 
Q_{s\sqrt{\varepsilon+\mathcal{h}(x'_{0})}}(x'_{0})$. 
By using \eqref{lapu=F} and 
\eqref{pointesti grad u} with $d=\frac{s-t}{10}\sqrt{\varepsilon+\mathcal{h}(x'_{0})}$, we obtain
\begin{align}\label{supDutdelta}
&\sup_{Q_{t\sqrt{\varepsilon+\mathcal{h}(x'_{0})} }(x'_{0})}|\nabla_{x'}\mathcal{U}|=|\nabla_{x'}\mathcal{U}(z'_{0})| \nonumber\\
&\le \frac{40}{(s-t)\sqrt{\varepsilon+\mathcal{h}(x'_{0})}}\sup_{\partial Q_{\frac{s-
t}{10}\sqrt{\varepsilon+\mathcal{h}(x'_{0})}}(z'_{0})}|\mathcal{U}|+\frac{(s-t)\sqrt{\varepsilon+\mathcal{h}(x'_{0})}}{10}
\sup_{Q_{\frac{s-t}{10}\sqrt{\varepsilon+\mathcal{h}(x'_{0})}}(z'_{0})}|F| \nonumber\\
&\le \frac{40}{(s-t)\sqrt{\varepsilon+\mathcal{h}(x'_{0})}}
\sup_{Q_{s\sqrt{\varepsilon+\mathcal{h}(x'_{0})}}(x'_{0})}|\mathcal{U}|+\frac{(s-t)\sqrt{\varepsilon+\mathcal{h}(x'_{0})}}
{10}\sup_{Q_{s\sqrt{\varepsilon+\mathcal{h}(x'_{0})}}(x'_{0})}|F|\nonumber\\
&\le \frac{40}{(s-t)\sqrt{\varepsilon+\mathcal{h}(x'_{0})}}\sup_{Q_{s\sqrt{\varepsilon+\mathcal{h}(x'_{0})}}
(x'_{0})}|\mathcal{U}| +\frac{3(s-t)\sqrt{\varepsilon+\mathcal{h}(x'_{0})}}{10}\Big\|
\frac{\nabla_{x'}\mathcal{h}} {\varepsilon+\mathcal{h}}\Big\|_{L^{\infty}
(Q_{s\sqrt{\varepsilon+\mathcal{h}(x'_{0})}}(x'_{0}))} \nonumber\\ 
&\quad\cdot\sup_{Q_{s\sqrt{\varepsilon+\mathcal{h}(x'_{0})}}(x'_{0})}|\nabla_{x'}\mathcal{U}|+\frac{(s-t)
\sqrt{\varepsilon+\mathcal{h}(x'_{0})}}{10}\sup_{Q_{s\sqrt{\varepsilon+\mathcal{h}(x'_{0})}}(x'_{0})}|\mathcal{F}|.
\end{align} 
By \eqref{festi}, \eqref{uesti}, and Lemma \ref{delta-center}, we have 
\begin{equation}\label{supuesti}
\sup_{Q_{s\sqrt{\varepsilon+\mathcal{h}(x'_{0})}}(x'_{0})}|\mathcal{U}|\le\|(\varepsilon+\mathcal{h})^{\gamma+1}\|_{L^{\infty}
(Q_{s\sqrt{\varepsilon+\mathcal{h}(x'_{0})}}(x'_{0}))}\le C(\varepsilon+\mathcal{h}(x'_{0}))^{\gamma+1}
\end{equation}
and 
\begin{equation}\label{supfest}
\sup_{Q_{s\sqrt{\varepsilon+\mathcal{h}(x'_{0})}}(x'_{0})}|\mathcal{F}|\le\|(\varepsilon+\mathcal{h})^{\gamma}\|
_{L^{\infty}(Q_{s\sqrt{\varepsilon+\mathcal{h}(x'_{0})}}(x'_{0}))}\le C(\varepsilon+\mathcal{h}(x'_{0}))^{\gamma}.
\end{equation}
It follows from \eqref{partdelta/del esti}, \eqref{supDutdelta}, 
\eqref{supuesti}, and \eqref{supfest} that 
\begin{align*}
\sup_{Q_{t\sqrt{\varepsilon+\mathcal{h}(x'_{0})}}(x'_{0})}|\nabla_{x'}\mathcal{U}|\le\frac{3}{5}
\sup_{Q_{s\sqrt{\varepsilon+\mathcal{h}(x'_{0})}}(x'_{0})}|\nabla_{x'}\mathcal{U}|+\frac{C}{s-t}
(\varepsilon+\mathcal{h}(x'_{0}))^{\gamma+\frac12}.
\end{align*}
By Lemma \ref{itLemma}, for any $0<t<s\le\frac{1}{2(C_{a}+C_{a+1}+\mathcal{C})}
$, there holds
\begin{align*}
 \sup_{Q_{t\sqrt{\varepsilon+\mathcal{h}(x'_{0})}}(x'_{0})}|\nabla_{x'}\mathcal{U}|\le 
\frac{C(\varepsilon+\mathcal{h}(x'_{0}))^{\gamma+\frac12}}{s-t},
\end{align*}
which implies 
$$|\nabla_{x'}\mathcal{U}(x'_{0})|\le C(\varepsilon+\mathcal{h}(x'_{0}))^{\gamma+\frac12}\quad\quad\text{for}\quad x'_{0}
\in B_{2R}(0').$$
For $l\ge2$, continuing to differentiate both side of equation \eqref{zyfc/} and 
repeating the previous process, we can obtain \eqref{highestofu}. We thus complete the proof.
\end{proof}

At the end of this section, we rewrite
$${\boldsymbol\psi}_{1}=\begin{pmatrix}
1 \\
0\\
0
\end{pmatrix},~
{\boldsymbol\psi}_{2}=\begin{pmatrix}
0\\
1\\
0
\end{pmatrix},~
{\boldsymbol\psi}_{3}=\begin{pmatrix}
0\\
0\\
1
\end{pmatrix},~
{\boldsymbol\psi}_{4}=\begin{pmatrix}
x_{2}\\
-x_{1}\\
0
\end{pmatrix},~
{\boldsymbol\psi}_{5}=\begin{pmatrix}
x_{3}\\
0\\
-x_{1}
\end{pmatrix},~
{\boldsymbol\psi}_{6}=\begin{pmatrix}
0\\
x_{3}\\
-x_{2}
\end{pmatrix}.$$ 
Let $({\bf u}_{1}^\alpha,p_1^\alpha)$, $\alpha=1,2,\dots,6$, be the solution to \eqref{u,peq1}. 
In the next three sections, we establish the higher derivative estimates of the solutions by applying 
Proposition \ref{keyprop1}. For 
each $({\bf u}_1^\alpha, \bar{p}_1^\alpha)$, we will construct a series of auxiliary function pairs $({\bf v}_\alpha^l,\bar{p}_l)$, $l=1,2,\dots,m+1$, such 
that the quantity \(\Big|\sum_{l=1}^{m+1} (\mu \Delta {\bf v}_{\alpha}^{l} - 
\nabla \bar{p}_{l})(x)\Big|\) satisfies a suitable upper bound in \(\Omega_{2R}
\), as required by Proposition \ref{keyprop2}. As mentioned before, we will use the modified Keller-type function $k(x)\in C^{m+1}(\mathbb{R}^3)$ to construct $
({\bf v}_\alpha^l,\bar{p}_l)$,
\begin{equation}\label{def_k}
k(x)=\frac{x_3-\frac12(h_1-h_2)(x')}{\delta(x')}\quad\mbox{in}~\Omega_{2R},
\end{equation}
which satisfy $
k(x)+\frac12=1$ on $\Gamma_{2R}^+$, $k(x)+\frac12=0$ on
$\Gamma_{2R}^-$, and $k(x)^2-\frac14=0$ on $\Gamma_{2R}^+
\cup\Gamma_{2R}^-.$
A direct calculation shows that
\begin{equation}\label{kdds}
\partial_{x_{i}}k(x)=-\frac{\partial_{x_i}(h_1-h_2)(x')}{2\delta(x')}-\frac{\partial_{x_i}(h_1+h_2)(x')}{\delta(x')}k(x),\quad \partial_{x_{3}}k(x)=\frac{1}{\delta(x')}\quad\mbox{in}~\Omega_{2R}.
\end{equation}

\section{Estimates for \texorpdfstring{$({\bf u}_{1}^{\alpha},p_{1}^{\alpha}),\,\, \alpha=1,2$}{}}\label{sec_estu11,12}
In this section, we will establish high-order derivative estimates for the solution to equation \eqref{u,peq1} with \(\boldsymbol{\psi}_{1}\) and \(\boldsymbol{\psi}_{2}\) on the boundary \(\partial D_{1}\) by utilizing Proposition \ref{keyprop1}. For general domains \(D_1\) and \(D_2\), we construct auxiliary functions to accurately represent the leading terms up to the seventh-order derivatives. However, for eighth- and higher-order derivatives, we are only able to obtain an upper bound, as Proposition \ref{U_esti} is not applicable when $\gamma\ge -1$.

\begin{prop}\label{u11-general}
Under the same assumption as in Theorem \ref{main thm1}, let ${\bf u}_{1}
^{\alpha}$ and $p_{1}^{\alpha}$ be the solution to \eqref{u,peq1} with $
\alpha=1,2$. Then for sufficiently small $0<\varepsilon<1/2$, we have
\begin{equation*}
\mbox{(i)}~|\nabla{\bf u}_{1}^{\alpha}(x)|\le C(\varepsilon+|x'|^2)^{-1},\quad 
|p_1^\alpha(x)-p_1^\alpha(z',0)|\leq C(\varepsilon+|x'|^2)^{-\frac32}\quad\text{in}~
\Omega_{R}
\end{equation*}
for some point $z=(z',0)\in \Omega_{R}$ with $|z'|=R/2$; and

(ii) for high-order derivatives, we have
\begin{equation*}
\begin{aligned}
&|\nabla^{m+1}{\bf u}_{1}^{\alpha}(x)|+|\nabla^{m}{p}_{1}^{\alpha}(x)| \le 
C(\varepsilon+|x'|^2)^{-\frac{m+3}{2}}\quad\text{for}~m\le6, \\
&|\nabla^{m+1}{\bf u}_{1}^{\alpha}(x)|+|\nabla^{m}p_{1}^{\alpha}(x)|\le C
(\varepsilon+|x'|^2)^{\frac32-m}~~\quad\text{for}~m\ge7
\end{aligned}
\quad\text{in}~\Omega_{R}.
\end{equation*}
\end{prop}
\begin{proof}[Proof of Proposition \ref{u11-general}]
We only prove the case when $\alpha=1$ for instance, since the other case is the same. We will construct a series of auxiliary function pairs \(({\bf v}
_{1}^{l}(x), \bar{p}_{l}(x))\), $l=1,2,3$ with appropriate pointwise upper bound estimates 
in $\Omega_{2R}$. Unlike the two-dimensional case \cite{DLTZ1}, these functions are not explicit. However, we prove that \(\nabla^{m+1}
\sum_{l=1}^{m+1}{\bf v}_1^l\) and \(\nabla^m\sum_{l=1}^{m+1}{\bar{p}}_{l}\) 
 capture the leading terms in \(\nabla^{m+1}{\bf u}_{1}^{1}\) and 
\(\nabla^{m}{p}_{1}^1\) for $m\leq 6$, respectively.

\subsection{Gradient estimates}
In order to apply Proposition \ref{keyprop2} 
to derive the estimates of $\nabla{\bf u}_{1}^{1},$ we first construct 
a divergence-free function ${\bf v}_{1}^{1}$ satisfying ${\bf v}_{1}^{1}
={\bf u}_{1}^{1}={\boldsymbol\psi}_{1}$ on $\Gamma^{+}_{2R}$, ${\bf v}_{1}
^{1}={\bf u}_{1}^{1}=0$ on $
\Gamma^{-}_{2R}$, and suitable $\bar{p}_1$ such that 
\begin{equation}\label{yjds}
 |\mu\Delta {\bf v}_{1}^{1}- \nabla \bar{p}_1|\leq C\delta(x')^{-1}.
\end{equation}
Then, their difference
\begin{equation*}
{\bf w}:={\bf u}_{1}^{1}-{\bf v}_{1}^{1}\quad\mbox{and}\quad 
q:=p_{1}^{1}-\bar{p}_{1}
\end{equation*}
satisfy the following boundary value problem in the narrow region
\begin{equation*}
\begin{cases}
-\mu\Delta{\bf w}+\nabla q={\bf f}_{1}^{1}\quad&\mathrm{in}\ \Omega_{2R},\\
\nabla\cdot {\bf w}=0\quad&\mathrm{in}\ \Omega_{2R},\\
{\bf w}=0\quad&\mathrm{on}\ \Gamma^+_{2R}\cup\Gamma^-_{2R},
\end{cases}
\end{equation*}
where
\begin{equation}\label{f11_def}
{\bf f}_{1}^{1}:=\mu\Delta{\bf v}_{1}^{1}-\nabla\bar{p}_{1}.
\end{equation}
By Proposition \ref{keyprop2}, we have
$|\nabla {\bf u}_1^1|\leq |\nabla{\bf v}_{1}^{1}|+C.$

{\bf Step I.1. Construction of ${\bf v}_{1}^{1}$ and $\bar{p}_1$ and their 
estimates.}
We use the Keller-type function \eqref{def_k} to construct an auxiliary 
function ${\bf v}_{1}^{1}$ in the following form:
\begin{align}\label{v11def}
{\bf v}_{1}^{1}=\boldsymbol\psi_{1}\Big(k(x)+\frac{1}{2}\Big)+\begin{pmatrix}
F(x)\\
G(x)\\
 H(x)
\end{pmatrix}\big(k(x)^2-\frac14\big)\quad\hbox{in}\ \Omega_{2R},
\end{align}
where $F(x)$, $G(x)$, and $H(x)$ are to be determined in the following way. 
By a direct calculation, we have
\begin{equation}\label{bds1}
\begin{split}
&\partial_{x_{1}}({\bf v}_{1}^{1})^{(1)}=\partial_{x_{1}} F(x)\Big(k(x)^2-
\frac{1}{4}\Big)+2F(x)k(x)\partial_{x_{1}} k(x)+\partial_{x_{1}} k(x),\\
&\partial_{x_{2}}({\bf v}_{1}^{1})^{(2)}=\partial_{x_{2}} G(x)\Big(k(x)^2-
\frac{1}{4}\Big)+2G(x)k(x)\partial_{x_{2}} k(x),\\
&\partial_{x_{3}}({\bf v}_{1}^{1})^{(3)}=\partial_{x_{3}} H(x)\Big(k(x)^2-
\frac{1}{4}\Big)+\frac{2k(x)}{\delta(x')}H(x).
\end{split}
\end{equation}
In order to achieve $\nabla\cdot {\bf v}_1^1=0$ in $\Omega_{2R}$, we set 
\begin{equation}\label{ffz1}
\frac14(\partial_{x_{1}}F+\partial_{x_{2}}G+\partial_{x_{3}}H)(x)=-
\frac{\partial_{x_1}(h_1-h_2)(x')}{2\delta(x')}.
\end{equation}
By \eqref{bds1} and \eqref{ffz1}, we have 
\begin{align*}
&-\frac{2\partial_{x_1}(h_1-h_2)(x')}{\delta(x')}\big(k(x)^2-\frac14\big)+
\big(2F(x)\partial_{x_{1}} k(x)+2G(x)\partial_{x_{2}} k(x)\big)k(x)+\partial_{x_{1}} k(x)\\
&\qquad+\frac{2k(x)}{\delta(x')}H(x)=0. 
\end{align*}
By using \eqref{kdds}, we derive
\begin{equation}\label{ffz2}
H(x)=\partial_{x_{1}}(h_1-h_2)(x')k(x)+\frac12\partial_{x_{1}}(h_1+h_2)(x')-
\delta(x')\big(F(x)\partial_{x_1}k(x)-G(x)\partial_{x_2}k(x)\big).
\end{equation}
Combining this with \eqref{ffz1} and using \eqref{kdds} again, we see that $F(x)
$ and $G(x)$ satisfy the following first-order partial differential equation:
\begin{equation*}
\begin{aligned}
\partial_{x_{1}}F(x)+\partial_{x_{2}}G(x)&+\frac{\partial_{x_{1}}\delta(x')}
{\delta(x')}F(x)+\frac{\partial_{x_{2}}\delta(x')}{\delta(x')}G(x)\\
&-\delta(x')\big(\partial_{x_1}k(x)\partial_{x_{3}}F(x)-\partial_{x_2}k(x)
\partial_{x_{3}}G(x)\big)=-\frac{3\partial_{x_{1}}(h_1-h_2)(x')}{\delta(x')}.
\end{aligned}
\end{equation*}
In particular, if we assume $F(x)=F(x')$ and $G(x)=G(x')$, then we have
\begin{equation}\label{2by1}
\partial_{x_{1}}F(x')+\partial_{x_{2}}G(x')+\frac{\partial_{x_{1}}\delta(x')}{\delta(x')}F(x')+\frac{\partial_{x_{2}}\delta(x')}{\delta(x')}G(x')=-\frac{3\partial_{x_{1}}(h_1-h_2)(x')}{\delta(x')}.
\end{equation}

In order to cancel the leading terms in $\mu\Delta{\bf v}_1^1$ and derive \eqref{yjds}, we first choose 
$\hat{p}_{1}(x)=\hat{p}_{1}(x')$ 
so that 
\begin{equation}\label{tj1}
\partial_{x_{1}}\hat{p}_{1}(x')=\mu\partial_{x_3x_3}({\bf v}_1^1)^{(1)}=\frac{2\mu F(x')}{\delta(x')^2},\quad
\partial_{x_{2}}\hat{p}_{1}(x')=\mu\partial_{x_3x_3}({\bf v}_1^1)^{(2)}=\frac{2\mu G(x')}{\delta(x')^2}.
\end{equation}
Substituting \eqref{tj1} into \eqref{2by1}, we derive
\begin{equation}\label{zyfc1}
\Delta_{x'} \hat{p}_{1}(x')+\frac{3}{\delta(x')}
\nabla_{x'} \delta(x')\nabla_{x'}\hat{p}_{1}(x')=-\frac{6\mu\partial_{x_{1}}(h_1-h_2)(x')}
{\delta(x')^3}:=\mathcal{F}_1(x').
\end{equation}
By \eqref{h1h2prop}, we have
\begin{equation*}
|\mathcal{F}_1(x')|\leq\,C\delta(x')^{-\frac52}.
\end{equation*}
Similarly,
$$|\nabla_{x'}^s\mathcal{F}_1(x')|\leq\,C\delta(x')^{-\frac{5+s}{2}},\quad 1\le s\le m.$$
By Proposition \ref{U_esti} with $\gamma=-\frac52$, there exists $\hat{p}_{1}$ which satisfies \eqref{zyfc1}, with the following estimates
\begin{equation}\label{esthatp1}
 |\hat{p}_{1}|\leq C\delta(x')^{-\frac32}, \quad|\nabla_{x'}^k\hat{p}_{1}|\le 
 C\delta(x')^{-\frac{3+k}{2}},\quad\partial_{x_3}\hat{p}_{1}=0,\quad 0\le k\le 
 m+1.
\end{equation}
Then, by \eqref{tj1}, we have
\begin{equation}\label{u11-FG-form}
F(x')=\frac{\delta(x')^2}{2\mu}\partial_{x_1}\hat{p}_{1}(x'),\quad 
G(x')=\frac{\delta(x')^2}{2\mu}\partial_{x_2}\hat{p}_{1}(x'),
\end{equation}
and by \eqref{h1h2prop} and \eqref{esthatp1}, we get for $0\le k\le m+1$,
\begin{equation}\label{fgdegj1}
|\nabla_{x'}^k F|,~|\nabla_{x'}^k G|\leq C\delta(x')^{-\frac k2},\quad 
\partial_{x_3}F=\partial_{x_3}G=0.
\end{equation}

For $H(x)$, by \eqref{h1h2prop}, \eqref{def_k}, \eqref{ffz2}, and \eqref{fgdegj1}, one can see that $H(x)$ has the following form
\begin{equation}\label{hdxsu11}
H(x)=H_1(x')x_3+H_2(x'),
\end{equation}
with the estimates for $k\ge0,$ 
\begin{equation}\label{ys1}
|\nabla_{x'}^kH_1|\leq C\delta(x')^{-\frac{1+k}{2}},\quad |\nabla_{x'}^kH_2|\leq C\delta(x')^{\frac{1-k}{2}}.
\end{equation}
Thus, we obtain
\begin{equation}\label{estv11h}
 |\nabla_{x'}^k \partial_{x_3}^s H|\leq C\delta(x')^{\frac12-s-\frac k2},\quad s\leq 1,\quad\partial_{x_3}^2 H=0.
\end{equation}
Hence, by \eqref{v11def}, \eqref{fgdegj1}, and \eqref{estv11h}, we have, for 
$k\ge1$ and $i=1,2$,
\begin{equation}\label{v11highesti}
\begin{split}
&|\nabla_{x'}^{k}\partial_{x_3}^{s}({\bf v}_{1}^{1})^{(i)}|\le 
C\delta(x')^{-\frac{k}{2}-s}\quad\text{for}~s\le2,
\quad\partial_{x_3}^{s}({\bf v}_{1}^{1})^{(i)}=0\quad\text{for}~s\ge3,\\
&|\nabla_{x'}^{k}\partial_{x_3}^{s}({\bf v}_{1}^{1})^{(3)}|\le 
C\delta(x')^{\frac{1}{2}-\frac{k}{2}-s}\quad\text{for}~s\le3,
\quad\partial_{x_3}^{s}({\bf v}_{1}^{1})^{(3)}=0\quad\text{for}~s\ge4.
\end{split}
\end{equation}
In addition, it follows from \eqref{v11def}, \eqref{tj1}, and \eqref{v11highesti} that the 
leading term in $\mu\Delta{\bf v}_{1}^{1}-\nabla\hat{p}_{1}$ is $
\mu\partial_{x_{3}x_3}({\bf v}_1^1)^{(3)}$, which is of order $\delta(x')^{-
\frac32}$. To further cancel out this leading term, we choose
\begin{equation}\label{defhatp_1}
\tilde{p}_1=\frac{3\mu H_1(x')}{\delta(x')^2}x_3^2+2\mu \frac{H_2(x')-
(h_{1}-h_{2})H_{1}(x')}{\delta(x')^2}x_3\quad\text{in}~\Omega_{2R},
\end{equation}
so that
\begin{equation}\label{xzp_1}
\partial_{x_{3}}\tilde{p}_1=\mu\partial_{x_{3}x_3}({\bf v}_1^1)^{(3)}.
\end{equation}
Moreover, by \eqref{ys1} and \eqref{defhatp_1}, we have, for $k\ge 0$,
\begin{equation}\label{ys2}
|\nabla_{x'}^k\partial_{x_{3}}^s \tilde{p}_{1}|\leq C\delta(x')^{-\frac{1+k}
{2}-s}, ~s\le 2, \quad \partial_{x_{3}}^3\tilde{p}_{1}=0.
\end{equation}
Let $\bar{p}_{1}=\hat{p}_{1}+\tilde{p}_{1}.$ Combining \eqref{esthatp1}, 
\eqref{v11highesti}, and \eqref{ys2}, we have, for $m\ge1$,
\begin{equation}\label{estinabv11barp}
\begin{split}
 &|\nabla{\bf v}_{1}^{1}|\le C\delta(x')^{-1},\quad|\nabla^{m+1}{\bf v}_{1}^{1}|\le C 
\delta(x')^{-\frac{m+3}{2}},\\& |\bar{p}_{1}|\le C\delta(x')^{-\frac32},\quad|
\nabla^m\bar{p}_{1}|\le C\delta(x')^{-\frac{m+3}{2}}.
\end{split}
\end{equation}

{\bf Step I.2. Estimates of ${\bf f}_{1}^{1}$.} By virtue of \eqref{f11_def}, 
\eqref{tj1}, and \eqref{xzp_1}, we have, for $i=1,2$,
\begin{equation*}
({\bf f}_{1}^{1})^{(i)}=\mu\Delta_{x'}({\bf v}_1^1)^{(i)}-\partial_{x_i}
\tilde{p}_1,\quad({\bf f}_{1}^{1})^{(3)}=\mu\Delta_{x'}({\bf v}_{1}^{1})^{(3)}
\quad \text{in}~\Omega_{2R}.
\end{equation*}
By \eqref{v11def}, \eqref{u11-FG-form}, and \eqref{hdxsu11}, we notice that $({\bf v}
_1^1)^{(i)}$ is a quadratic polynomial in $x_3$, while $({\bf v}_1^1)^{(3)}$ is 
a cubic polynomial in $x_3$. Therefore, by \eqref{v11highesti} and \eqref{ys2}, 
we can rewrite ${\bf f}_{1}^{1}$ as polynomials in $x_{3}$:
\begin{align}\label{f11form}
\begin{split}
&({\bf f}_{1}^{1})^{(b)}=S_{21}^b(x')x_{3}^{2}+S_{11}^b(x')x_{3}
+S_{01}^{b}(x'), \quad b=1,2,\\
&({\bf f}_{1}^{1})^{(3)}=G_{31}(x')x_3^{3}+G_{21}(x')x_3^{2}+G_{11}
(x')x_3+G_{01}(x'),
\end{split}
\end{align}
where for any $s\ge 0$, $k=0,1,2$ and $j=0,1,2,3$, 
\begin{equation}\label{SGesti}
|\nabla_{x'}^sS_{k1}^{b}|\le C\delta(x')^{-k-\frac{s+2}{2}},\quad|\nabla_{x'}
^sG_{j1}|\le C\delta(x')^{-j-\frac{s+1}{2}}.
\end{equation}
Then, it is immediate from \eqref{f11form} and \eqref{SGesti} that
\begin{equation}\label{1-f11-1,2,3esti}
|({\bf f}_{1}^{1})^{(1)}|,|({\bf f}_{1}^{1})^{(2)}|\le C\delta(x')^{-1}, \quad|
({\bf f}_{1}^{1})^{(3)}|\le C\delta(x')^{-1/2}.
\end{equation}
Hence, 
\begin{equation*}
|{\bf f}_{1}^{1}|\le C (|({\bf f}_{1}^{1})^{(1)}|+|({\bf f}_{1}^{1})^{(2)}|+|
({\bf f}_{1}^{1})^{(3)}|)\le C\delta(x')^{-1}\quad\text{in}~\Omega_{2R}.
\end{equation*}
By Proposition \ref{keyprop2}, we obtain $
\|\nabla{\bf u}_{1}^{1}-\nabla{\bf v}_{1}^{1}\|_{L^{\infty}(\Omega_{\delta(x')/2}(x'))}\le C~\mbox{for}~ x\in\Omega_{R}.$
Thus, combining this with \eqref{estinabv11barp} yields
$$|\nabla{\bf u}_{1}^{1}(x)|\le C|\nabla {\bf v}_{1}^{1}(x)|+C \le 
C\delta(x')^{-1}\quad \text{in}~\Omega_{R}.$$
Thus, Proposition \ref{u11-general} holds true for $m=0$.

\subsection{Second-order derivatives estimates.}
 In order to derive second-order derivative estimates for the solutions of \eqref{u,peq1}, we construct new auxiliary functions to further reduce the upper bounds of the non-homogeneous term.

{\bf Step II.1. Construction of ${\bf v}_{1}^{2},\, \bar{p}_{2}$ and their estimates.} By \eqref{1-f11-1,2,3esti}, one can derive that the leading terms in ${\bf f}_1^1$ are $({\bf f}_{1}^{1})^{(1)}$ and $
({\bf f}_{1}^{1})^{(2)}$. To cancel these leading terms, we construct a divergence-free function 
$${\bf v}_{1}^{2}=\Big(({\bf v}_{1}^{2})^{(1)}(x),({\bf v}_{1}^{2})^{(2)}(x),
({\bf v}_{1}^{2})^{(3)}(x)\Big)^{T}$$ 
satisfying 
\begin{equation}\label{u11-v12-1,2}
({\bf v}_{1}^{2})^{(b)}(x)=\Big(\sum_{i=0}^{2}F_{b2}^{i}(x')x_{3}^{i}+\tilde{F}
_{b2}(x') \Big)\big(k(x)^2-\frac14\big),\,\, b=1,2,
\end{equation} 
 and 
\begin{equation*}
 ({\bf v}_{1}^{2})^{(3)}(x)=\Big(\sum_{i=0}^{3}F_{32}^{i}(x')x_{3}^{i}+\tilde{F}
 _{32}(x) \Big)\big(k(x)^2-\frac14\big).
 \end{equation*} 
 
 We emphasize that \(F_{b2}^{i}(x')\), \(i=0,1,2\), are chosen to cancel out the leading terms in ${\bf f}_1^{1}$ in the following way
\begin{equation}\label{cancelf12}
\mu\partial_{x_3x_3}\Big(\sum_{i=0}^{2}F_{b2}^{i}(x')x_3^{i}\Big(k(x)^2-
\frac14\Big)\Big)=-({\bf f}_{1}^{1})^{(b)},\quad b=1,2.
\end{equation}
By comparing the coefficients of each term of $x_3$, we obtain
\begin{equation*}
\begin{aligned}
&F_{b2}^{2}(x')=-\frac{\delta(x')^2}{12\mu}S_{21}^{b}(x'),\,\,F_{b2}^{1}(x')=-
\frac{\delta(x')^2}{6\mu}S_{11}^{b}(x')+(h_1-h_{2})(x')F_{b2}^{2}(x'),\\
&F_{b2}^{0}(x')=-\frac{\delta(x')^2}{2\mu}S_{01}^{b}(x')+(h_1-h_2)(x')F_{b2}^{1}
(x')+\frac14(\varepsilon+2h_{1}(x'))(\varepsilon+2h_{2}(x'))F_{b2}^{2}(x').
\end{aligned}
\end{equation*}
To ensure $\nabla\cdot {\bf v}_1^2=0$, we first take \(F_{32}^{i}, i=0,1,2,3\), so that
\begin{equation*}
\nabla\cdot
\begin{pmatrix}
\Big(\sum_{i=0}^{2}F_{12}^{i}(x')x_{3}^{i}\Big)\big(k(x)^2-\frac14\big)\\
\Big(\sum_{i=0}^{2}F_{22}^{i}(x')x_{3}^{i}\Big)\big(k(x)^2-\frac14\big)\\
\Big(\sum_{i=0}^{3}F_{32}^{i}(x')x_{3}^{i}\Big)\big(k(x)^2-\frac14\big)
\end{pmatrix}:=R(x')
\end{equation*}
is the $0$-th order term in the polynomial of $x_3$. By comparing the 
coefficients of each term of $x_3$, we have 
\begin{align}
F_{32}^{3}(x')=&\,-\frac{\delta(x')^{2}}{5}\sum_{i=1}^{2}\partial_{x_i}
\Big(\frac{F_{i2}^{2}(x')}{\delta(x')^{2}}\Big), \nonumber \\ 
F_{32}^{2}(x')=&\,-\frac{\delta(x')^2}{4}\sum_{i=1}^{2}\partial_{x_i}
\Big(\frac{F_{i2}^{1}(x')-(h_1-h_2)(x')F_{i2}^{2}(x')}{\delta(x')^2}\Big)+ (h_1-
h_2)(x')F_{32}^{3}(x'), \nonumber \\
F_{32}^{1}(x')=&\,-\frac{\delta(x')^2}{3}\sum_{i=1}^{2}\partial_{x_i}\Big( 
\frac{F_{i2}^{0}(x')}{\delta(x')^2}-\frac{(h_1-h_2)(x')F_{i2}^{1}(x')}
{\delta(x')^2}- \frac{(\varepsilon+2h_1)(\varepsilon+2h_2)}{4\delta(x')^2}F_{i2}
^{2}(x')\Big) \label{F12i22iform} \\ 
&\,+\frac14(\varepsilon+2h_1(x'))(\varepsilon+2h_2(x'))F_{32}^{3}(x')+
(h_1-h_2)(x')F_{32}^{2}(x'),\nonumber \\
F_{32}^{0}(x')=&\,\frac{\delta(x')^2}{2}\sum_{i=1}^{2}\partial_{x_i}
\Big(\frac{(h_1-h_2)(x')F_{i2}^{0}(x')}{\delta(x')^2}+\frac{(\varepsilon+2h_1)
(\varepsilon+2h_2)}{4\delta(x')^2}F_{i2}^{1}(x')\Big) \nonumber \\
&\,+\frac14(\varepsilon+2h_1(x'))(\varepsilon+2h_2(x'))F_{32}^{2}(x')+(h_1-
h_2)F_{32}^{1}(x'), \nonumber
\end{align}
and
\begin{equation}\label{1Rform} 
\begin{aligned}
R(x')=&\, -\sum_{i=1}^{2}\partial_{x_i}\Big(\frac{(\varepsilon+2h_1(x'))
(\varepsilon+2h_2(x'))}{4\delta(x')^2}F_{i2}^{0}(x') \Big) \\
&\,-\frac{(\varepsilon+2h_1(x'))(\varepsilon+2h_2(x'))}{4\delta(x')^2}F_{32}^{1}
(x')-\frac{h_1-h_2}{\delta(x')^2}F_{32}^{0}(x').
\end{aligned} 
\end{equation}
By \eqref{SGesti}, \eqref{F12i22iform}, and \eqref{1Rform}, we 
have for $k\ge 0$ and $i=0,1,2$, 
\begin{equation}\label{F1222esti}
\begin{aligned}
&|\nabla_{x'}^kF_{b2}^{i}|\le C\delta(x')^{1-i-\frac{k}{2}},b=1,2,\,\,|
\nabla_{x'}^kF_{32}^{i}|\le C\delta(x')^{\frac32-i-\frac{k}{2}},\\
&|\nabla_{x'}^kF_{32}^{3}|\le C\delta(x')^{-\frac32-\frac{k}{2}},\quad
|\nabla_{x'}^{k}R|\le C \delta(x')^{\frac12-\frac{k}{2}}.
\end{aligned}
\end{equation}
Finally, we choose 
$\tilde{F}_{12}(x')$, $\tilde{F}_{22}(x')$, and $\tilde{F}_{32}(x)$ such that
\begin{equation}\label{tilFb2help}
\nabla \cdot\begin{pmatrix}
\tilde{F}_{12}(x')\big(k(x)^2-\frac14\big) \\
\tilde{F}_{22}(x')\big(k(x)^2-\frac14\big)\\
\tilde{F}_{32}(x)\big(k(x)^2-\frac14\big)
\end{pmatrix}=-R(x'), 
\end{equation}
so that $\nabla\cdot{\bf v}_1^2=0$. By a direct calculation, \eqref{tilFb2help} is equivalent to
\begin{align}\label{tilFb2eqs}
\big(\partial_{x_{1}}\tilde{F}_{12}(x')+\partial_{x_{2}}\tilde{F}_{22}(x')&+
\partial_{x_{3}}\tilde{F}_{32}(x)\big)\big(k(x)^2-\frac14\big)+
\frac{2k(x)}{\delta(x')}\tilde{F}_{32}(x)\nonumber\\
&+2\big(\tilde{F}_{12}(x')\partial_{x_{1}}k(x)+\tilde{F}_{22}(x')
\partial_{x_{2}}k(x)\big)k(x)=-R(x').
\end{align}
To determine $\tilde{F}_{12}(x')$ and $\tilde{F}_{22}(x)$, we set 
\begin{equation}\label{tilFbeq1}
\partial_{x_{1}}\tilde{F}_{12}(x')+\partial_{x_{2}}\tilde{F}_{22}(x')+
\partial_{x_{3}}\tilde{F}_{32}(x)=4R(x').
\end{equation}
It follows from \eqref{tilFb2eqs} and \eqref{tilFbeq1} that
\begin{equation}\label{tilF32eqs}
\tilde{F}_{32}(x)=-2\delta(x')R(x')k(x)-\delta(x')\partial_{x_1}k(x)
\tilde{F}_{12}(x')-\delta(x')\partial_{x_2}k(x)\tilde{F}_{22}(x').
\end{equation}
By \eqref{kdds}, \eqref{tilFbeq1}, and \eqref{tilF32eqs}, we have
\begin{equation}\label{TilFb2andR}
\partial_{x_{1}}\tilde{F}_{12}(x')+\partial_{x_{2}}\tilde{F}_{22}(x')+
\frac{\partial_{x_{1}}\delta(x')}{\delta(x')}\tilde{F}_{12}(x')+
\frac{\partial_{x_{2}}\delta(x')}{\delta(x')}\tilde{F}_{22}(x')=6R(x').
\end{equation}
Similar to the argument that led to \eqref{tj1}, we select $\hat{p}_{2}(x')$ 
such that
\begin{equation}\label{hatp2eq}
\partial_{x_1}\hat{p}_{2}=\frac{2\mu\tilde{F}_{12}(x')}{\delta(x')^2},\quad
\partial_{x_2}\hat{p}_{2}=\frac{2\mu\tilde{F}_{22}(x')}{\delta(x')^2}.
\end{equation}
Substituting \eqref{hatp2eq} into \eqref{TilFb2andR} yields
\begin{equation}\label{Deltahatp2eq}
\Delta_{x'}\hat{p}_{2}+\frac{3}{\delta(x')}\nabla_{x'}\delta(x')
\nabla_{x'}\hat{p}_{2}=\frac{12\mu R(x')}{\delta(x')^2}:=\mathcal{F}_2(x').
\end{equation}
By \eqref{F1222esti}, we have
\begin{equation*}
|\nabla_{x'}^s\mathcal{F}_2(x')|\leq\,C\delta(x')^{-\frac{3+s}{2}},\quad 0\le s\le m.
\end{equation*}
By using Proposition \ref{U_esti} with $\gamma=-\frac32$, there exists $\hat{p}_{2}$ that satisfies \eqref{Deltahatp2eq} and for any $0\le l\le m+1,$ 
\begin{equation}\label{hatp12esti} 
|\nabla_{x'}^{l}\hat{p}_{2}|\le C\delta(x')^{-\frac12-\frac{l}{2}}.
\end{equation}
Combining \eqref{tilF32eqs}, \eqref{hatp2eq} and \eqref{hatp12esti}, we have,
for any $k\ge0,$ 
\begin{equation}\label{tilF122232esti}
|\nabla_{x'}^{k}\tilde{F}_{b2}|\le C\delta(x')^{1-\frac{k}{2}},b=1,2,\quad|
\nabla_{x'}^{k}\tilde{F}_{32}|\le C\delta(x')^{\frac32-\frac{k}{2}}.
\end{equation}
By \eqref{f11form}, \eqref{cancelf12}, and \eqref{hatp2eq}, a calculation gives 
\begin{align}\label{f11+Delv12-nabhatp2}
 {\bf f}_{1}^{1}+\mu\Delta{\bf v}_{1}^{2}-\nabla\hat{p}_{2}=\begin{pmatrix}
 \mu \Delta_{x'}({\bf v}_{1}^{2})^{(1)}\\
\mu\Delta_{x'}({\bf v}_{1}^{2})^{(2)}\\
({\bf f}_{1}^{1})^{(3)}+\mu\Delta({\bf v}_{1}^{2})^{(3)}
 \end{pmatrix}
=\begin{pmatrix}
\sum_{i=0}^{4}H_{i2}^{1}(x')x_3^{i}\\
\sum_{i=0}^{4}H_{i2}^{2}(x')x_3^{i}\\
\sum_{i=0}^{5}G_{i2}(x')x_3^{i}+\sum_{i=0}^{3}\tilde{S}_{i2}(x')x_3^{i}
\end{pmatrix},
\end{align}
where 
$\sum_{i=0}^{5}G_{i2}(x')x_3^{i}=\mu\Delta_{x'}({\bf v}_{1}^{2})^{(3)}
~\mbox{and}~ \sum_{i=0}^{3}\tilde{S}_{i2}(x')x_3^{i}= ({\bf f}_{1}
^{1})^{(3)}+\mu\partial_{x_3x_3}({\bf v}_{1}^{2})^{(3)}.$
By \eqref{h1h2prop}, \eqref{F1222esti}, \eqref{hatp12esti} and 
\eqref{tilF122232esti}, we obtain, for $k\ge0,$
\begin{equation}\label{HGSi2esti}
|\nabla_{x'}^{k}H_{i2}^{b}|\le C \delta(x')^{-i-\frac{k}{2}},b=1,2,\,\,|
\nabla_{x'}^{k}G_{i2}|\le\delta(x')^{\frac12-i-\frac{k}{2}},\,\,\,|\nabla_{x'}
^{k}\tilde{S}_{i2}|\le C\delta(x')^{-\frac12-i-\frac{k}{2}}.
\end{equation}
Then from \eqref{f11+Delv12-nabhatp2} and \eqref{HGSi2esti}, we have
\begin{align*}
&|\Delta_{x'}({\bf v}_{1}^{2})^{(b)}|\leq C,~b=1,2,\quad|\Delta_{x'}({\bf v}_{1}
^{2})^{(3)}|\leq C\delta(x')^{\frac12}, \\
& |({\bf f}_{1}^{1})^{(3)}+\mu\partial_{x_3x_3}({\bf v}_{1}^{2})^{(3)}|\leq 
C\delta(x')^{-\frac12}.
\end{align*}
Notice that the leading term in ${\bf f}_{1}^{1}+\mu\Delta{\bf v}_{1}
^{2}-\nabla\hat{p}_{2}$ is $\sum_{i=0}^{3}\tilde{S}_{i2}(x')x_3^{i}$.
Thus, to further reduce the upper bound of the inhomogeneous term, we choose 
\begin{equation}\label{tilp12}
\tilde{p}_{2}=\sum_{i=0}^{3}\frac{\tilde{S}_{i2}(x')}{i+1}x_3^{i+1}
\quad\text{in}~\Omega_{2R},
\end{equation}
such that
$\partial_{x_{3}}\tilde{p}_{2}=\sum_{i=0}^{3}\tilde{S}_{i2}(x')x_3^{i}.$

{\bf Step II.2. Estimates of ${\bf f}_{1}^{2}$.}
Let $\bar{p}_{2}:=\hat{p}_{2}+\tilde{p}_{2}$ and ${\bf f}_{1}^{2}:={\bf f}_{1}^{1}+\mu\Delta{\bf v}_{1}^{2}-\nabla\bar{p}
_{2}.$ By \eqref{f11+Delv12-nabhatp2} and \eqref{tilp12}, we can write ${\bf f}
_{1}^{2}$ as polynomials in $x_3$:
\begin{equation*}
\begin{split}
({\bf f}_{1}^{2})^{(b)}=&\,\sum_{i=1}^{4}\Big(H_{i2}^{b}-\frac1i\partial_{x_{b}}
\tilde{S}_{(i-1)2}\Big)(x')x_3^{i}+H_{02}^{b}(x'):=\sum_{i=0}^{4}S_{i2}^{b}
(x')x_3^{i},~b=1,2,\\
({\bf f}_{1}^{2})^{(3)}=&\,\sum_{i=0}^{5}G_{i2}(x')x_3^{i},
\end{split}
\end{equation*} 
where $S_{i2}^{b}(x')=H_{i2}^{b}(x')-\frac1i\partial_{x_{b}}\tilde{S}_{(i-1)2}
(x')$, $1\le i\le4$, and $S_{02}^{b}(x')=H_{02}^{b}(x').$ 
In view of \eqref{HGSi2esti}, we have for $k\ge0,$
\begin{equation}\label{Si2bGi2besti}
|\nabla_{x'}^{k}S_{i2}^{b}|\le C\delta(x')^{-i-\frac{k}{2}},\quad|\nabla_{x'}
^{k}G_{i2}|\le C\delta(x')^{\frac12-i-\frac{k}{2}}.
\end{equation}
It is immediate that $|({\bf f}_{1}^{2})^{(1)}|+ |({\bf f}_{1}^{2})^{(2)}|\leq C 
$ and $|({\bf f}_{1}^{2})^{(3)}|\le C\delta(x')^{\frac12}.$ Thus, we have 
\begin{equation}\label{f12esti} 
|{\bf f}_{1}^{2}|\le |({\bf f}_{1}^{2})^{(1)}|+ |({\bf f}_{1}^{2})^{(2)}|+|({\bf 
f}_{1}^{2})^{(3)}| \le C.
\end{equation}

{\bf Step II.3. Estimates of $p_{1}^{1}(x),\nabla p_{1}^{1}(x)$, and $
\nabla^2{\bf u}_{1}^{1}(x)$.} By using \eqref{u11-v12-1,2}, \eqref{F1222esti}, \eqref{hatp12esti},
\eqref{tilF122232esti}, and \eqref{tilp12}, a calculation gives for 
$k\ge0$ and $b=1,2,$
\begin{align}\label{v12hatp2tilp2est}
\begin{split}
&|\nabla_{x'}^{k}\partial_{x_3}^{s}({\bf v}_{1}^{2})^{(b)}|
\le\delta(x')^{1-s-\frac{k}{2}}\quad\text{for}~ s\le 4,\quad
\partial_{x_3}^{s}({\bf v}_{1}^{2})^{(b)}=0\quad\text{for}~s\ge5, \\
&|\nabla_{x'}^{k}\partial_{x_3}^{s}({\bf v}_{1}^{2})^{(3)}|
\le\delta(x')^{\frac32-s-\frac{k}{2}}\quad\text{for}~ s\le 5,\quad
\partial_{x_3}^{s}({\bf v}_{1}^{2})^{(3)}=0\quad\text{for}~s\ge 
6, \\
&|\nabla_{x'}^{k}\partial_{x_3}^{s}\tilde{p}_{2}|\le 
C\delta(x')^{\frac12-s-\frac{k}{2}}\quad\text{for}~0\leq s\le4,\\
&\partial_{x_3}^{s}\tilde{p}_{2}=0\quad\text{for}~s\ge 5,\quad \quad\quad |
\nabla_{x'}^{k}\hat{p}_{2}|\le C\delta(x')^{-\frac12-\frac{k}{2}}.
\end{split}
\end{align}
 By \eqref{f12esti} and Proposition 
\ref{keyprop2}, it holds that
\begin{equation}\label{u11-v2est}
\|\nabla^{2}({\bf u}_{1}^{1}-({\bf v}_{1}^{1}+{\bf v}_{1}^{2}))\|_{L^{\infty}
(\Omega_{\delta(x')/2}(x'))}+\|\nabla(p_{1}^{1}-\bar{p}_{1}-\bar{p}_{2})\|
_{L^{\infty}(\Omega_{\delta(x')/2}(x'))} \le C.
\end{equation}
By \eqref{v12hatp2tilp2est}, we obtain
\begin{equation}\label{v12est}
\begin{split}
&|\nabla^{m+1}{\bf v}_{1}^{2}|\le C\delta(x')^{-\frac{m+3}{2}}
\quad\text{for}~m\ge3 ,\quad
|\nabla^{m+1}{\bf v}_{1}^{2}|\le C\delta(x')^{-m}
\quad\text{for}~0\le m\le 2,\\
&|\nabla^{m}\bar{p}_{2}|\le C \delta(x')^{-\frac{m+3}{2}}\quad\text{for}
~ m\ge5,\quad\quad\, |\nabla\bar{p}_{2}|\leq C\delta(x')^{-1},\\ 
&|\nabla^{m}\bar{p}_{2}|\le C\delta(x_1)^{\frac12-m}\quad\text{for}~2\le 
m \le 4 ,~~~\,~ |\bar{p}_{2}|\leq C\delta(x_1)^{-\frac12}.
\end{split}
\end{equation}
It follows from \eqref{estinabv11barp} and \eqref{v12est} that for 
$x\in\Omega_{R}$,
\begin{equation}\label{gjv2barp2}
|\nabla^{2}({\bf v}_{1}^{1}(x)+ {\bf v}_{1}^{2}(x))|+|\nabla(\bar{p}_{1}(x)+
\bar{p}_{2}(x))|\leq C\delta(x')^{-2}.
\end{equation}
By \eqref{u11-v2est} and \eqref{gjv2barp2}, we have for $x\in\Omega_{R},$
\begin{equation*}
|\nabla^{2}{\bf u}_{1}^{1}(x)|+|\nabla{p}_{1}^{1}(x)|\le C \delta(x')^{-2} 
\quad\text{in}~\Omega_{R}.
\end{equation*}
For the estimates of $p_1^1(x)$, by the mean value theorem, 
\eqref{estinabv11barp}, \eqref{u11-v2est}, and \eqref{v12est}, we have
\begin{align*}
|p_{1}^{1}(x)-p_{1}^{1}(z',0)|\leq&\,| p_{1}^{1}-\bar{p}_1-\bar{p}_2-
(p_{1}^{1}-\bar{p}_1-\bar{p}_2)(z',0)|+|\bar{p}_1+\bar{p}_2|+C\\
\leq&\, C\|\nabla (p_{1}^{1}-\bar{p}_1-\bar{p}_2)\|_{L^\infty(\Omega_{R})}+|
\bar{p}_1+\bar{p}_2|+C
\leq\, C \delta(x')^{-3/2}
\end{align*}
for a fixed point $(z',0)\in\Omega_{R}$ with $|z'|=R/2$.

\subsection{High-order derivatives estimates.}
 Next, we inductively derive 
the estimates of $\nabla^{m+1}{\bf u}_{1}^{1}$ and $\nabla^{m}p_{1}^{1}$ 
based on Propositions \ref{keyprop1} and \ref{keyprop2}.
Denote
\begin{equation*}
{\bf f}_{1}^{j}:={\bf f}_{1}^{j-1}(x)+\mu\Delta{\bf v}_{1}^{j}(x)-\nabla\bar{p}
_{j}(x)=\sum_{l=1}^{j}(\mu\Delta{\bf v}_{1}^{l}-\nabla\bar{p}_{l})
(x)\quad\text{for}~ 2\le j\le m+1.
\end{equation*}
To use Proposition \ref{keyprop2} to derive higher derivative estimates, we successively construct auxiliary functions to make the non-homogeneous term as small as possible. To cancel out the singularity of $({\bf f}_{1}^{2})^{(b)},b=1,2,$ we choose ${\bf v}_{1}^{3}(x)$ satisfying
\begin{equation}\label{v13-12}
\begin{aligned}
 ({\bf v}_{1}^{3}(x))^{(b)}=\Big(\sum_{i=0}^{4}F_{b3}^{i}(x')x_{3}
 ^{i} \Big)\big(k(x)^2-\frac14\big),~b=1,2,
\end{aligned}
\end{equation}
where 
\begin{equation}\label{F1323form}
 F_{b3}^{i}(x')=-\frac{\delta(x')^2S_{i2}^{b}(x')}{\mu(i+1)(i+2)}+
 (h_1-h_2)(x')F_{b3}^{i+1}(x')+\frac14(\varepsilon+2h_1)
 (\varepsilon+2h_2)F_{b3}^{i+2}(x'),
\end{equation}
such that 
\begin{equation}\label{cancel1f12}
\mu\partial_{x_3x_3}({\bf v}_{1}^{3})^{(b)} =-({\bf f}_{1}^{2})^{(b)},\quad 
b=1,2.
\end{equation}
Here we use the convention that $F_{b3}^{i}(x')\equiv0$ if $i\notin\{0,\dots,4\}
$ and $F_{33}^{i}=0$ if $i\notin\{0,\dots,5\}$. We note that once $F_{b3}^i$, 
$i=0,1,\dots,4$ is determined, in general, there does not exist $F_{33}^i$, 
$i=0,1,\dots,5$ such that $\nabla\cdot{\bf v}_1^3=0$.
By comparing the coefficients of the polynomial of $x_3$, we take 
\begin{equation}\label{v13-3form}
 ({\bf v}_{1}^{3})^{(3)}=\sum_{i=0}^{5}F_{33}^{i}(x')x_{3} ^{i}\big(k(x)^2-\frac14\big),
\end{equation}
where 
\begin{equation}
\begin{aligned}
&F_{33}^{i}(x')=(h_1-h_2)(x')F_{33}^{i+1}(x')+\frac14(\varepsilon+2h_1(x'))
(\varepsilon+2h_2(x'))F_{33}^{i+2}(x') \\
&\,-\frac{\delta(x')^2}{i+2}\sum_{j=1}^{2}\partial_{x_j}\Big(\delta(x')^{-2}
\Big(F_{13}^{i-1}(x')-(h_1-h_2)F_{13}^{i}(x')-\frac14(\varepsilon+2h_1)
(\varepsilon+2h_2)F_{i3}^{i+1}(x')\Big)\Big) \label{F33form},
\end{aligned}
\end{equation}
such that 
$\nabla\cdot{\bf v}_{1}^{3}=R(x'),$
where 
\begin{equation}\label{3Rx'}
\begin{aligned}
R(x')=&\,-\sum_{i=1}^{2}\partial_{x_i}\Big(\frac{(\varepsilon+2h_1(x'))
(\varepsilon+2h_2(x'))}{4\delta(x')^2}F_{i3}^{0}(x')\Big)\\
&\,-\frac{(\varepsilon+2h_1(x'))(\varepsilon+2h_2(x'))}{4\delta(x')^2}F_{33}
^{1}(x')-\frac{h_1(x')-h_2(x')}{\delta(x')^2}F_{33}^{0}(x').
\end{aligned}
\end{equation}
By \eqref{h1h2prop}, \eqref{Si2bGi2besti}, \eqref{F1323form}, \eqref{F33form}, 
and \eqref{3Rx'}, 
we have for $k\ge 0$, $i=0,1,2,3,4$, and $b=1,2$,
\begin{equation}\label{F1323esti}
\begin{aligned}
&|\nabla_{x'}^kF_{b3}^{i}|\le C\delta(x')^{2-i-\frac{k}{2}},\,\,|
\nabla_{x'}^kF_{32}^{i}|\le C\delta(x')^{\frac52-i-\frac{k}{2}},\,\,|\nabla_{x'}^kF_{32}^{5}|\le C\delta(x')^{-\frac52-\frac{k}{2}},
\end{aligned}
\end{equation}
and 
\begin{equation}\label{rdgj24}
 |\nabla_{x'}^{k}R|\le C \delta(x')^{\frac32-\frac{k}{2}}.
\end{equation}
From \eqref{v13-12}, \eqref{cancel1f12}, and \eqref{v13-3form}, a calculation gives 
\begin{equation}\label{f12+Delv13}
\begin{aligned}
&{\bf f}_{1}^{2}+\mu\Delta{\bf v}_{1}^{3}=\begin{pmatrix}
\mu \Delta_{x'}({\bf v}_{1}^{3})^{(1)}\\
\mu\Delta_{x'}({\bf v}_{1}^{3})^{(2)}\\
({\bf f}_{1}^{2})^{(3)}+\mu\Delta({\bf v}_{1}^{3})^{(3)}
\end{pmatrix}
=\begin{pmatrix}
\sum_{i=0}^{6}H_{i3}^{1}(x')x_3^{i}\\
\sum_{i=0}^{6}H_{i3}^{2}(x')x_3^{i}\\
\sum_{i=0}^{7}G_{i3}(x')x_3^{i}+\sum_{i=0}^{5}\tilde{S}_{i3}(x')x_3^{i}
\end{pmatrix},
\end{aligned}
\end{equation}
where 
$$\sum_{i=0}^{7}G_{i3}(x')x_3^{i}=\mu\Delta_{x'}({\bf v}_{1}^{3})^{(3)}
\quad\mbox{and}\quad\sum_{i=0}^{5}\tilde{S}_{i2}(x')x_3^{i}= ({\bf f}_{1}
^{2})^{(3)}+\mu\partial_{x_3x_3}({\bf v}_{1}^{3})^{(3)}.$$
By \eqref{F1323esti} and \eqref{rdgj24}, we obtain for $k\ge0,$
\begin{equation}\label{HGSi3esti}
|\nabla_{x'}^{k}H_{i3}^{b}|\le C\delta(x')^{1-i-\frac{k}{2}},~b=1,2,\,\,|
\nabla_{x'}^{k}G_{i2}|\le\delta(x')^{\frac32-i-\frac{k}{2}},\,\,\,|
\nabla_{x'}^{k}\tilde{S}_{i3}|\le C\delta(x')^{\frac12-i-\frac{k}{2}}.
\end{equation}
By \eqref{f12+Delv13} and \eqref{HGSi3esti}, we derive, for $b=1,2$,
\begin{equation*}
|\Delta_{x'}({\bf v}_{1}^{3})^{(b)}|\leq C\delta(x'),\,\, 
|\Delta_{x'}({\bf v}_{1}^{3})^{(3)}|\leq C\delta(x')^{\frac32}, \,\,|({\bf f}_{1}^{2})^{(3)}+\mu\partial_{x_3x_3}({\bf v}_{1}^{3})^{(3)}|
\leq C\delta(x')^{\frac12}.
\end{equation*}
Thus the leading term in ${\bf f}_{1}^{2}+\mu\Delta{\bf v}_{1}^{3}$ is $
\sum_{i=0}^{5}\tilde{S}_{i3}(x')x_3^{i}.$ To cancel this leading term, we take 
\begin{equation}\label{barp3form}
\bar{p}_{3}=\sum_{i=0}^{5}\frac{\tilde{S}_{i3}(x')}{i+1}x_3^{i+1}
\quad\text{in}~\Omega_{2R},
\end{equation}
such that
$\partial_{x_{3}}\bar{p}_{3}=\sum_{i=0}^{5}\tilde{S}_{i3}(x')x_3^{i}.$

From \eqref{f12+Delv13} and \eqref{barp3form}, we can write ${\bf f}_{1}^{3}$ as polynomials in $x_3$:
\begin{equation*}
\begin{split}
({\bf f}_{1}^{3})^{(b)}=&\,\sum_{i=1}^{6}\Big(H_{i3}-
\frac1i\partial_{x_{b}}
\tilde{S}_{(i-1)3}\Big)(x')x_3^{i} +H_{03}(x'):=\sum_{i=0}^{6}S_{i3}^{b}
(x')x_3^{i},\,\, b=1,2,
\\({\bf f}_{1}^{3})^{(3)}=&\,\sum_{i=0}^{7}G_{i3}(x')x_3^{i},
\end{split}
\end{equation*} 
where $S_{i3}^{b}(x')=H_{i3}(x')-\frac1i\partial_{x_{b}}\tilde{S}_{(i-1)3}(x')$, $1\le i\le6$, and $S_{03}^{b}(x')=H_{03}(x').$
In view of \eqref{HGSi3esti}, we have for $k\ge0,$
\begin{equation*}
|\nabla_{x'}^{k}S_{i3}^{b}|\le C\delta(x')^{1-i-\frac{k}{2}},\quad|
\nabla_{x'}^{k}G_{i3}|\le C\delta(x')^{\frac32-i-\frac{k}{2}}.
\end{equation*}
Thus, we have
\begin{equation*}
|({\bf f}_{1}^{3})^{(b)}|\leq C\delta(x'),~b=1,2,\quad |({\bf f}_{1}^{3})^{(3)}|
\leq C\delta(x')^{\frac32},
\end{equation*}
and
\begin{equation}\label{f13esti}
|{\bf f}_{1}^{3}|\le |({\bf f}_{1}^{3})^{(1)}|+|({\bf f}_{1}^{3})^{(2)}|+|({\bf 
f}_{1}^{3})^{(3)}|\le C\delta(x').
\end{equation}
By using \eqref{v13-12}, \eqref{v13-3form}, \eqref{F1323esti}, \eqref{HGSi3esti}, and \eqref{barp3form}, a calculation gives for 
$k\ge0$ and $b=1,2,$
\begin{align}\label{v13barp3est}
\begin{split}
&|\nabla_{x'}^{k}\partial_{x_3}^{s}({\bf v}_{1}^{3})^{(b)}|
\le\delta(x')^{2-s-\frac{k}{2}}\quad\text{for}~ s\le 6,\quad
\partial_{x_3}^{s}({\bf v}_{1}^{3})^{(b)}=0\quad\text{for}~s\ge7, \\
&|\nabla_{x'}^{k}\partial_{x_3}^{s}({\bf v}_{1}^{3})^{(3)}|
\le\delta(x')^{\frac52-s-\frac{k}{2}}\quad\text{for}~ s\le 7,\quad
\partial_{x_3}^{s}({\bf v}_{1}^{3})^{(3)}=0\quad\text{for}~s\ge 
8, \\
&|\nabla_{x'}^{k}\partial_{x_3}^{s}\bar{p}_{3}|\le 
C\delta(x')^{\frac32-s-\frac{k}{2}}\quad\text{for}~0\leq s\le6,\quad 
\partial_{x_3}^{s}\bar{p}_{3}=0\quad\text{for}~s\ge 7.
\end{split}
\end{align}
Denote ${\bf v}^{3}={\bf v}_{1}^{1}+{\bf v}_{1}^{2}+{\bf v}_{1}^{3},\,{\bar{p}
^3}=\bar{p}_{1}+\bar{p}_{2}+\bar{p}_{3}.$ Recall that $\nabla\cdot{\bf v}^{3}
=R(x')$ with $|R(x')|\le C\delta(x')^{3/2}.$ By virtue of \eqref{f13esti} 
and Proposition \ref{keyprop1}, we have for any $m\ge1,$
\begin{align}\label{u11-v3est}
&\|\nabla^{m+1}({\bf u}_{1}^{1} -{\bf v}^{3})\|_{L^{\infty}(\Omega_{\delta(x')/
2}(x'))} +\|\nabla^m(p_{1}^{1}-\bar{p}^{3})\|_{L^{\infty}(\Omega_{\delta(x')/2}
(x'))} \nonumber\\
&\le C(\delta(x')^{2-m} + \delta(x')^{\frac32-m})\le C\delta(x')^{\frac{3}
{2}-m}.
\end{align}
By \eqref{v13barp3est}, we summarize
\begin{equation}\label{v13est}
\begin{split}
&|\nabla^{m+1}{\bf v}_{1}^{3}|\le C\delta(x')^{1-m}
\quad\text{for}~0\le m\le 4,\,\,|\nabla^{m+1}{\bf v}_{1}^{3}|\le C\delta(x')^{-\frac{m+3}{2}}
\quad\text{for}~m\ge5 ,\\
&|\nabla^{m}\bar{p}_{3}|\le C\delta(x')^{\frac32-m}\quad\quad\text{for}~0\le 
m \le 6,\,\,|\nabla^{m}\bar{p}_{3}|\le C \delta(x')^{-\frac{m+3}{2}}\qquad\text{for}
~ m\ge7 .
\end{split}
\end{equation}
By \eqref{estinabv11barp}, \eqref{v12est}, and \eqref{v13est}, we have
\begin{equation}\label{v3m+1deri}
|\nabla^{m+1}{\bf v}^{3}|+|\nabla^{m}\bar{p}^3|\le C \delta(x')^{-\frac{m+3}
{2}}.
\end{equation}
It follows from \eqref{u11-v3est} and \eqref{v3m+1deri} that for $m\ge1$,
\begin{equation*}
|\nabla^{m+1}{\bf u}_{1}^{1}|+|\nabla^{m}p_{1}^{1}|\le C \delta(x')^{-\frac{m+3}
{2}}+\delta(x')^{\frac32-m},
\end{equation*}
which implies that in $\Omega_{R},$
\begin{equation*}
\begin{aligned}
&|\nabla^{m+1}{\bf u}_{1}^{1}|+|\nabla^{m}p_{1}^{1}|\le C \delta(x')^{-
\frac{m+3}{2}}\quad\text{for}~m\le6,\\
&|\nabla^{m+1}{\bf u}_{1}^{1}|+|\nabla^{m}p_{1}^{1}|\le C
\delta(x')^{\frac32-m}\quad\,\,\text{for}~m\ge7.
\end{aligned}
\end{equation*}
This completes the proof of Proposition \ref{u11-general}.
\end{proof}

\begin{remark}
We would like to emphasize that, unlike in \cite{DLTZ1} for the two-dimensional case, Proposition \ref{U_esti} is essential for studying three or higher-dimensional problems.
For example, $F(x')=-\frac{3(h_1-h_2)(x')}{\delta(x')}$ and  $G(x')=0$ is the solution of \eqref{2by1}, but it is hard to find an $\bar{p}_1$ such that $|\mu\Delta{\bf v}_1^1-\nabla\bar{p}_1|\leq C\delta(x')^{-1}$.
\end{remark}

\section{Estimates for \texorpdfstring{$({\bf u}_{1}^{3},p_{1}^{3})$}{}}\label{sec4}
In this section, we derive higher-order derivatives of \({\bf u}_1^3\) and \(p_1^3\) using the method outlined in Proposition \ref{u11-general}. We reduce the construction of auxiliary functions to studying a class of ordinary differential equations, thereby fully capture the singular behavior of \(\nabla^{m+1}{\bf u}_1^3\) and \(\nabla^mp_1^3\). 

\begin{prop}\label{u13-esti}
Under the same assumption as in Theorem \ref{main thm1}, let ${\bf u}_{1}^{3}$ 
and $p_{1}^{3}$ be the solution to \eqref{u,peq1} with $\alpha=3$. Then for 
sufficiently small $0<\varepsilon<1/2$, we have
$$\mbox{(i)}\quad~|\nabla{\bf u}_{1}^{3}(x)|\le C(\varepsilon+|x'|^2)^{-
\frac32},\quad |p_1^3(x)-p_1^3(z',0)|\leq C(\varepsilon+|x'|^2)^{-2},\quad 
x\in\Omega_{R},\quad$$ 
for some point $z=(z',0)\in \Omega_{R}$, $|z'|=R/2$; and 

(ii) for any $m\ge 1$, 
\begin{equation*}
|\nabla^{m+1}{\bf u}_{1}^{3}(x)|+|\nabla^{m}{p}_{1}^{3}(x)| \le C(\varepsilon+|
x'|^2)^{-\frac{m+4}{2}},\quad x\in\Omega_{R}.
\end{equation*}
\end{prop}

\begin{proof}[Proof of Proposition \ref{u13-esti}]
The proof is similar to that of Proposition \ref{u11-general}. To show the main idea more clearly, we divide the proof into three steps.

\subsection{Gradient estimate.} 
In order to establish the gradient estimate 
of ${\bf u}_1^3$, we first construct a divergence-free auxiliary function $
{\bf v}_{3}^{1}$ which satisfies the boundary condition ${\bf v}_{3}^{1}(x)={\bf 
u}_{1}^{3}(x)=\boldsymbol{\psi}_{3}$ on $\Gamma^{+}_{2R}$, ${\bf v}_{3}^{1}
(x)={\bf u}_{1}^{3}(x)=0$ on $\Gamma^{-}_{2R}$, and a corresponding $\bar{p}_{1}
$, such that
$$|\mu\Delta{\bf v}_{3}^{1}-\nabla\bar{p}_{1}|\leq C\delta(x')^{-1}.$$
However, $\boldsymbol{\psi}_{3}=(0,0,1)^{T}$ introduce in a higher singularity 
than $\boldsymbol{\psi}_{1}=(1,0,0)^{T}$ and $\boldsymbol{\psi}_{2}=(0,1,0)^{T}.
$ We divide the process into two steps to construct
$${\bf v}_{2}^{1}:=\tilde{\bf v}_{2}^{1}+\hat{\bf v}_{2}^{1}\quad\mbox{and}\quad 
\bar{p}_{1}:=\tilde{p}_{1}+\hat{p}_{1}.$$

{\bf Step I.1. Construction of $\tilde{\bf v}_{3}^{1}$ and $\tilde{p}_{1}$.} 
We begin by constructing an auxiliary function $\tilde{\bf v}_{3}^{1}$ that 
satisfies the boundary conditions $\tilde{\bf v}_{3}^{1}(x)={\bf u}_{1}^{3}
(x)=\boldsymbol{\psi}_{3}$ on $\Gamma^{+}_{2R}$ and $\tilde{\bf v}_{3}^{1}
(x)={\bf u}_{1}^{3}(x)=0$ on $\Gamma^{-}_{2R}$. The construction process of $
\tilde{\bf v}_{3}^{1}(x)$ is similar to that of ${\bf v}_{1}^{1}(x)$ in 
Proposition \ref{u11-general}, with some minor modifications. Specifically, we 
assume that $\tilde{\bf v}_{3}^{1}(x)$ takes the following form:
\begin{equation*}
\tilde{\bf v}_{3}^{1}(x)=
\begin{pmatrix}
0\\0\\k(x)+\frac{1}{2}
\end{pmatrix}+ 
\begin{pmatrix}
F(x)\\
G(x)\\
H(x)\\
\end{pmatrix}\big(k(x)^2-\frac14\big)\quad\text{in}~\Omega_{2R},
\end{equation*}
where $k(x)=\frac{x_3-\frac12(h_1-h_2)(x')}{\delta(x')}$ is defined in 
\eqref{def_k}.
By a direct calculation, we have
\begin{equation}\label{u12fz1}
\begin{split}
&\partial_{x_{1}}(\tilde{\bf v}_{3}^{1})^{(1)}=\partial_{x_{1}} F(x)\big(k(x)^2-\frac14\big)+2F(x)k(x)\partial_{x_{1}} k(x),\\
&\partial_{x_{2}}(\tilde{\bf v}_{3}^{1})^{(2)}=\partial_{x_{2}} G(x)\big(k(x)^2-\frac14\big)+2G(x)k(x)\partial_{x_{2}} k(x),\\
&\partial_{x_{3}}(\tilde{\bf v}_{3}^{1})^{(3)}=\frac{1}{\delta(x')}+
\partial_{x_{3}} H(x)\Big(k(x)^2-
\frac{1}{4}\Big)+\frac{2k(x)}{\delta(x')}H(x).
\end{split}
\end{equation}
To determine $F(x)$, $G(x)$, and $H(x)$, instead of \eqref{ffz1}, we set 
\begin{equation}\label{u12fz2}
\frac14(\partial_{x_{1}}F+\partial_{x_{2}}G+\partial_{x_{3}}H)(x)=\frac{1}
{\delta(x')}.
\end{equation}
By using \eqref{u12fz1}, \eqref{u12fz2}, and the divergence-free condition $\nabla\cdot {\bf v}=0$, we obtain
\begin{equation}\label{u12fz3}
H(x)=-2k(x)-F(x)\delta(x')\partial_{x_{1}}k(x)-G(x)\delta(x')\partial_{x_{2}}k(x).
\end{equation}

By \eqref{u12fz2} and \eqref{u12fz3}, we conclude that $F(x)$ and $G(x)$ 
satisfy the following first-order differential equation:
\begin{align}\label{1-st order eq31}
\begin{split}
\partial_{x_{1}}F(x)+\partial_{x_{2}}G(x)&+\frac{\partial_{x_{1}}\delta(x')}
{\delta(x')}F(x)+\frac{\partial_{x_{2}}\delta(x')}{\delta(x')}G(x)\\
&-\delta(x')\big(\partial_{x_1}k(x)\partial_{x_{3}}F(x)+\partial_{x_2}k(x)
\partial_{x_{3}}G(x)\big)=\frac{6}{\delta(x')}.
\end{split}
\end{align}
We assume $F(x)=F(x')$ and $ G(x)=G(x')$. Then \eqref{1-st order eq31} becomes
\begin{equation}\label{u13-FG-eq}
\partial_{x_{1}}F(x)+\partial_{x_{2}}G(x) +\frac{\partial_{x_{1}}\delta(x')}
{\delta(x')}F(x)+\frac{\partial_{x_{2}}\delta(x')}{\delta(x')}G(x)=\frac{6}
{\delta(x')}.
\end{equation}
Similar to \eqref{tj1} and \eqref{zyfc1}, we choose ${\tilde{p}_1^1}(x')$ such that 
\begin{equation}\label{cancel41}
\partial_{x_{1}}\tilde{p}_{1}^1(x')=\frac{2\mu F(x')}{\delta(x')^2},\quad
\partial_{x_{2}}\tilde{p}_{1}^1(x')=\frac{2\mu G(x')}{\delta(x')^2}.
\end{equation}
Substituting \eqref{cancel41} into \eqref{u13-FG-eq} yields
\begin{equation}\label{zyfc41}
\Delta_{x'} \tilde{p}_{1}^1(x')+\frac{3\partial_{x_1} \delta(x')}
{\delta(x')}\partial_{x_1}\tilde{p}_{1}^1(x')+\frac{3\partial_{x_2} 
\delta(x')}
{\delta(x')}\partial_{x_2}\tilde{p}_{1}^1(x')=\frac{12\mu}{\delta(x')^3}.
\end{equation}
We rewrite \eqref{zyfc41} as 
\begin{equation}\label{tilp41radeq}
\tilde{p}_{1}^{1}(r)''+\Big(\frac{3\delta'(r)}{\delta(r)}+\frac{1}{r}\Big)
\tilde{p}_{1}^{1}(r)'=\frac{12\mu}{\delta(r)^3}.
\end{equation}
Then 
\begin{equation}\label{tilpform}
\tilde{p}_{1}^{1}(x')=-6\mu\int_{r}^{R}\frac{s}{\delta(s)^3}\ ds
\end{equation}
is a particular solution, which is a radial function. By \eqref{cancel41}, we have 
\begin{equation}\label{FG31form}
F(x')=\frac{\delta(x')^2}{2\mu}\partial_{x_1}\tilde{p}_{1}^{1}=\frac{3x_1}
{\delta(x')},\quad G(x')=\frac{\delta(x')^2}{2\mu}\partial_{x_2}\tilde{p}
_{1}^{1}=\frac{3x_2}{\delta(x')}.
\end{equation}
Substituting them into \eqref{u12fz3}, we have
\begin{equation*}
H(x)=-2k(x)-3x_1\partial_{x_1}k(x)-3x_2\partial_{x_2}
k(x)=-2k(x)-3r\partial_{r}k(x),
\end{equation*}
which is radial about $x'.$
Therefore, we have 
\begin{equation}\label{tildev21form}
\tilde{\bf v}_{3}^{1}(x)=
\begin{pmatrix}
0\\0\\k(x)+\frac{1}{2}
\end{pmatrix}+ 
\begin{pmatrix}
\frac{3x_1}{\delta(x')}\\
\frac{3x_2}{\delta(x')}\\
-2k(x)-3x_1\partial_{x_1}k(x)-3x_2\partial_{x_2}k(x)\\
\end{pmatrix}\big(k(x)^2-\frac14\big)\quad\text{in}~\Omega_{2R},
\end{equation}
which satisfies the divergence-free condition condition $\nabla\cdot\tilde{\bf 
v}_{3}^{1}=0$ in $\Omega_{2R}$.

By \eqref{h1h2prop} and a direct calculation, we get 
\begin{align*}
&|\nabla_{x'}(\tilde{\bf v}_{3}^{1})^{(i)}|\le C\delta(x')^{-1},\quad|
\partial_{x_{3}}(\tilde{\bf v}_{3}^{1})^{(i)}|\le C\delta(x')^{-3/2}, \,\, i=1,2,
\nonumber\\
&|\nabla_{x'}(\tilde{\bf v}_{3}^{1})^{(3)}|\le C\delta(x')^{-1/2},\quad|
\partial_{x_{3}}(\tilde{\bf v}_{3}^{1})^{(3)}|\le C \delta(x')^{-1},
\end{align*}
and for any $k\ge 1,$ and $ i=1,2$,
\begin{align}\label{tilv31highesti}
\begin{split}
&|\nabla_{x'}^{k}\partial_{x_{3}}^{s}(\tilde{\bf v}_{3}^{1})^{(i)}|\le 
C\delta(x')^{-\frac{1+k}{2}-s} \,\,\, \text{for}\,\,0\le s\le 2, \quad
\partial_{x_{3}}^{s}(\tilde{\bf v}_{3}^{1})^{(i)}=0 \,\,\, \text{for}\,
\,s\ge3,\\ 
&|\nabla_{x'}^{k}\partial_{x_{3}}^{s}(\tilde{\bf v}_{3}^{1})^{(3)}|\le 
C\delta(x')^{-\frac{k}{2}-s} \,\,\,\,\,\,\,\, \text{for}\,\,0\le s\le 3, \quad
\partial_{x_{3}}^{s}(\tilde{\bf v}_{3}^{1})^{(3)}=0 \,\,\, \text{for}\,
\,s\ge4.
\end{split}
\end{align}
By \eqref{tilv31highesti}, one can see that the leading term in $\Delta(\tilde{\bf v}_{3}^{1})^{(3)}$ is 
$\partial_{x_3x_3}(\tilde{\bf v}_{3}^{1})^{(3)}.$ To further cancel this leading term, we take
\begin{equation}\label{tilp31def}
{\tilde{p}_1}(x):={\tilde{p}_1^1}(x')+\mu\partial_{x_{3}}\Big(H(x)\big(k(x)^2-\frac14\big)\Big)={\tilde{p}_1^1}(x')+\mu\partial_{x_{3}}(\tilde{\bf v}_3^1)^{(3)}-\frac{\mu}{\delta(x')}.
\end{equation}
By \eqref{tildev21form} and the definition of $k(x)$, we see that $(\tilde{\bf v}_{3}^{1})^{(i)},i=1,2$, is a quadratic polynomial in $x_3$, while $(\tilde{\bf v}_{3}^{1})^{(3)}$ is a cubic polynomial in $x_3$. By \eqref{cancel41}, 
\eqref{FG31form}, \eqref{tildev21form}, \eqref{tilv31highesti}, and 
\eqref{tilp31def}, a direct calculation gives
\begin{align}\label{muDelta-gradp1}
&\mu\Delta\tilde{\bf v}_{3}^{1}-\nabla\tilde{p_{1}}\nonumber\\
&=
\mu\begin{pmatrix}
\Delta_{x'}(\tilde{\bf v}_{3}^{1})^{(1)}-\partial_{x_{1}x_{3}}
(\tilde{\bf{v}}_{3}^{1})^{(3)}+\partial_{x_{1}}\frac{1}{\delta(x')}\\\\
\Delta_{x'}(\tilde{\bf v}_{3}^{1})^{(2)}-\partial_{x_{2}x_{3}}
(\tilde{\bf{v}}_{3}^{1})^{(3)}+\partial_{x_{2}}\frac{1}{\delta(x')}\\\\
\Delta_{x'}(\tilde{\bf v}_{3}^{1})^{(3)}
\end{pmatrix}:= \begin{pmatrix}
A_{1}(x')x_{3}^{2}+B_{1}(x')x_{3}+C_{1}(x') \\\\
{A}_{2}(x')x_{3}^{2}+{B}_{2}(x')x_{3}+{C}_{2}(x')\\\\
\hat{A}(x')x_{3}^{3}+\hat{B}(x')x_{3}^{2}+\hat{C}(x')x_{3}+\hat{D}(x')
\end{pmatrix},
\end{align}
where for $0\le s\le m,$ in $\Omega_{2R},$
\begin{equation}\label{ABCDest}
\begin{split}
&|\nabla_{x'}^{s}A_{i}|\le C\delta(x')^{-\frac{7+s}{2}},\,\,|\nabla_{x'}^{s}
B_{i}|\le C\delta(x')^{-\frac{5+s}{2}},\,\,|\nabla_{x'}^{s}C_{i}|\le 
C\delta(x')^{-\frac{3+s}{2}},~i=1,2, \\
&|\nabla_{x'}^{s}\hat{A}|\le C\delta(x')^{-\frac{8+s}{2}},\,\,|
\nabla_{x'}^{s}\hat{B}|\le C\delta(x')^{-\frac{6+s}{2}},\\
&|\nabla_{x'}^{s}\hat{C}|\le C\delta(x')^{-\frac{4+s}{2}},\,\,|
\nabla_{x'}^{s}\hat{D}|\le C\delta(x')^{-\frac{2+s}{2}}.
\end{split}
\end{equation}
 Recalling the definition of $\tilde{\bf v}_{3}^{1}$ in \eqref{tildev21form}, by \eqref{muDelta-gradp1}, a direct calculation shows, for $i=1,2$, 
\begin{equation}\label{ys1c2}
 A_{i}=x_{i}A(r),\quad B_{i}=x_{i}B(r), \quad C_{i}=x_{i}C(r),
\end{equation}
for some radial functions $A(r)$, $B(r)$, and $C(r)$.
By using \eqref{muDelta-gradp1} and \eqref{ABCDest}, we have 
\begin{equation}\label{zj2}
|(\mu\Delta\tilde{\bf v}_{3}^{1}-\nabla\tilde{p}_{1})^{(i)}|\leq 
C\delta(x')^{-3/2},~i=1,2,\quad |(\mu\Delta\tilde{\bf v}_{3}^{1}-
\nabla\tilde{p}_{1})^{(3)}|\leq C\delta(x')^{-1}.
\end{equation}
Thus, we have $|\mu\Delta\tilde{\bf v}_{3}^{1}-\nabla\tilde{p}_{1}|\le 
C\delta(x')^{-3/2}$. 

However, this bound is too large to apply Proposition 
\ref{keyprop1}. Therefore, in order to fully capture the singularity in $
\nabla{\bf u}_1^3$, we need to introduce additional divergence-free correctors 
to reduce the upper bound, $\delta(x')^{-3/2}$.%

{\bf Step I.2. Construction of $\hat{\bf v}_{3}^{1}$ and $\hat{p}_{1}$.}
It follows from \eqref{zj2} that the leading terms in $\mu\Delta{\bf \tilde{v}}_{3}^1-\nabla\tilde{p}_{1}$ are $(\mu\Delta{\bf \tilde{v}}_{3}^1-\nabla\tilde{p}_{1})^{(i)},i=1,2$. To cancel out them, similar to the process of Step II in the proof of Proposition \ref{u11-general}, we choose the corrector $\hat{{\bf v}}_{3}^{1}$ in the following form in $\Omega_{2R}$:
\begin{align}\label{hatv21def}
\hat{{\bf v}}_{3}^{1}(x)&=
\begin{pmatrix}
F_{11}^{2}(x')x_3^{2}+F_{11}^{1}(x')x_3+F_{11}^{0}(x')+\tilde{F}
_{11}(x') \\\\
F_{21}^{2}(x')x_3^2+F_{21}^{1}(x')x_3+F_{21}^{0}(x')+\tilde{F}_{21}(x')\\\\
F_{31}^{3}(x')x_3^{3}+F_{31}^{2}(x')x_3^2+F_{31}^{1}(x')x_3+F_{31}^{0}(x')+
\tilde{F}_{31}(x)
\end{pmatrix}\big(k(x)^2-\frac14\big) \nonumber \\
 &:= \begin{pmatrix}
 F_{11}(x)+\tilde{F}_{11}(x') \\
 F_{21}(x)+\tilde{F}_{21}(x')\\
 F_{31}(x)+\tilde{F}_{31}(x)
 \end{pmatrix}\big(k(x)^2-\frac14\big).
\end{align}

We first construct $F_{11}(x)$, $F_{21}(x)$, and $F_{31}(x)$. Here $F_{11}(x)$ and $F_{21}(x)$ are chosen to cancel out $(\mu\Delta{\bf \tilde{v}}_{3}^1-\nabla\tilde{p}_{1})^{(i)},i=1,2,$ in the following way,
\begin{equation}\label{part22cancel}
\mu \partial_{x_3x_3}\Big(F_{i1}\Big(k(x)^2-
\frac14\Big)\Big)=-(\mu\Delta\tilde{\bf v}_{3}^{1}-\nabla\tilde{p}
_{1})^{(i)},\, i=1,2.
\end{equation}
By comparing the coefficients of each term of $x_3$, we obtain
\begin{equation*}
\begin{split}
&F_{i1}^{2}(x')=-\frac{\delta(x')^2}{12\mu}A_{i}(x'),\,\, F_{i1}^{1}(x')=-
\frac{\delta(x')^{2}}{6\mu}B_{i}(x')+(h_{1}-h_{2})F_{i1}^{2}(x'), \\ 
&F_{i1}^{0}(x')=-\frac{\delta(x')^2}{2\mu}C_{i}(x')+(h_{1}-h_{2})F_{i1}
^{1}(x')+\frac14(\varepsilon+2h_1(x'))(\varepsilon+2h_2(x'))F_{i1}^{2}
(x'),~i=1,2.
\end{split}
\end{equation*}
 Combining this and \eqref{ys1c2}, we have, for $j=0,1,2$,
\begin{equation}\label{ys2c2}
 F_{i1}^{j}(x')=x_iF^j(r)
\end{equation}
for some radial functions $F^j(r)$. We further take 
\begin{align*}
&F_{31}^{3}(x')=-\frac{\delta(x')^{2}}{5} \sum_{i=1}^{2}\partial_{x_i}
\Big(\frac{F_{i1}^{2}(x')}{\delta(x')^{2}}\Big) ,\\ 
&F_{31}^{2}(x')=\frac{\delta(x')^2}{4}
\sum_{i=1}^{2}\partial_{x_{i}}\Big(\frac{(h_1-h_2)(x')F_{i1}^{2}(x')-
F_{i1}^{1}(x')}{\delta(x')^2}\Big)+(h_1-h_2)(x')F_{31}^{3}(x'),\\
&F_{31}^{1}(x')=\frac{\delta(x')^2}{3}\Big(\sum_{i=1}^{2}\partial_{x_{i}}
\Big(\frac{(h_1-h_2)(x')F_{i1}^{1}(x')-F_{i1}^{0}(x') }{\delta(x')^2}\Big) 
\nonumber\\
&\qquad\qquad+\sum_{i=1}^{2}\partial_{x_{i}}
\Big(\frac{(\varepsilon+2h_1(x'))(\varepsilon+2h_2(x'))}{4\delta(x')^2}
F_{i1}^{2}(x')\Big)+(h_1-h_2)(x')F_{31}^{2}(x')\nonumber\\ &\qquad\qquad+
\frac14(\varepsilon+2h_1(x'))(\varepsilon+2h_2(x'))F_{31}^{3}(x'),\\
&F_{31}^{0}(x')=\frac{\delta(x')^2}{2}\sum_{i=1}^{2}\Big(\partial_{x_{i}}
\Big(\frac{(h_1-h_2)(x')F_{i1}^{0}(x')}{\delta(x')^2}\Big)+\partial_{x_{i}}
\Big(\frac{(\varepsilon+2h_1(x'))(\varepsilon+2h_2(x'))}{4\delta(x')^2}
F_{i1}^{1}(x')\Big) \Big) \\
&\qquad\qquad+(h_1-h_2)(x')F_{31}^{1}(x')+\frac14(\varepsilon+2h_1(x'))
(\varepsilon+2h_2(x'))F_{31}^{2}(x'),
\end{align*} 
satisfying
\begin{equation*}
\nabla\cdot
\begin{pmatrix}
F_{11}\big(k(x)^2-\frac14\big) \\
F_{21}\big(k(x)^2-\frac14\big)\\
F_{31}\big(k(x)^2-\frac14\big)
\end{pmatrix}:=R(x'),
\end{equation*}
where
\begin{align*}
R(x')=-\frac{h_{1}-h_{2}}{\delta(x')^2}F_{31}^{0}(x')-\sum_{i=1}^{2}
\partial_{x_{i}}\Big(\frac{(\varepsilon+2h_{1})(\varepsilon+2h_{2})}
{4\delta(x')^2}F_{i1}^{0}(x')\Big)-\frac{(\varepsilon+2h_{1})
(\varepsilon+2h_{2})}{4\delta(x')^2}F_{31}^{1}(x') 
\end{align*}
is independent of $x_3.$ It follows from \eqref{ys2c2} that $F_{31}^{j}(x')=F_{31}^{j}(r)$. Recalling the definition of $F_{i1}$, $i=1,2$, by using \eqref{ys2c2} again, we have $\partial_{x_1}F_{11}+\partial_{x_2}F_{21}$ is a radial function. Hence $R(x')$ is radial. By using \eqref{h1h2prop} and \eqref{ABCDest}, we have, for 
$k\ge0, i=0,1,2,$ and $j=0,1,2,3,$
\begin{equation}\label{F1121est}
|\nabla_{x'}^{k}F_{11}^{i}|,|\nabla_{x'}^{k}F_{21}^{i}|\le C\delta(x')^{\frac12-
i-\frac{k}{2}},\,\,|\nabla_{x'}^kF_{31}^{j}|\le C\delta(x')^{1-j-\frac{k}{2}},
\,|\nabla_{x'}^{k}R|\le C\delta(x')^{-\frac{k}{2}}.
\end{equation}

To ensure $\nabla\cdot\hat{{\bf v}}_{3}^{1}=0,$ we now choose $\tilde{F}
_{11}(x'),\tilde{F}_{21}(x'),$ and $\tilde{F}_{31}(x)$ such that 
\begin{align}\label{tilF11divfree}
 \nabla\cdot
\begin{pmatrix}
\tilde{F}_{11}(x')\big(k(x)^2-\frac14\big) \\
\tilde{F}_{21}(x')\big(k(x)^2-\frac14\big)\\
\tilde{F}_{31}(x)\big(k(x)^2-\frac14\big)
\end{pmatrix}=-R(x').
\end{align}
After a direct calculation, \eqref{tilF11divfree} becomes
\begin{align}\label{tilF11eq}
 &\partial_{x_1}\tilde{F}_{11}\big(k(x)^2-\frac14\big) +\partial_{x_2}\tilde{F}
 _{21}\big(k(x)^2-\frac14\big)+2k(x)\partial_{x_{1}}k(x')\tilde{F}
 _{11}(x')\nonumber\\
&+2k(x')\partial_{x_{2}}k(x)\tilde{F}_{21}(x') + \partial_{x_3}
 \tilde{F}_{31}\big(k(x)^2-\frac14\big)+ \frac{2k(x)\tilde{F}_{31}(x)}
 {\delta(x')} +R(x')=0.
\end{align}
To determine $\tilde{F}
_{11}(x'),\tilde{F}_{21}(x')$, and $\tilde{F}_{31}(x)$, we let 
\begin{equation}\label{divfgheqR}
\partial_{x_1}\tilde{F}_{11}(x')+\partial_{x_2}\tilde{F}_{21}(x')+
\partial_{x_3}\tilde{F}_{31}(x)=4R(x').
\end{equation}
It follows from \eqref{tilF11eq} and \eqref{divfgheqR} that 
\begin{equation}\label{tilF31form}
\tilde{F}_{31}(x)=-2\delta(x')k(x)R(x')-\delta(x')\partial_{x_1}k(x)\tilde{F}
_{11}(x')-\delta(x')\partial_{x_2}k(x)\tilde{F}_{21}(x').
\end{equation}
Subtituting \eqref{tilF31form} into \eqref{divfgheqR}, we get 
\begin{equation}\label{tilFrelR}
\partial_{x_1}\tilde{F}_{11}(x')+\frac{\partial_{x_1}\delta(x')}{\delta(x')}
\tilde{F}_{11}(x')+\partial_{x_2}\tilde{F}_{21}(x')+\frac{\partial_{x_2}
\delta(x')}{\delta(x')}\tilde{F}_{21}(x')= 6R(x').
\end{equation}
Similar to \eqref{2by1}-\eqref{zyfc1}, to determine $\tilde{F}_{11}(x')$ and $
\tilde{F}_{21}(x')$, we choose $\hat{p}_{1}(x')$ such that 
\begin{equation}\label{hatp1retilF}
\partial_{x_1}\hat{p}_{1}=\frac{2\mu\tilde{F}_{11}(x')}{\delta(x')^2},
\quad\partial_{x_2}\hat{p}_{1}=\frac{2\mu\tilde{F}_{21}(x')}{\delta(x')^2}.
\end{equation}
Combining \eqref{tilFrelR} and \eqref{hatp1retilF}, we have
\begin{equation}\label{hatp1eq}
\Delta_{x'}\hat{p}_{1}+\frac{3\partial_{x_1}\delta(x')}{\delta(x')}
\partial_{x_1}\hat{p}_{1} + \frac{3\partial_{x_2}\delta(x')}{\delta(x')}
\partial_{x_2}\hat{p}_{1}=\frac{12\mu R(x')}{\delta(x')^2}.
\end{equation}
Since $R(x')$ is radial, by a similar argument that led to \eqref{tilp41radeq}, 
we obtain
\begin{equation}\label{hatp1form}
\hat{p}_{1}(x')=-12\mu\int_{r}^{R}\frac{1}{s\delta(s)^3}\int_{0}^{s}
R(t)t\delta(t)\, dtds.
\end{equation}
Thus by \eqref{hatp1retilF}, we have
\begin{align}\label{tilF11}
\begin{split}
\tilde{F}_{11}(x')=\frac{\delta(x')^2}{2\mu}\partial_{x_1}\hat{p}_{1}
=\frac{6x_{1}}{r^2\delta(r)}\int_{0}^{r}R(t)t\delta(t)\, dt, \\
\tilde{F}_{21}(x')=\frac{\delta(x')^2}{2\mu}\partial_{x_2}\hat{p}_{1}
=\frac{6x_2}{r^2\delta(r)}\int_{0}^{r}R(t)t\delta(t) \, dt.
\end{split}
\end{align}
Substituting them into \eqref{tilF31form} yields
\begin{align}\label{tilF31form1}
\tilde{F}_{31}(x)&=-2\delta(x')k(x)R(x')-\frac{6x_{1}\partial_{x_1}k(x)}{r^2}
\int_{0}^{r}R(t)t\delta(t)\, dt-\frac{6x_{1}\partial_{x_1}k(x)}{r^2}\int_{0}^{r}
R(t)t\delta(t)\, dt \nonumber\\
&=-2\delta(x')k(x)R(x')-\frac{6\partial_{r}k(x)}{r}\int_{0}^{r}
R(t)t\delta(t)\, dt:=P(x')x_3+Q(x'),
\end{align}
which is radial with respect to $x'$. 
By \eqref{F1121est}, \eqref{hatp1form}, \eqref{tilF11}, and \eqref{tilF31form1}, 
a calculation gives, for any $k\ge0,$
\begin{align}\label{tilF11esti}
|\nabla_{x'}^{k}\tilde{F}_{11}|,|\nabla_{x'}^{k}\tilde{F}_{21}|\le 
C\delta(x')^{\frac{1-k}{2}},\,\,|\nabla_{x'}^kP|\le C\delta(x')^{-\frac{k}{2}},|
\nabla_{x'}^kQ|\le C\delta(x')^{1-\frac{k}{2}}\quad\text{in}
~\Omega_{2R}.
\end{align}
Now we define ${\bf v}_{3}^{1}(x)=\tilde{{\bf v} }_{3}^{1}(x)+\hat{{\bf v}}_{3}
^{1}(x)$ and $\bar{p}_{1}=\tilde{p}_{1}+\hat{p}_{1},$ which satisfy $
\nabla\cdot{\bf v}_{3}^{1}(x)=\nabla\cdot\tilde{{\bf v} }_{3}^{1}(x)+
\nabla\cdot\hat{{\bf v}}_{3}^{1}(x)=0 $ and $\partial_{x_3}\hat{p}_{1}=0$, respectively.

{\bf Step I.3. Estimates of ${\bf f}_3^1$.}
 Denote 
 $${\bf f}_{3}^{1}:=\mu\Delta{\bf v}_{3}^{1}-\nabla\bar{p}_{1}\quad\text{in}
 ~\Omega_{2R}.$$ 
By virtue of \eqref{muDelta-gradp1}, \eqref{hatv21def}, \eqref{part22cancel}, 
and \eqref{hatp1retilF}, we have 
\begin{equation}\label{f31form}
{\bf f}_{3}^{1}=\begin{pmatrix}
\mu\Delta_{x'}(\hat{\bf v}_{3}^{1})^{(1)}(x)\\\\
\mu\Delta_{x'}(\hat{\bf v}_{3}^{1})^{(2)}(x)\\\\
\mu \Delta_{x'}(\hat{\bf v}_{3}^{1})^{(3)}(x)+\mu\Delta_{x'}(\tilde{\bf v}
_{3}^{1})^{(3)}(x)+\mu\partial_{x_3x_3}
(\hat{\bf v}_{3}^{1})^{(3)}(x)
\end{pmatrix}.
\end{equation} 
It follows from \eqref{hatv21def} that $(\hat{\bf v}_3^1)^{(i)},i=1,2$, is a 
polynomial of order $4$ in $x_3$, while $(\hat{\bf v}_3^1)^{(3)}$ is a polynomial of order $5$ in $x_3$. In view of \eqref{muDelta-gradp1}, we can rewrite 
\begin{equation}\label{f311,f312}
({\bf f}_{3}^{1})^{(j)}=\sum_{i=0}^{4}H^{j}_{i1}(x')x_{3}^{i},\, j=1,2,
\quad({\bf f}_{3}^{1})^{(3)}=\sum_{i=0}^{3}S_{i1}(x')x_{3}^{i}+
\sum_{i=0}^{5}G_{i1}(x')x_{3}^{i}\quad \text{in}~\Omega_{2R},
\end{equation}
where 
\begin{equation*}
\sum_{i=0}^{3}S_{i1}(x')x_3^{i}:= \mu\Delta_{x'}(\tilde{\bf v}
_{3}^{1})^{(3)}+\mu\partial_{x_3x_3}(\hat{\bf v}
_{3}^{1})^{(3)},\quad
\sum_{i=0}^{5}G_{i1}(x')x_{3}^{i}:=\mu\Delta_{x'}
(\hat{\bf v}_{3}^{1})^{(3)}.
\end{equation*}
By using \eqref{ABCDest}, \eqref{F1121est}, and \eqref{tilF11esti}, for 
$s\ge0,$ $1\le i\le 4,$ $0\le j\le3$, and $0\le k\le5,$
\begin{equation}\label{GH1esti}
|\nabla_{x'}^{s}H_{i1}^{1}|, |\nabla_{x'}^{s}H_{i1}^{2}| \le C\delta(x')^{-
\frac12-i-\frac{s}{2}},\quad|\nabla_{x'}^{s}S_{j1}|\le C \delta(x')^{-1-j-
\frac{s}{2}},\, |\nabla_{x'}^sG_{k1}|\le C \delta(x')^{-k-\frac{s}{2}}.
\end{equation}
Therefore, by \eqref{f311,f312} and \eqref{GH1esti}, we have
\begin{equation*}
|({\bf f}_{3}^{1})^{(1)}|,|({\bf f}_{3}^{1})^{(2)}| \le C\delta(x')^{-1/2}, 
\quad|({\bf f}_{3}^{1})^{(3)}|\le C\delta(x')^{-1}\quad \text{in}
~\Omega_{2R}.
\end{equation*}
Thus,
\begin{equation}\label{f31esti}
|{\bf f}_{3}^{1}|\le C (|({\bf f}_{3}^{1})^{(1)}|+ |({\bf f}_{3}^{1})^{(2)}| 
+|({\bf f}_{3}^{1})^{(3)}|) \le C \delta(x')^{-1}\quad \text{in}
~\Omega_{2R}.
\end{equation}

By applying Proposition \ref{keyprop2} and using \eqref{f31esti}, we have for 
any 
$x\in \Omega_{R}$,
$$
\| \nabla {\bf u}_{1}^{3}-\nabla{\bf v}_{1}^{3}\|_{L^{\infty}
(\Omega_{\delta(x_1)/2}(x'))}\le C. 
$$
By virtue of \eqref{tilp31def}, \eqref{tilpform}, \eqref{tilv31highesti}, 
\eqref{F1121est}, \eqref{hatp1form}, and \eqref{tilF11esti}
we have
\begin{equation*}
\begin{split}
&|\nabla^l {\bf v}_3^1|\le C\delta(x')^{-\frac{2l+1}{2}}\quad\text{for}
~l\le 2,
\quad |\nabla^l {\bf v}_3^1|\le C\delta(x_1)^{-\frac{l+3}{2}}\quad 
\text{for}~ l\ge 3,\\
& |\bar{p}_{1}|\le C\delta(x')^{-2},\quad|\nabla^l \bar{p}_1|\le 
C\delta(x')^{-\frac{l+4}{2}}\quad \text{for}~ l\ge 1.
\end{split}
\end{equation*}
Thus, Proposition \ref{u13-esti} holds for $m=0$, with
$$|\nabla {\bf u}_{1}^{3}|\le C|\nabla {\bf v}_{3}^{1}(x)| +C\le 
C\delta(x')^{-3/2}\quad\text{in}~\Omega_{R}.$$

\subsection{Higher-order derivatives estimates.}
Denote
\begin{equation}\label{def3j}
{\bf f}_{3}^{j}:={\bf f}_{3}^{j-1}+\mu\Delta{\bf v}_{3}^{j}(x)-\nabla\bar{p}_{j}
(x)=\sum_{l=1}^{j}(\mu\Delta{\bf v}_{3}^{l}-\nabla\bar{p}_{l})(x)\quad\text{for}
~ 2\le j\le m+1.
\end{equation}
We will inductively prove that for $j\ge 1$, ${\bf f}_{2}^{j}$ can be 
represented as polynomials in $x_3$ as follows:
\begin{equation}\label{f3j-form}
({\bf f}_{3}^{j})^{(b)}=\sum_{i=0}^{2j+2}H_{ij}^{b}(x')x_3^{i},\, b=1,2, \quad
({\bf f}_{3}^{j})^{(3)}=\sum_{i=0}^{2j+1}S_{ij}(x')x_3^{i}+\sum_{i=0}
^{2j+3}G_{ij}(x')x_{3}^{i}\quad\text{in}~\Omega_{2R},
\end{equation}
where for any $k\ge 0$, 
\begin{equation}\label{HSG-esti}
|\nabla_{x'}^{k}H_{ij}^{b}(x')|\le C\delta(x')^{j-i-\frac{k}{2}-\frac{3}{2}},\,|
\nabla_{x'}^kS_{ij}|\le C\delta(x')^{j-i-\frac{k}{2}-2},\,|\nabla_{x'}^kG_{ij}|
\le C \delta(x')^{j-i-\frac{k}{2}-1}.
\end{equation}
Indeed, by \eqref{f311,f312} and \eqref{GH1esti}, we have that \eqref{f3j-form} and 
\eqref{HSG-esti} hold for $j=1$. Assuming that \eqref{f3j-form} and \eqref{HSG-esti} hold for $j=l-1$ with $l\ge 2$, that is, in $\Omega_{2R}$, 
\begin{equation}\label{fzgn2}
\begin{aligned}
&({\bf f}_{3}^{(l-1)})^{(b)}=\sum_{i=0}^{2l}H_{i(l-1)}^{b}(x')x_3^{i}=\mu 
\Delta_{x'}({\bf v}_{3}^{l-1})^{(b)},~b=1,2,\\ 
&({\bf f}_{3}^{(l-1)})^{(3)}=\sum_{i=0}^{2l-1}S_{i(l-1)}(x')x_3^{i}+
\sum_{i=0}^{2l+1}G_{i(l-1)}(x')x_{3}^{i},
\end{aligned}
\end{equation}
where for any $k\ge 0$,
\begin{equation}\label{fzgn3}
|\nabla_{x'}^kH_{i(l-1)}^{b}|\le C\delta(x')^{l-i-\frac{k}{2}-\frac{5}{2}},\,|
\nabla_{x'}^kS_{i(l-1)}|\le C\delta(x')^{l-i-\frac{k}{2}-3},\,|\nabla_{x'}^{k}
G_{i(l-1)}|\le C \delta(x')^{l-i-\frac{k}{2}-2}.
\end{equation}

{\bf Step II.1. Construction of ${\bf v}_{3}^{l}$ and $\bar{p}_{l}$.} From \eqref{fzgn2} and \eqref{fzgn3}, we see that the leading term in ${\bf f}_{3}
^{l-1}$ is $\sum_{i=0}^{2l-1}S_{i(l-1)}(x')x_3^{i}$. To cancel this bad term, we take
\begin{equation*}
\tilde{p}_{l}=\sum_{i=0}^{2l-1}\frac{S_{i(l-1)}(x')}{i+1}x_{3}^{i+1}
\quad\text{in}~\Omega_{2R}
\end{equation*} 
such that 
\begin{equation}\label{f3l-1-nambatilp2}
{\bf f}_{3}^{l-1}-\nabla\tilde{p}_{l} =
\begin{pmatrix}
\mu\Delta_{x'}({\bf v}_{3}^{l-1})^{(1)}(x)-\partial_{x_1}\tilde{p}_{l}\\\\
\mu\Delta_{x'}({\bf v}_{3}^{l-1})^{(2)}(x)-\partial_{x_2}\tilde{p}_{l} \\\\
\sum_{i=0}^{2l+1}G_{i(l-1)}(x')x_{3}^{i}
\end{pmatrix}:=
\begin{pmatrix}
\sum_{i=0}^{4}A_{i1}^{1}(x')x_3^{i}\\\\
\sum_{i=0}^{4}A_{i1}^{2}(x')x_3^{i}\\\\
\sum_{i=0}^{5}B_{i1}(x')x_{3}^{i}
\end{pmatrix}.
\end{equation}
By \eqref{fzgn3}, we have for any $k\ge0,$
\begin{equation}\label{Ai1esti}
|\nabla_{x'}^{k}A_{i1}^{1}|,|\nabla_{x'}^kA_{i1}^{2}|\le C\delta(x')^{l-
i-\frac{k}{2}-\frac52},\quad|\nabla_{x'}^kB_{i1}|\le C \delta(x')^{l-i-
\frac{k}{2}-2}\quad\text{in}~\Omega_{2R}.
\end{equation}
It follows from \eqref{f3l-1-nambatilp2} and \eqref{Ai1esti} that
\begin{equation*}
|({\bf f}_3^{l-1})^{(j)}-\partial_{x_{j}}\tilde{p}_l|\leq C\delta(x')^{l-
\frac52},~j=1,2,\quad|({\bf f}_3^{l-1})^{(3)}-\partial_{x_{3}}\tilde{p}_l|
\leq C\delta(x')^{l-2}.
\end{equation*}

To further reduce the upper bound, $\delta(x')^{l-\frac52}$, we choose a divergence-free auxiliary function ${\bf v}_{3}^{l}(x)$ that satisfies ${\bf v}_{3}^{l}(x)=0$ on $\Gamma^+_{2R}\cup\Gamma^-_{2R}$, specifically,
\begin{equation}\label{v3lform}
\begin{split}
&({\bf v}_{3}^{l})^{(j)}=\Big(\sum_{i=0}^{2l}F_{jl}^{i}(x')x_{3}^{i}
+\tilde{F}_{jl}(x')\Big)\big(k(x)^2-\frac14\big),~j=1,2,\\
&({\bf v}_{3}^{l})^{(3)}=\Big(\sum_{i=0}^{2l+1}F_{3l}^{i}(x')x_{3}^{i}+\tilde{F}
_{3l}(x)\Big)\big(k(x)^2-\frac14\big)
\end{split}\quad\text{in}~\Omega_{2R}.
\end{equation}
Here we take $F_{1l}^i(x')$ and $F_{2l}^i(x')$, $i=0,1,\dots,2l$, so that
\begin{equation}\label{F12delef3l-1-nabtilpl}
\mu\partial_{x_3x_3}\Big(\sum_{i=0}^{2l}F_{jl}^{i}(x')x_{3}^{i}\Big(k(x)^2-
\frac{1}{4}\Big)\Big)+ ({\bf f}_{3}^{l-1})^{(j)}-\partial_{x_j}\tilde{p}_{l}
=0,\,\, j=1,2.
\end{equation}
By comparing the coefficients of different order of $x_3$, we have
\begin{equation}\label{Fjlform}
\begin{aligned}
F_{jl}^{i}(x')=&\,-\frac{\delta(x')^2A_{i(l-1)}^{j}(x')}{\mu(i+1)(i+2)}+(h_1-
h_2)(x')F_{jl}^{i+1}(x') +\frac14(\varepsilon+2h_1)(\varepsilon+2h_2)F_{jl}^{i+2}(x').
\end{aligned}
\end{equation}
Here we used the convention that $F_{1l}^{i}, F_{2l}^{i}=0$ if $i\notin\{0,\ldots,
2l\}$ and $F_{3l}^{i}=0$ if $i\notin\{0,\ldots,2l+1\}.$
To ensure $\nabla\cdot{\bf v}_{3}^{2}=0,$ we set
\begin{align}\label{F32iform}
F_{3l}^{i}(x')=&\sum_{j=1}^{2}\frac{\delta(x')^2}{i+2}\Big(\partial_{x_j}
\Big(\frac{(h_1-h_2)F_{jl}^{i}(x')-F_{jl}^{i-1}(x')}{\delta(x')^2}+\frac{(\varepsilon+2h_1)
(\varepsilon+2h_2)}{4\delta(x')^2}F_{jl}^{i+1}(x')\Big)\Big)\nonumber\\
&+(h_1-h_2)(x')F_{3l}^{i+1}(x')+\frac14(\varepsilon+2h_1(x'))
(\varepsilon+2h_2(x'))F_{3l}^{i+2}(x'),
\end{align}
such that
\begin{equation*}
\nabla\cdot
\begin{pmatrix}
F_{1l}\big(k(x)^2-\frac14\big) \\
F_{2l}\big(k(x)^2-\frac14\big) \\
F_{3l}\big(k(x)^2-\frac14\big)
\end{pmatrix} :=\nabla \cdot\begin{pmatrix}
\Big(\sum_{i=0}^{2l}F_{1l}^{i}(x')x_{3}
^{i}\Big)\big(k(x)^2-\frac14\big)\\
\Big(\sum_{i=0}^{2l}F_{2l}^{i}(x')x_{3}
^{i}\Big)\big(k(x)^2-\frac14\big)\\
\Big(\sum_{i=0}^{2l+1}F_{3l}^{i}(x')x_{3}
^{i}\Big)\big(k(x)^2-\frac14\big)
\end{pmatrix}
 :=R(x'),
\end{equation*}
where
\begin{align*}
R(x')=&\,-\frac{h_1-h_2}{\delta(x')^2}F_{3l}^{0}(x')-
\frac{(\varepsilon+2h_1(x'))(\varepsilon+2h_2(x'))}{4\delta(x')^2}F_{3l}^{1}(x')
\\
&\, -\sum_{j=1}^{2}\partial_{x_j}\Big(\frac{(\varepsilon+2h_1(x'))
(\varepsilon+2h_2(x'))}{4\delta(x')^2}F_{jl}^{0}(x')\Big).
\end{align*} 
By \eqref{Ai1esti}, \eqref{Fjlform}, and \eqref{F32iform}, a calculation gives, 
for $k\ge0,$ $0\le i\le 2l$, and $0\le j\le 2l+1,$
\begin{equation}\label{F1222est}
|\nabla_{x'}^{k}F_{1l}^{i}|,|\nabla_{x'}^{k}F_{2l}^{i}|\le 
C\delta(x')^{l-i-\frac{k}{2}-\frac12},\,\,|\nabla_{x'}^kF_{3l}^{i}|\le 
C\delta(x')^{l-i-\frac{k}{2}},\,|\nabla_{x'}^{k}R|\le C\delta(x')^{l-1-\frac{k}
{2}}.
\end{equation}

As in Step I.2, $R(x')$ is radial. In addition, we can choose $\tilde{F}_{1l}$, 
$\tilde{F}_{2l}$, $\tilde{F}_{3l}$, and $\hat{p}_{2}$ such that 
\begin{align*}
\nabla\cdot
\begin{pmatrix}
\tilde{F}_{1l}(x')\big(k(x)^2-\frac14\big) \\
\tilde{F}_{2l}(x')\big(k(x)^2-\frac14\big)\\
\tilde{F}_{3l}(x)\big(k(x)^2-\frac14\big)
\end{pmatrix}=-R(x'),
\quad
\partial_{x_1}\hat{p}_{l}=\frac{2\mu\tilde{F}_{1l}(x')}{\delta(x')^2},\quad
\mbox{and}~~\partial_{x_2}\hat{p}_{l}=\frac{2\mu\tilde{F}_{2l}(x')}{\delta(x')^2}.
\end{align*}

As in \eqref{tilF11divfree}--\eqref{tilF11}, we derive
\begin{equation}\label{hatp2tilF12}
\begin{aligned}
\hat{p}_{l}(x')&=-12\mu\int_{r}^{R}\frac{1}{s\delta(s)^3}\int_{0}^{s}
R(t)t\delta(t)\, dt \ ds,\\
\tilde{F}_{il}(x')&=\frac{\delta(x')^2}{2\mu}\partial_{x_i}\hat{p}_{l}
=\frac{6x_{i}}{r^2\delta(r)}\int_{0}^{r}R(t)t\delta(t)\, dt,\quad i=1,2,
\end{aligned}
\end{equation}
and 
\begin{align}\label{tilF32form1}
\tilde{F}_{3l}(x)
=-2\delta(x')k(x)R(x')-\frac{6\partial_{r}k(x)}{r}\int_{0}^{r}
R(t)t\delta(t)\, dt:=P(x')x_3+Q(x').
\end{align}
By \eqref{F1222est}, \eqref{hatp2tilF12}, and \eqref{tilF32form1}, 
a calculation gives, for any $k\ge0,$
\begin{align}\label{tilF12esti}
|\nabla_{x'}^{k}\tilde{F}_{1l}|,|\nabla_{x'}^{k}\tilde{F}_{2l}|\le 
C\delta(x')^{l-\frac{k+1}{2}},\,\,|\nabla_{x'}^kP|\le C\delta(x')^{l-\frac{k+2}
{2}},\,\,|\nabla_{x'}^kQ|\le C\delta(x')^{l-\frac{k}{2}}.
\end{align}
Now we have $\nabla\cdot{\bf v}_{3}^{l}=0$.

{\bf Step II.2. Estimates of ${\bf f}_{3}^{l}.$ }
Let $\bar{p}_{l}=\tilde{p}_{l}+\hat{p}_{l}.$ By \eqref{def3j}, \eqref{f3l-1-nambatilp2}, 
\eqref{F12delef3l-1-nabtilpl}, and \eqref{hatp2tilF12}, we have 
\begin{equation*}
 {\bf f}_{3}^{l}={\bf f}_{3}^{l-1}+\mu\Delta{\bf v}_{3}^{l}-\nabla\bar{p}_{l}=
\begin{pmatrix}
\mu\Delta_{x'}({\bf v}_{3}^{l})^{(1)}(x)\\\\
 \mu\Delta_{x'}({\bf v}_{3}^{l})^{(2)}(x)\\\\
\mu \Delta_{x'}({\bf v}_{3}^{l})^{(3)}(x)+\mu\Delta_{x'}
({\bf v}_{3}^{l-1})^{(3)}(x)+\mu\partial_{x_3x_3}({\bf v}_{3}^{l})^{(3)}(x)
\end{pmatrix}.
\end{equation*}
From \eqref{v3lform}, it follows that $({\bf v}_3^l)^{(i)},i=1,2$ is a 
polynomial in $x_3$ of order $2l+2$, while $({\bf v}_3^1)^{(3)}$ is a 
polynomial in $x_3$ of order $2l+3$.

In view of \eqref{f3l-1-nambatilp2}, we can rewrite 
\begin{equation}\label{f3l1,f3l2}
({\bf f}_{3}^{l})^{(j)}=\sum_{i=0}^{2l+2}H^{j}_{il}(x')x_{3}^{i},\, j=1,2,
\quad({\bf f}_{3}^{l})^{(3)}=\sum_{i=0}^{2l+1}S_{il}(x')x_{3}^{i}+
\sum_{i=0}^{2l+3}G_{il}(x')x_{3}^{i}\quad \text{in}~\Omega_{2R},
\end{equation}
where $
\sum_{i=0}^{2l+1}S_{il}(x')x_3^{i}:= \mu\Delta_{x'}(\hat{\bf v}
_{3}^{l-1})^{(3)}+\mu\partial_{x_3x_3}({\bf v}_{3}^{l})^{(3)}$ and $
\sum_{i=0}^{5}G_{il}(x')x_{3}^{i}:=\mu\Delta_{x'}
({\bf v}_{3}^{l})^{(3)}.$
By using \eqref{Ai1esti}, \eqref{F1222est}, and \eqref{tilF12esti}, for 
$k\ge0,$ in $\Omega_{2R},$ 
\begin{equation}\label{GH2esti}
|\nabla_{x'}^{k}H_{il}^{1}|, |\nabla_{x'}^{k}H_{il}^{2}| \le 
C\delta(x')^{l-i-\frac{k}{2}-\frac32},\,\,|\nabla_{x'}^{k}S_{il}|\le C 
\delta(x')^{l-i-\frac{k}{2}-2},\,\, |\nabla_{x'}^kG_{il}|\le C \delta(x')^{l-i-
\frac{k}{2}-1}.
\end{equation}
Thus, \eqref{f3j-form} and \eqref{HSG-esti} hold for any $j\ge1.$

By \eqref{v3lform}, \eqref{F1222est}, and \eqref{tilF12esti}, we have, for 
$i=1,2$ and $k\ge0$, in $\Omega_{2R}$,
\begin{equation}\label{v3lest}
\begin{split}
&|\nabla_{x'}^{k}\partial_{x_3}^{s}({\bf v}_{3}^{l})^{(i)}|\le C \delta(x')^{l-
s-\frac{k+1}{2}}\quad\text{for}~s\le2l+2,\quad\partial_{x_3}^{s}
({\bf v}_{3}^{l})^{(i)}=0\quad\text{for}~s\ge2l+3, \\&|\nabla_{x'}^{k}\partial_{x_3}^{s}({\bf v}_{3}^{l})^{(3)}|\le C 
\delta(x')^{l-s-\frac{k}{2}}\quad\text{for}~s\le2l+3,
\quad\partial_{x_3}^{s}({\bf v}_{3}^{l})^{(3)}=0\quad\text{for}~s\ge2l+4,\\
&|\nabla_{x'}^{k}\partial_{x_3}^{s}\tilde{p}_{l}|\le C\delta(x')^{l-
s-\frac{k+4}{2}}\, \text{for}~s\le2l,\, \partial_{x_3}^{s}\tilde{p}_{l}
=0\,\,\,\text{for}~s>2l, \,|\nabla_{x'}^{k}\hat{p}_{l}|\le 
C\delta(x')^{l-\frac{k+4}{2}}.
\end{split}
\end{equation}
From \eqref{f3l1,f3l2} and \eqref{GH2esti}, we have
\begin{equation*}
|({\bf f}_{3}^{m+1})^{(1)}|,|({\bf f}_{3}^{m+1})^{(2)}| \le C\delta(x')^{m-
\frac12}, \quad|({\bf f}_{3}^{m+1})^{(3)}|\le C\delta(x')^{m-1} \quad \text{in}
~\Omega_{2R}.
\end{equation*}
Thus,
\begin{equation*}
|{\bf f}_{3}^{m+1}|\le C (|({\bf f}_{3}^{m+1})^{(1)}|+|({\bf f}_{3}
^{m+1})^{(2)}| +|({\bf f}_{3}^{m+1})^{(3)}|) \le C\delta(x')^{m-1} 
\quad\text{in}~\Omega_{2R}.
\end{equation*}
Similarly, for $1\le s\le m+1$,
\begin{equation*}
|\nabla^{s}{\bf f}_{3}^{m+1}|\le C \delta(x')^{m-1-s}\quad\text{in}~\Omega_{2R}.
\end{equation*}

{\bf Step II.3 Estimates of $\nabla^{m+1}{\bf u}_{1}^{3}$ and $\nabla^{m}p_{1}^3$.}
Denote
\begin{equation*}
{\bf v}^{m+1}(x)=\sum_{l=1}^{m+1}{\bf v}_{3}^{l}(x),\quad \bar{p}^{m+1}
(x)=\sum_{l=1}^{m+1}\bar{p}_{l}(x)\quad\text{in}~\Omega_{2R}.
\end{equation*}
Then we have 
${\bf v}^{m+1}(x)={\bf u}_1^3(x)$ on $\Gamma_{2R}^{\pm}$, $\nabla\cdot{\bf v}
^{m+1}(x)=0$ in $\Omega_{2R}$, and for $1\le s\le m$,
\begin{equation*}
|{\bf f}^{m+1}(x)|=|\mu\Delta{\bf v}^{m+1}(x)-\nabla \bar{p}^{m+1}(x)|\leq C 
\delta(x')^{m-1},\quad |\nabla^s{\bf f}^{m+1}(x)|\leq C \delta(x')^{m-s-1}.
\end{equation*}
Thus, by Proposition \ref{keyprop2}, it holds that
\begin{equation}\label{u13minusest}
\|\nabla^{m+1}({\bf u}_{1}^{3} -{\bf v}^{m+1})\|_{L^{\infty}(\Omega_{\delta(x')/
2}(x'))} +\|\nabla^{m} (p_{1}^{3}-\bar{p}^{m+1})\|_{L^{\infty}
(\Omega_{\delta(x')/2}(x'))} \le C.
\end{equation}
Moreover, by \eqref{v3lest}, we have, for $2\le l\le m+1,$
\begin{equation}\label{gradv3lest}
\begin{split}
&|\nabla^{m+1}{\bf v}_{3}^{l}|\le C\delta(x')^{-\frac{m+4}{2}}\,\,\text{for}
~l\le \frac{m-2}{2},\quad
|\nabla^{m+1}{\bf v}_{3}^{l}|\le C\delta(x')^{l-m-\frac{3}{2}}\quad\text{for}
~\frac{m-1}{2}\le l,\\
&|\nabla^{m}\bar{p}_{l}|\le C \delta(x')^{-\frac{m+4}{2}}\quad\,\,\,\text{for}~ 
l\le \frac{m-1}{2},\quad |\nabla^{m}\bar{p}_{l}|\le C\delta(x')^{l-m-2}
\quad\text{for}~ \frac{m}{2}\le l.
\end{split}
\end{equation} 
It follows from \eqref{u13minusest} and \eqref{gradv3lest} that for $m\ge1$,
\begin{equation*}
|\nabla^{m+1} {\bf u}_{1}^{3}| +|\nabla^{m} {p}_{1}^{3}|\le C \delta(x')^{-
\frac{m+4}{2}}\quad\text{in}~\Omega_{R}.
\end{equation*}
For the estimates of $p_{1}^{3}$, for a fixed point $z=(z'_1,0)\in\Omega_{R}$ 
with $|z'_1|=R$, by the mean value theorem, \eqref{u13minusest}, and 
\eqref{gradv3lest}, we have
\begin{equation*}
|p_{1}^{3}(x)-p_{1}^{3}(z'_1,0)|\le C\|\nabla(p_1^3-\bar{p}_1-\bar{p}_2)\|
_{L^\infty}+C|\bar{p}_1(x)+\bar{p}_2(x)|+C\leq C\delta(x')^{-2}.
\end{equation*}
This finishes the proof of Proposition \ref{u13-esti}.

\end{proof}

\section{Estimates for \texorpdfstring{$({\bf u}_{1}^{\alpha},p_{1}^{\alpha}), \,\,\alpha=4,5,6$}{}}\label{sec5}
In this section, we will derive the estimates of ${\bf u}_{1}^{\alpha}$ and $p_{1}^{\alpha}$, $\alpha=4,5,6$, by
using the same method as in the proof of Proposition \ref{u11-general}.
\begin{prop}\label{u14-esti}
Under the same assumption as in Theorem \ref{main thm1}, let ${\bf u}_{1}^{4}$ 
and $p_{1}^{4}$ be the solution to \eqref{u,peq1} with $\alpha=4$. 
Then for sufficiently small $0<\varepsilon<1/2$, we have 
\begin{equation*}
\mbox{(i)}~|\nabla{\bf u}_{1}^{4}(x)|\le C(\varepsilon+|x'|^2)^{-
\frac12},\quad|p_1^{4}(x)-p_1^{4}(z',0)|\leq C \quad\text{in}
~\Omega_{R}\qquad\quad\qquad
\end{equation*}
for some point $z=(z',0)\in\Omega_{R}$, $|z'|=R/2,$ and

(ii) for any $m\ge 1$, 
\begin{equation*}
|\nabla^{m+1}{\bf u}_{1}^{4}(x)|+|\nabla^{m}p_{1}^{4}(x)|\le C (\varepsilon+|x'|
^2)^{-\frac{m+1}{2}}\quad \text{for}~x\in\Omega_{R}.
\end{equation*}
\end{prop}
\begin{proof}[Proof of Proposition \ref{u14-esti}]
 The proof essentially follows that of Proposition \ref{u11-general} with some modifications. We only list the key steps here.
 
{\bf Step I. Gradient estimates.} 
We first choose 
\begin{equation}\label{v41form}
{\bf v}_{4}^{1}(x)=
\begin{pmatrix}
x_2\\-x_1\\0
\end{pmatrix}\Big(k(x)+\frac{1}{2}\Big) \quad\text{in}~\Omega_{2R},
\end{equation}
which has the same boundary conditions as ${\bf u}_1^4$ on $\Gamma_{2R}^{\pm}.
$ Since $k(x)$ is radial in $x',$ we have
\begin{equation}\label{x2par1keqx1par2k}
x_2\partial_{x_1}k(x)=\frac{x_1x_2}{r}\partial_{r}k(x)=x_1\partial_{x_2}k(x).
\end{equation}
By \eqref{h1h2prop} and \eqref{x2par1keqx1par2k}, a simple calculation yields
\begin{equation}\label{Deltav41}
\nabla\cdot{\bf v}_{4}^{1}=0,\quad
\Delta({\bf v}_{4}^{1})^{(1)}=\Delta_{x'}k(x)x_{2}+2\partial_{x_{2}}k(x),
\quad\Delta({\bf v}_{4}^{1})^{(2)}=-\Delta_{x'}k(x)x_{1}-2\partial_{x_{1}}
k(x),
\end{equation}
and for any $i=1,2,k\ge0$,
\begin{equation*}
|\nabla_{x'}^{k}\partial_{x_{3}}^s({\bf v}_{4}^{1})^{(i)}|\leq 
C\delta(x')^{\frac{1-k}{2}-s},~s\le 1,\quad \partial_{x_{3}}^2({\bf v}_{4}
^{1})^{(i)}=0, \quad \nabla^k ({\bf v}_{4}^{1})^{(3)}=0.
\end{equation*}
Hence, we have for any $m\ge0$
\begin{equation}\label{gradv41esti}
|\nabla^{m+1}{\bf v}_{4}^{1}|\le C\delta(x')^{-\frac{m+1}{2}}.
\end{equation} 
Denote ${\bf f}_4^1:=\mu\Delta{\bf v}_{4}^{1}-\nabla\bar{p}_1$ and take $\bar{p}
_1(x)=0$. By using \eqref{h1h2prop}, \eqref{def_k}, and \eqref{Deltav41}, 
 ${\bf f}_4^1$ can be rewritten as
\begin{equation}\label{f41}
({\bf f}_4^1)^{(b)}=S_{11}^b(x')x_3+S_{01}^b(x'), \,b=1,2,\quad({\bf f}
_4^1)^{(3)}=0,
\end{equation}
and for any $k\ge0$,
\begin{equation}\label{f41xsesti}
|\nabla_{x'}^kS_{11}^b|\leq C\delta(x')^{-\frac{k+3}{2}},\quad |\nabla_{x'}
^kS_{01}^b|\leq C\delta(x')^{-\frac{k+1}{2}}.
\end{equation}
By \eqref{v41form}, one can see that
\begin{equation}\label{xj1}
 S_{i1}^1(x')=x_2 g_{i1}(r),\quad S_{i1}^2(x')=-x_1 g_{i1}(r),~i=0,1,
\end{equation}
where $ g_{11}(r)=\mu\partial_{rr}\Big(\frac{1}{\delta(x')}\Big)+\frac{3\mu}{r}\partial_{r}\Big(\frac{1}{\delta(x')}\Big)$ and $ g_{01}(r)=\mu\partial_{rr}\Big(\frac{\varepsilon+2h_2(x')}{2\delta(x')}\Big)+\frac{3\mu}{r}\partial_{r}\Big(\frac{\varepsilon+2h_2(x')}{2\delta(x')}\Big)$. 

It follows from \eqref{f41} and \eqref{f41xsesti} that
\begin{equation}\label{f41esti}
|{\bf f}_{4}^{1}|\le C (|({\bf f}_{4}^{1})^{(1)}|+|({\bf f}_{4}^{1})^{(2)}|) \le 
C \delta(x')^{-\frac12}\quad \text{in}~\Omega_{2R}.
\end{equation}
By applying Proposition \ref{keyprop2} and using \eqref{f41esti}, we have for 
any 
$x\in \Omega_{R}$,
\begin{equation}\label{u14-v41est}
\| \nabla {\bf u}_{1}^{4}-\nabla{\bf v}_{4}^1\|_{L^{\infty}
(\Omega_{\delta(x')/2}(x'))}\le C\delta(x')^{\frac12}.
\end{equation}
It follows from \eqref{gradv41esti} and \eqref{u14-v41est} that 
$$|\nabla {\bf u}_{1}^{4}|\le C|\nabla {\bf v}^{1}_{4}(x)| 
+C\delta(x')^{\frac12}\le C\delta(x')^{-\frac12}\quad\text{in}~\Omega_{R}.$$
Thus, Proposition \ref{u14-esti} holds for $m=0.$

{\bf Step II. Higher-order derivatives estimates.} 
Denote
\begin{equation*}
{\bf f}_{4}^{j}:={\bf f}_{4}^{j-1}(x)+\mu\Delta{\bf v}_{4}^{j}(x)-\nabla\bar{p}_{j}(x)=\sum_{l=1}^{j}(\mu\Delta{\bf v}_{4}^{l}-\nabla\bar{p}_{l})(x)
\quad\text{for}~ 2\le j\le m+1.
\end{equation*}
In order to use Proposition \ref{keyprop2} to estimate $\nabla^{m+1}{\bf u}_1^4$ and $\nabla \bar{p}_{m}^{4}$ for $m\ge1$, we need to construct a series of auxiliary functions ${\bf v}_{4}^{l}$ and $\bar{p}_{l}$ to successively reduce the upper bounds of the non-homogeneous term. 

To cancel the leading terms 
introduced by $({\bf f}_4^1)^{(b)}$, $b=1,2$, we first choose ${\bf v}_4^2$ 
satisfying 
\begin{equation}\label{v41-12form}
({\bf v}_4^2)^{(b)}= \Big(F_{b2}^1(x')x_3+F_{b2}^0(x')\Big)\Big(k(x)^2-
\frac14\Big),~b=1,2,
\end{equation}
where 
\begin{equation}\label{4F1222form}
F_{b2}^{0}(x')=-\frac{\delta(x')^2}{2\mu}S_{01}^{b}(x')+(h_1-h_{2})F_{b2}^{1}
(x'),\quad F_{b2}^{1}(x')=-\frac{\delta(x')^2}{6\mu}S_{11}^{b}(x'),
\end{equation}
such that 
\begin{equation}\label{cancelf421}
\mu\partial_{x_3x_3}({\bf v}_{4}^2)^{(b)}=-({\bf f}_{4}^{1})^{(b)},\quad b=1,2.
\end{equation}
We also take $({\bf v}_{4}^{2})^{(3)}=0.$ Then, substituting \eqref{xj1} into \eqref{4F1222form}, yields
\begin{equation}\label{xj2}
\begin{aligned}
&F_{12}^1(x')=x_2f_{2}^1(r), ~ F_{22}^1(x')=-x_1f_{2}^1(r),~ 
F_{12}^{0}(x')=x_2f_{2}^{0}(r),~ F_{22}^{0}(x')=-x_1f_{2}^{0}(r),
\end{aligned} 
\end{equation}
where $f_{2}^{1}(r)=-\frac{\delta(x')^2}{2\mu}g_{01}(r)-\frac{\delta(x')^2}
{6\mu}(h_1-h_2)g_{11}(r)$ and $f_{2}^{0}(r)=-\frac{\delta(x')^2}{6\mu}g_{11}(r)
$ are radial. 
By \eqref{v41-12form} and \eqref{xj2}, a direct calculation gives $
\partial_{x_1}({\bf v}_{4}^{2})^{(1)}+\partial_{x_2}({\bf v}_{4}
^{2})^{(2)}=0.$ Thus $\nabla\cdot{\bf v}_{4}^{2}=0.$ 

In addition, by \eqref{f41xsesti} and \eqref{4F1222form}, a calculation gives, for 
any $k\ge0,$
\begin{equation}\label{v42xsesti}
|\nabla_{x'}^{k}F_{b2}^{1}|\le C\delta(x')^{\frac12-\frac{k}{2}},\quad|
\nabla_{x'}^{k}F_{b2}^{0}|\le C\delta(x')^{\frac32-\frac{k}{2}}.
\end{equation}
Hence, by \eqref{h1h2prop}, \eqref{v41-12form}, and \eqref{v42xsesti}, we get for $i=1,2$ and $k\ge0$,
\begin{equation*}
|\nabla_{x'}^{k}\partial_{x_{3}}^s({\bf v}_{4}^{2})^{(i)}|\leq 
C\delta(x')^{\frac{3-k}{2}-s},~s\le3 ,\quad \partial_{x_{3}}^4({\bf v}_{4}
^{2})^{(i)}=0,\quad\nabla^k ({\bf v}_{4}^{2})^{(3)}=0.
\end{equation*}
Therefore, for any $m\ge0,$
\begin{equation*}
|\nabla^{m+1}{\bf v}_{4}^{2}|\le C\delta(x')^{\frac12-m}\quad\text{for}
~0\le m\le2,\quad |\nabla^{m+1}{\bf v}_{4}^{2}|\le C\delta(x')^{-\frac{m+1}
{2}}\quad\text{for}~m\ge3.
\end{equation*}
Let $\bar{p}_{2}=0.$ By \eqref{v41-12form}, \eqref{cancelf421}, and 
\eqref{v42xsesti}, we can write ${\bf f}_{4}^{2}$ as polynomials in $x_3$ 
\begin{equation}\label{f42form}
({\bf f}_{4}^{2})^{(b)}=\mu\Delta_{x'}({\bf v}_{4}^{2})^{(b)}=\sum_{i=0}^{3}
S_{i2}^b(x') x_{3}^{i},\quad({\bf f}_{4}^{2})^{(3)}=0,
\end{equation}
with the estimates for $k\ge0$,
\begin{equation}\label{S1222esti}
|\nabla_{x'}^{k}S^{b}_{i2}|\le C\delta(x')^{\frac12-i-\frac{k}{2}}.
\end{equation}
Thus, we have
\begin{equation*}
|({\bf f}_{4}^{2})^{(b)}|=|\mu\Delta_{x'}({\bf v}_{4}^{2})^{(b)}|\le 
C\delta(x')^{\frac12},\quad({\bf f}_{4}^{2})^{(3)}=0\quad\text{in}~\Omega_{2R},
\end{equation*}
and
\begin{equation*}
|{\bf f}_{4}^{2}|\le |({\bf f}_{4}^{2})^{(1)}| +|({\bf f}_{4}^{2})^{(2)}|+|({\bf 
f}_{4}^{2})^{(3)}|\le C\delta(x')^{\frac12}.
\end{equation*}
In fact, by \eqref{v41-12form} and \eqref{xj2}, one can write
\begin{equation}\label{xj3}
S^{1}_{i2}(x')=x_2g_{i2}(r),\quad S^{2}_{i2}(x')=-x_1g_{i2}(r),~i=0,1,2,3,
\end{equation}
where $g_{i2}(x')$ are radial functions. 

Now we will inductively prove that for $2\le j\le m+1,$ we can successively 
construct a sequence of divergence-free ${\bf v}_{4}^{j}$ and $\bar{p}_{j}=0$ to 
make ${\bf f}$ in \eqref{weqs1} as small as possible. Besides, it can be 
inductively proved that ${\bf f}_{4}^{j}$ can be represented as polynomials in 
$x_{3}$ as follows:
\begin{equation}\label{f2jform}
({\bf f}_{4}^{j})^{(b)}=\mu\Delta_{x'}({\bf v}_{4}^{j})^{(b)}=\sum_{i=0}^{2j-1}
S^{b}_{ij}(x')x_3^{i},~b=1,2,\quad({\bf f}_{4}^{j})^{(3)}=0\quad\text{in}
~\Omega_{2R},
\end{equation} 
where for any $k\ge0$, $b=1,2,$
 \begin{equation}\label{Sijesti}
 |\nabla_{x'}^{k}S^{b}_{ij}|\le C \delta(x')^{j-i-\frac{k}{2}-\frac32},\quad
 \nabla\cdot{\bf v}_{4}^{j}=0.
 \end{equation}
Indeed, by \eqref{f42form} and \eqref{S1222esti}, we have \eqref{f2jform} and 
\eqref{Sijesti} hold for $j=2.$ Assume that \eqref{f2jform} and \eqref{Sijesti} 
hold for $j=l-1$ with $l\ge3,$ that is, in $\Omega_{2R},$
\begin{equation}\label{f4l-1form}
({\bf f}_{4}^{l-1})^{(b)}=\mu\Delta_{x'}({\bf v}_{4}^{l-1})^{(b)}=\sum_{i=0}
^{2l-3}S^{b}_{i(l-1)}(x')x_3^{i},~b=1,2,\quad({\bf f}_{4}^{l-1})^{(3)}
=0\quad\text{in}~\Omega_{2R},
\end{equation}
where for any $k\ge0$ and $b=1,2,$
\begin{equation}\label{Sil-1form}
|\nabla_{x'}^{k}S^{b}_{i(l-1)}|\le C \delta(x')^{l-i-\frac{k}{2}-\frac52},\quad
\nabla\cdot{\bf v}_{4}^{l-1}=0.
\end{equation}
Then we will construct ${\bf v}_{4}^{l}$ and $\bar{p}_{l}=0$ so that ${\bf f}
_{4}^{l}$ satisfies \eqref{f2jform} and \eqref{Sijesti}. 

It follows from \eqref{f4l-1form} and \eqref{Sil-1form} that 
$$|({\bf f}_{4}^{l-1})^{(b)}|\le C\delta(x')^{l-\frac52},~b=1,2,\quad({\bf f}_{4}
^{l-1})^{(3)}=0\quad\text{in}~\Omega_{2R}.$$
To reduce the upper bound, $\delta(x')^{l-\frac52}$, we choose a divergence-free 
auxiliary function ${\bf v}_{4}^{l}$ satisfying ${\bf v}_{4}^{l}(x)=0$ on $
\Gamma_{2R}^{+}\cup\Gamma_{2R}^{-},$ specifically,
\begin{equation}\label{v4l-12form}
({\bf v}_{4}^{l})^{(b)}=\sum_{i=0}^{2l-3}F_{bl}^{i}(x')x_3^{i}\Big(k(x)^2-
\frac14\Big),~b=1,2\quad\text{in}~\Omega_{2R}.
\end{equation}
We take $F_{bl}^{i}(x'),$ $i=0,\ldots,2l-3$, such that 
\begin{equation}\label{par33delef4l-1}
\mu\partial_{x_3x_3}({\bf v}_{4}^{l})^{(b)}=-({\bf f}_{4}^{l-1})^{(b)}.
\end{equation}
By comparing the coefficients of different order of $x_{3}$, we have
\begin{equation}\label{Fblform}
F_{bl}^{i}(x')=-\frac{\delta(x')^2S_{i(l-1)}^{b}(x')}{\mu(i+1)(i+2)}+
(h_1-h_2)(x')F_{bl}^{i+1}(x')+\frac14(\varepsilon+2h_1(x'))
(\varepsilon+2h_2(x'))F_{bl}^{i+2}(x').
\end{equation}
Here we use the convention that $F_{bl}^{i}(x')\equiv0$ if $i\notin\{0,\dots,2l-3\}.$ Moreover, by \eqref{h1h2prop}, \eqref{Sil-1form} and \eqref{Fblform}, we derive 
for $k\ge0,0\le i\le 2l-3,$ in $\Omega_{2R},$
\begin{equation}\label{Fblesti}
|\nabla_{x'}^{k}F_{bl}^{i}|\le C\delta(x')^{l-i-\frac{k}{2}-\frac12}.
\end{equation}
Let $({\bf v}_{4}^{l})^{(3)}=0 $ and $\bar{p}_{l}=0.$ By using \eqref{h1h2prop}, \eqref{v4l-12form}, 
\eqref{par33delef4l-1}, and \eqref{Fblesti}, we have
\begin{equation}\label{f41form}
({\bf f}_{4}^{l})^{(b)}=\mu\Delta_{x'}({\bf v}_{4}^{l})^{(b)}=\sum_{i=0}
^{2l-1}S^{b}_{il}(x')x_3^{i},~b=1,2,\quad({\bf f}_{4}^{l})^{(3)}
=0\quad\text{in}~\Omega_{2R},
\end{equation}
where for $k\ge0$ and $b=1,2$,
\begin{equation}\label{Silesti}
|\nabla_{x'}^{k}S^{b}_{ij}|\le C \delta(x')^{l-i-\frac{k}{2}-\frac32}.
\end{equation}
Therefore, by \eqref{f41form} and \eqref{Silesti}, we have \begin{equation*}
|({\bf f}_{4}^{l})^{(b)}|=|\mu\Delta_{x'}({\bf v}_{4}^{l})^{b}|\le 
C\delta(x')^{l-\frac32},\quad({\bf f}_{4}^{l})^{(3)}=0.
\end{equation*}
Hence \eqref{f2jform} and \eqref{Sijesti} hold for $j=l.$ By using 
\eqref{h1h2prop}, \eqref{v41-12form}, and \eqref{Fblesti}, a calculation gives 
for 
$l\ge2,~b=1,2,$ and $k\ge0,$ 
\begin{equation}\label{v41esti}
|\nabla_{x'}^{k}\partial_{x_3}^{s}({\bf v}_{4}^{l})^{(b)}|\le C\delta(x')^{l-s-
\frac{k+1}{2}}~\text{for}~s\le2l-1,\,\,\partial_{x_3}^s({\bf v}_{4}
^{l})^{(b)}=0~\text{for}~s\ge2l,\,\,\nabla^{k}({\bf v}_{4}^{l})^{(3)}=0.
\end{equation}
From \eqref{f2jform} and \eqref{Sijesti}, we have $|({\bf f}_{1}^{m+1})^{(1)}|,|
({\bf f}_{1}^{m+1})^{(1)}|\le C\delta(x')^{m-\frac12}$, $({\bf f}_{1}
^{m+1})^{(3)}=0$, and then 
$$|{\bf f}_{4}^{m+1}|\le |({\bf f}_{1}^{m+1})^{(1)}|+|({\bf f}_{1}^{m+1})^{(1)}|
\le C\delta(x')^{m-\frac12}.$$ 
Similarly, for $1\le s\le m,$
$|\nabla^s{\bf f}_{4}^{m+1}|\le C\delta(x')^{m-s-\frac12}.$

We now prove that ${\bf v}_{4}^{l}$ constructed above satisfies the divergence-free condition (i.e., $\nabla\cdot{\bf v}_{4}^{l}=0$). To this end, we inductively show that for $j\ge2$, $S_{ij}^{b}$ can be written as 
\begin{equation}\label{xj4}
S_{ij}^{1}(x')=x_2g_{ij}(r),\quad S_{ij}^{2}(x')=-x_1g_{ij}(r),~i=0,\ldots,2j-1,
\end{equation}
where $g_{ij}(x')$ are radial. Indeed, it follows from \eqref{xj3} that 
\eqref{xj4} holds for $j=2.$ Assume \eqref{xj4} holds true for $j=l-1$ with 
$l\ge3,$ that is, in $\Omega_{2R},$
 \begin{equation}\label{xj5}
 S_{i(l-1)}^{1}(x')=x_2g_{i(l-1)}(r),\quad S_{i(l-1)}^2(x')=-x_1g_{i(l-1)}
 (r),~i=0,\ldots,2l-3.
 \end{equation}
Substituting \eqref{xj5} into \eqref{Fblform}, we have
\begin{equation}\label{xj6}
F_{1l}^{i}(x')=x_2f_{l}^{i}(r),\,\,F_{2l}^{i}(x')=-x_1f_{l}^{i}(r),~i=0,\ldots,
2l-3,
\end{equation} 
for radial functions $f_{l}^{i}(r)$. By \eqref{v4l-12form} and \eqref{xj6}, 
after a simple calculation, we have \eqref{xj4} holds for any $j\ge2$. By using \eqref{v4l-12form} and \eqref{xj6} again, we have $\partial_{x_{1}}({\bf v}_{4}^{l})^{(1)}+
\partial_{x_2}({\bf v}_{4}^{l})^{(2)}=0.$ Thus $\nabla\cdot{\bf v}_{4}
^{l}=0.$ Denote 
${\bf v}^{m+1} :=\sum_{l=0}^{m+1}{\bf v}_{4}^{l}(x),~ \bar{p}^{m+1}
:=\sum_{l=0}^{m+1}\bar{p}_{l}=0. $
It is easy to verify that ${\bf v}^{m+1}(x)={\bf u}_{1}^{4}(x)=0$ on $
\Gamma_{2R}^{\pm},\nabla\cdot{\bf v}^{m+1}(x)=0$ in $\Omega_{2R},$ and
\begin{align*}
|{\bf f}^{m+1}(x)|=|\mu\Delta{\bf v}^{m+1}(x)|\le C\delta(x')^{m-\frac12},
 \quad
|\nabla^{s}{\bf f}^{m+1}(x)|\le C\delta(x')^{m-s-\frac12},\,\,1\le s\le m.
\end{align*}
Thus, by applying Proposition \ref{keyprop2}, there holds
\begin{equation}\label{u14m-v4mest}
\| \nabla^{m+1} ({\bf u}_{1}^{4}-{\bf v}_{4}^1)\|_{L^{\infty}
(\Omega_{\delta(x')/2}(x'))} +\| \nabla^m (p_{1}^{4}- \bar{p}^{m+1})\|
_{L^{\infty}(\Omega_{\delta(x')/2}(x'))} \le C\delta(x')^{\frac12}.
\end{equation}
By \eqref{v41esti}, we obtain for $l\ge 2,$
\begin{align}\label{nabmv4l}
&|\nabla^{m+1}{\bf v}_{4}^{l}|\le C\delta(x')^{-\frac{m+1}{2}}\quad\text{for}
~l\le \frac{m+1}{2},\nonumber\\
&|\nabla^{m+1}{\bf v}_{4}^{l}|\le C\delta(x')^{l-m-\frac32}\quad\text{for}
~\frac{m+2}{2}\le l\le m+1.
\end{align}
It follows from \eqref{gradv41esti} and \eqref{nabmv4l} that for $m\ge 1$,
\begin{equation}\label{namvm+1}
|\nabla^{m+1}{\bf v}^{m+1}|\le C\delta(x')^{-\frac{m+1}{2}}.
\end{equation}
By \eqref{u14m-v4mest} and \eqref{namvm+1}, we have for $x\in\Omega_{R}$,
\begin{equation*}
|\nabla^{m+1}{\bf u}_{1}^4|\le C\delta(x')^{-\frac{m+1}{2}},\quad|
\nabla^mp_{4}^{1}|\le C\delta(x')^{\frac12}.
\end{equation*}
For the estimates of $p_1^4(x)$, by the mean value theorem, \eqref{u14m-v4mest}, 
and the fact $\bar{p}_{l}=0,$ $l\ge1$, we have
\begin{align*}
|p_{1}^{4}(x)-p_{1}^{4}(z',0)|\leq&\,|p_{1}^{4}-\bar{p}_1-\bar{p}_2-
(p_{1}^{4}-\bar{p}_1-\bar{p}_2)(z',0)|+|\bar{p}_1+\bar{p}_2|\\
\leq&\, C\|\nabla (p_{1}^{1}-\bar{p}_1-\bar{p}_2)\|
_{L^\infty(\Omega_{R})}\leq C
\end{align*}
for a fixed point $(z',0)\in\Omega_{R}$ with $|z'|=R/2.$
\end{proof}

\begin{prop}\label{u15-esti}
Under the same assumption as in Theorem \ref{main thm1}, let ${\bf u}_{1}
^{\alpha}$ and $p_{1}^{\alpha}$ be the solution to \eqref{u,peq1} with $
\alpha=5,6$. Then for sufficiently small $0<\varepsilon<1/2$, we have 
\begin{equation*}
 \mbox{(i) }~ |\nabla{\bf u}_{1}^{\alpha}(x)|\le C(\varepsilon+|x'|^2)^{-1},
 \quad |p_1^\alpha(x)-
 p_1^i(z',0)|\leq C(\varepsilon+|x'|^2)^{-\frac32}\quad\mbox{in}~
 \Omega_{R}
\end{equation*}
for some point $z=(z',0)\in \Omega_{R}$ with $|z'|=R/2$; and
 
(ii) for high-order derivatives, we have
\begin{equation*}
 \begin{aligned}
 &|\nabla^{m+1}{\bf u}_{1}^{\alpha}(x)|+|\nabla^{m}{p}_{1}^{\alpha}(x)| \le 
 C(\varepsilon+|x'|^2)^{-\frac{m+3}{2}}\quad\text{for}~m\le6, \\
 &|\nabla^{m+1}{\bf u}_{1}^{\alpha}(x)|+|\nabla^{m}p_{1}^{\alpha}(x)|\leq\,C(\varepsilon+|x'|^2)^{\frac32-m}\quad\text{for}~m\ge7
 \end{aligned}
 \quad\mbox{in}~\Omega_{R}.
\end{equation*}
\end{prop}
\begin{proof}[Proof of Proposition \ref{u15-esti}]
We only take the case $\alpha=5$ for instance, since the other case is the same. 

{\bf Step I. Gradient estimate.} Similar to the proof of Proposition \ref{u11-general}, we first set
 \begin{equation}\label{v51form}
 {\bf v}_{5}^{1}(x)={\boldsymbol{\psi}}_{5}\Big(k(x)+\frac{1}{2}\Big) +
 \begin{pmatrix}
 F(x)\\
 G(x)\\
 H(x)
 \end{pmatrix}\Big(k(x)^{2}-\frac{1}{4}\Big)\quad\text{in}~\Omega_{2R},
 \end{equation}
 which satisfies ${\bf v}_{5}^{1}(x)={\bf u}_1^5(x)$ on $\Gamma_{2R}^+
 \cup\Gamma_{2R}^-$. A calculation gives
 \begin{equation}\label{partv51}
 \begin{aligned}
 &\partial_{x_1}({\bf v}_{5}^{1})^{(1)}=x_{3}\partial_{x_{1}}k(x)+
 \partial_{x_1}F\big(k(x)^2-\frac14\big)+2k(x)\partial_{x_1}k(x)F(x),\\
 &\partial_{x_2}({\bf v}_{5}^{1})^{(2)}=\partial_{x_2}G(x)\Big(k(x)^2-
 \frac14\Big)+2k(x)\partial_{x_2}k(x)G(x), \\
 &\partial_{x_3}({\bf v}_{5}^{1})^{(3)}=-\frac{x_1}{\delta(x')}+
 \partial_{x_3}H\big(k(x)^2-\frac14\big)+H\big(k(x)^2-\frac14\big).
 \end{aligned}
 \end{equation}
To ensure that $\nabla\cdot{\bf v}_{5}^{1}(x)=0$ in $\Omega_{2R}$, we assume 
 \begin{equation}\label{partFGH}
 \frac14(\partial_{x_{1}}F+\partial_{x_{2}}G)(x)=x_2\partial_{x_{3}}k(x)-
 \frac{x_1}{\delta(x')}.
 \end{equation}
By \eqref{v51form}, \eqref{partv51}, and $\nabla\cdot{\bf v}_{5}^{1}=0,$ we 
obtain
\begin{equation}\label{Hform}
H(x)=-\delta(x')\partial_{x_1}k(x)F(x)-\delta(x')\partial_{x_2}k(x)G(x)+2k(x)
\Big(x_{1}-\delta(x')\partial_{x_1}k(x)x_{3}\Big). 
\end{equation}
Subtituting \eqref{Hform} into \eqref{partFGH}, we get
\begin{align}\label{5FGeq}
&\partial_{x_1}F+\frac{\partial_{x_1}\delta(x')}{\delta(x')}F+\partial_{x_2}G+
\frac{\partial_{x_2}\delta(x')}{\delta(x')}G -\delta(x')\partial_{x_1}k(x)
\partial_{x_3}F-\delta(x')\partial_{x_2}k(x)\partial_{x_3}G\nonumber\\
 &\quad\quad\quad\quad\quad\quad\quad\quad=-\frac{6x_{1}}
{\delta(x')}+6\partial_{x_1}k(x)x_{3}-\frac{2k(x)\partial_{x_1}\delta(x')x_{3}}
{\delta(x')}+2k(x)\partial_{x_1}k(x)\delta(x').
\end{align}
To determine $F(x)$ and $G(x)$, we consider 
\begin{equation}\label{F1G1eq}
\partial_{x_1}F_{1}+\frac{\partial_{x_1}\delta(x')}{\delta(x')}F_{1}+
\partial_{x_2}G_{1}+\frac{\partial_{x_2}\delta(x')}{\delta(x')}G_{1} -\delta(x')
\partial_{x_1}k(x)\partial_{x_3}F_{1}-\delta(x')\partial_{x_2}k(x)\partial_{x_3}
G_{1}= -\frac{6x_{1}}{\delta(x')},
\end{equation}
and
\begin{align}\label{F2G2eq}
&\partial_{x_1}F_{2}+\frac{\partial_{x_1}\delta(x')}{\delta(x')}F_{2}+
\partial_{x_2}G_{2}+
\frac{\partial_{x_2}\delta(x')}{\delta(x')}G_{2} -\delta(x')\partial_{x_1}k(x)
\partial_{x_3}F_{2}-\delta(x')\partial_{x_2}k(x)\partial_{x_3}G_{2}\nonumber\\
&= 6\partial_{x_1}k(x)x_{3}-\frac{2k(x)\partial_{x_1}\delta(x')x_{3}}
{\delta(x')}+2k(x)\partial_{x_1}k(x)\delta(x').
\end{align}
For \eqref{F2G2eq}, we take a pair of particular solution
\begin{equation}\label{F2G2form}
F_{2}(x)=-\frac{3x_{3}^{2}}{\delta(x')}-2k(x)x_{3},\quad G_{2}(x)=0.
\end{equation}
For \eqref{F1G1eq}, we assume $F_{1}(x)=F_{1}(x'),$ $
G_{1}(x)=G_{1}(x').$ In order to cancel the leading terms in $\mu\Delta{\bf v}_5^1$, we first choose $\hat{p}_{1}$ such that
\begin{equation}\label{F1F2-hatp1}
\partial_{x_1}\hat{p}_{1}(x')=\frac{2\mu F_{1}(x')}{\delta(x')^2},\quad 
\partial_{x_2}\hat{p}_{1}=\frac{2\mu G_{1}(x')}{\delta(x')^2}.
\end{equation}
Combining \eqref{F1G1eq} and \eqref{F1F2-hatp1}, we have 
\begin{equation*}
\Delta_{x'} \hat{p}_{1}(x')+\frac{3}{\delta(x')}\nabla_{x'} \delta(x')
\nabla_{x'}\hat{p}_{1}(x')=-\frac{12\mu x_1}{\delta(x')^3}:=\mathcal{F}_{1}(x').
\end{equation*}
Notice that $\big|\mathcal{F}_{1}(x')\big|\le\delta(x')^{-\frac52}$ and $|\nabla_{x'}^s\mathcal{F}_{1}(x')|\le\delta(x')^{-\frac{s+5}{2}}$ for $1\le s\le m.$ By Proposition \ref{U_esti} with $\gamma=-\frac52$, we obtain for $k\ge0$,
\begin{equation}\label{5hatp1esti}
|\hat{p}_{1}|\le C \delta(x')^{-\frac32},\quad|\nabla_{x'}^{k}\hat{p}_{1}|\le 
C\delta(x')^{-\frac{k+3}{2}},\quad \partial_{x_3}\hat{p}_{1}=0.
\end{equation}
Then by \eqref{h1h2prop}, \eqref{F1F2-hatp1}, and \eqref{5hatp1esti}, we have 
\begin{equation}\label{F1G1esti}
|\nabla_{x'}^{k} F_{1}|,|\nabla_{x'}^{k} G_{1}|\le C\delta(x')^{\frac{k}{2}},\,
\,\partial_{x_3}F_{1}=\partial_{x_3}G_{1}=0.
\end{equation}
Let $F=F_{1}+F_{2}$ and $G=G_{1}+G_{2}$ so that $F$ and $G$ satisfy
\eqref{5FGeq}. For $H(x),$ by \eqref{h1h2prop}, \eqref{Hform}, \eqref{F2G2form}, 
and \eqref{F1G1esti}, we can write 
\begin{equation}\label{H-polyform}
H=\sum_{i=1}^{3}H_{i}^{g}(x')x_{3}^{i}+\sum_{i=0}^{1}H_{i}^{b}(x')x_{3}^{i},
\end{equation}
where for any $k\ge0$
\begin{equation}\label{5Hesti}
|\nabla_{x'}^{k}H_{i}^{g}(x')|\le C\delta(x')^{\frac{3}{2}-i-\frac{k}{2}},
\quad|\nabla_{x'}^{k}H_{i}^{b}|\le C \delta(x')^{\frac12-i-\frac{k}{2}}.
\end{equation}
By \eqref{v51form}, \eqref{F2G2eq}, \eqref{F1F2-hatp1}, \eqref{Hform}, and 
\eqref{5Hesti}, we have
\begin{align}\label{5muDelta-gradhatp}
\begin{split}
&|\mu\Delta({\bf v}_{5}^{1})^{i}-\partial_{x_i}\hat{p}_{1}|\le 
C\delta(x')^{-1},~i=1,2,\quad|\Delta_{x'}({\bf v}_{5}^{1})^{(3)}|\le C 
\delta(x')^{-\frac12},\\
&\partial_{x_3x_3}({\bf v}_{5}^{1})^{(3)}=\partial_{x_3x_3}H\Big(k(x)^2-
\frac14\Big) \\
&=\partial_{x_3x_3}\Big(\sum_{i=1}^{3}H_{i}^{g}(x')x_{3}^{i}
\big(k(x)^2-\frac14\big)\Big)+ \partial_{x_3x_3}\Big(\sum_{i=0}^{1}H_{i}^{b}
(x')x_{3}^{i}\big(k(x)^2-\frac14\big)\Big):=I_{1}+I_{2},
\end{split}
\end{align}
where $|I_{1}|\le C \delta(x')^{-\frac12}$ and $|I_{2}|\le 
C\delta(x')^{-\frac32}$. Thus, the leading term in $\mu\Delta{\bf v}_{5}^{1}-\nabla\hat{p}_{1}$ is $
\mu\partial_{x_3x_3}\Big(\sum_{i=0}^{1}H_{i}^{b}(x')x_{3}^{i}\Big(k(x)^2-
\frac14\Big)\Big)$. To cancel this term, we take
\begin{equation}\label{5tilp1form}
\tilde{p}_{1}= \mu\partial_{x_3}\Big(\sum_{i=0}^{1}H_{i}^{b}(x')x_{3}^{i}
\big(k(x)^2-\frac14\big)\Big).
\end{equation}
Let $\bar{p}_{1}=\hat{p}_{1}+\tilde{p}_{1}$ and denote ${\bf f}_{5}^{1}
=\mu\Delta{\bf v}_{5}^{1}-\nabla\bar{p}_{1}$ in $\Omega_{2R}$. By 
\eqref{h1h2prop}, 
\eqref{v51form}, \eqref{F1F2-hatp1}, \eqref{H-polyform}, 
\eqref{5muDelta-gradhatp}, and \eqref{5tilp1form}, we can write ${\bf f}_{5}^{1}
$ in the following form
\begin{equation}\label{f51polyna}
\begin{aligned}
&({\bf f}_{5}^{1})^{(1)}=\sum_{i=0}^{2}S_{i1}^{1}(x')x_{3}^{i}+
\sum_{j=1}^{4}G_{j1}(x')x_{3}^{j}, \quad({\bf f}_{5}^{1})^{(2)}
=\sum_{i=0}^{2}S_{i1}^{2}(x')x_{3}^{i}, \\
&({\bf f}_{5}^{1})^{(3)}=\sum_{k=0}^{3}\tilde{S}_{k1}(x')x_{3}^{k}+
\sum_{i=1}^{5}\tilde{G}_{i1}(x')x_{3}^{i},
\end{aligned}
\end{equation}
where, for $0\le i\le 2,$ $1\le j\le 4$, $0\le k\le 3$, $1\le l\le5$, and $s\ge 
0$,
\begin{equation}\label{1-SilGilesti}
\begin{split}
&|\nabla_{x'}^{s}S_{i1}^{b}|\le C\delta(x')^{-i-\frac{s}{2}-1},~b=1,2,\quad|\nabla_{x'}^{s}G_{j1}|\le C\delta(x')^{-j-\frac{s}{2}},\\
&|\nabla_{x'}^{s}\tilde{S}_{k1}|\le C\delta(x')^{-k-\frac{s}{2}-\frac12},
\quad |\nabla_{x'}^{s}\tilde{G}_{l1}|\le C\delta(x')^{-l-\frac{s}{2}+
\frac12}.
\end{split}
\end{equation}
Then
\begin{equation}\label{u5fz1}
\begin{split}
&\Big|\sum_{i=0}^{2}S_{i1}^{b}(x')x_{3}^{i}\Big|\leq C\delta(x')^{-1},
\quad \Big|\sum_{j=1}^{4}G_{j1}(x')x_{3}^{j}\Big|\leq C,\\
&\Big|\sum_{k=0}^{3}\tilde{S}_{k1}(x')x_{3}^{k}\Big|\leq 
C\delta(x')^{-\frac12},\quad\Big|\sum_{l=1}^{5}\tilde{G}_{l1}(x')x_{3}
^{l}\Big|\leq C\delta(x')^{\frac12},
\end{split}
\end{equation}
which implies
$|({\bf f}_{5}^{1})^{(1)}|,~|({\bf f}_{5}^{1})^{(2)}|\le C \delta(x')^{-1}$ and
$|({\bf f}_{5}^{1})^{(3)}|\le C \delta(x')^{-1/2}$ in $\Omega_{2R}$. Therefore,
\begin{equation}\label{f51esti}
|{\bf f}_{5}^{1}|\le C\delta(x')^{-1}\quad \text{in}~\Omega_{2R}.
\end{equation}
By Proposition \ref{keyprop2} and \eqref{f51esti}, we have
$\|\nabla {\bf u}_{1}^{5}-\nabla{\bf v}_{5}^{1}\|_{L^{\infty}
(\Omega_{\delta(x')/2}(x'))}\le C. $
Moreover, a direct calculation gives that in $\Omega_{2R}$,
\begin{equation}\label{v51gradesti}
\begin{split}
&|\nabla{\bf v}_{5}^{1}|\le C \delta(x')^{-1},\quad|\nabla^{l}{\bf v}
_{5}^{1}|\le C\delta(x')^{-\frac{l+2}{2}}\quad \text{for}~l\ge2, \\ 
&|\bar{p}_1|\le C\delta(x')^{-3/2},\,\quad\quad|\nabla^{l}\bar{p}_{1}|
\le C\delta(x')^{-\frac{3+l}{2}}\quad\, \text{for}~l\ge1.
\end{split}
\end{equation}
Hence, it holds that
$|\nabla{\bf u}_{1}^{5}|\le |\nabla{\bf v}_{5}^{1}|+C \le C \delta(x')^{-1}.$

{\bf Step II. Second-order derivatives estimates.} Denote $${\bf f}_{5}^{l}:={\bf f}_{5}^{l-1}(x)+\mu \Delta{\bf v}_{5}^{l}(x)-
\nabla\bar{p}_{l}(x),\quad 2\le l\le m+1.$$ 
We first take
\begin{equation}\label{v52form}
{\bf v}_{5}^{2}:=
\begin{pmatrix}
F_{12}^{2}(x')x_3^{2}+F_{12}^{1}(x')x_{3}+F_{12}^{0}(x')+
\tilde{F}_{12}(x')\\\\ 
F_{22}^{2}(x')x_3^{2}+F_{22}^{1}(x')x_{2}+F_{22}^{0}(x')+
\tilde{F}_{22}(x')\\\\
F_{32}^{3}(x')x_{3}^{3}+F_{32}^{2}(x')x_{2}^{2}+F_{32}^{1}
(x')x_{3}+F_{32}^{0}(x')+\tilde{F}_{32}(x) 
\end{pmatrix}\big(k(x)^2-\frac14\big).
\end{equation} 
It follows from \eqref{u5fz1} that the leading term of $({\bf f}
_5^1)^{(b)},b=1,2$, is $\sum_{i=0}^{2}S_{i1}^{b}(x')x_{3}^{i}$, being of order $
\delta(x')^{-1}$. 
To cancel this leading term, we choose $F_{b2}^i$, $b=1,2$ and $i=0,1,2$, such that
\begin{equation}\label{fz5.51}
\mu\partial_{x_3x_3}\Big(\sum_{i=0}^{2}F_{b2}^ix_3^i\big(k(x)^2-\frac14\big)
\Big)=-\sum_{i=0}^{2}S_{i1}^b(x')x_{3}^{i},~b=1,2.
\end{equation}
By comparing the coefficients of each term of $x_3$, we have
\begin{equation}\label{u15-F1222}
\begin{aligned}
&F_{b2}^{2}(x')=-\frac{\delta(x')^2}{12\mu}S_{21}^{b}(x'),\,\,F_{b2}^{1}(x')=-
\frac{\delta(x')^2}{6\mu}S_{11}^{b}(x')+(h_1-h_{2})(x')F_{b2}^{2}(x'),\\
&F_{b2}^{0}(x')=-\frac{\delta(x')^2}{2\mu}S_{01}^{b}(x')+(h_1-h_2)(x')F_{b2}^{1}
(x')+\frac14(\varepsilon+2h_{1}(x'))(\varepsilon+2h_{2}(x'))F_{b2}^{2}(x').
\end{aligned}
\end{equation}
To ensure $\nabla\cdot {\bf v}_5^2=0$, we take
\begin{equation*}
\begin{aligned}
F_{32}^{3}(x')=&\,-\frac{\delta(x')^{2}}{5}\sum_{i=1}^{2}\partial_{x_i}
\Big(\frac{F_{i2}^{2}(x')}{\delta(x')^{2}}\Big),\\ 
F_{32}^{2}(x')=&\,-\frac{\delta(x')^2}{4}\sum_{i=1}^{2}\partial_{x_i}
\Big(\frac{F_{i2}^{1}(x')-(h_1-h_2)(x')F_{i2}^{2}(x')}{\delta(x')^2}\Big)+ (h_1-
h_2)(x')F_{32}^{3}(x'),\\
F_{32}^{1}(x')=&\,-\frac{\delta(x')^2}{3}\sum_{i=1}^{2}\partial_{x_i}
\Big(\frac{F_{i2}^{0}(x')-(h_1-h_2)(x')F_{i2}^{1}(x')}{\delta(x')^2}-
\frac{(\varepsilon+2h_1)(\varepsilon+2h_2)}{4\delta(x')^2}F_{i2}^{2}(x')\Big) \\
&\,+\frac14(\varepsilon+2h_1(x'))(\varepsilon+2h_2(x'))F_{32}^{3}(x')+
(h_1-h_2)(x')F_{32}^{2}(x'), \\
\end{aligned}
\end{equation*}
\begin{equation*}
\begin{aligned}
F_{32}^{0}(x')=&\,\frac{\delta(x')^2}{2}\sum_{i=1}^{2}\partial_{x_i}
\Big(\frac{(h_1-h_2)(x')F_{i2}^{0}(x')}{\delta(x')^2}+\frac{(\varepsilon+2h_1)
(\varepsilon+2h_2)}{4\delta(x')^2}F_{i2}^{1}(x')\Big) \\
&\,+\frac14(\varepsilon+2h_1(x'))(\varepsilon+2h_2(x'))F_{32}^{2}(x')+(h_1-
h_2)F_{32}^{1}(x'), 
\end{aligned}\qquad\quad\quad\quad\quad
\end{equation*}
such that 
\begin{align}\label{u15-1Rform} 
\nabla\cdot{\bf v}_{5}^{2}=&\, -\sum_{i=1}^{2}\partial_{x_i}\Big(\frac{(\varepsilon+2h_1(x'))
(\varepsilon+2h_2(x'))}{4\delta(x')^2}F_{i2}^{0}(x') \Big)\nonumber \\ 
&\,-\frac{(\varepsilon+2h_1(x'))(\varepsilon+2h_2(x'))}{4\delta(x')^2}F_{32}^{1}
(x')-\frac{h_1-h_2}{\delta(x')^2}F_{32}^{0}(x'):=R(x').
\end{align}
By \eqref{h1h2prop}, \eqref{1-SilGilesti}, \eqref{u15-F1222}, 
\eqref{u15-F1222}, and \eqref{u15-1Rform}, we have
\begin{equation}\label{5-F1121est}
|\nabla_{x'}^{k}F_{12}^{i}|,|\nabla_{x'}^{k}F_{22}^{i}|\le C\delta(x')^{1-
i-\frac{k}{2}},\,\,|\nabla_{x'}^kF_{31}^{j}|\le C\delta(x')^{\frac32-j-\frac{k}
{2}},\,|\nabla_{x'}^{k}R|\le C\delta(x')^{\frac12-\frac{k}{2}}.
\end{equation}

Furthermore, to make $\nabla\cdot{\bf v}_{5}^{2}=0,$ we choose $\tilde{F}_{12}
(x'),$ $\tilde{F}_{22}(x')$, and $\tilde{F}_{32}(x)$ such that
\begin{equation*}
\nabla \cdot\begin{pmatrix}
\tilde{F}_{12}(x')\big(k(x)^2-\frac14\big) \\
\tilde{F}_{22}(x')\big(k(x)^2-\frac14\big)\\
\tilde{F}_{32}(x)\big(k(x)^2-\frac14\big)
\end{pmatrix}=-R(x')\quad\text{in}~\Omega_{2R}.
\end{equation*}
Similar to the argument that led to \eqref{tilF11eq}--\eqref{hatp1eq}, we derive
\begin{equation*}
\begin{aligned}
&\tilde{F}_{12}(x')=\frac{\delta(x')^2}{2\mu}\partial_{x_1}\hat{p}_{2},
\quad\tilde{F}_{22}(x')=\frac{\delta(x')^2}{2\mu}\partial_{x_2}\hat{p}_{2},\\
&\tilde{F}_{31}(x)=-2\delta(x')k(x)R(x')-\delta(x')\partial_{x_1}k(x)\tilde{F}
_{11}(x')-\delta(x')\partial_{x_1}k(x)\tilde{F}_{21}(x'),
\end{aligned}
\end{equation*}
and
\begin{equation*}
\Delta_{x'}\hat{p}_{2}+\frac{3}{\delta(x')}
\nabla_{x'}\delta(x')\nabla_{x'}\hat{p}_{2}=\frac{12\mu R(x')}{\delta(x')^2}:=\mathcal{F}_{2}(x').
\end{equation*}
By \eqref{5-F1121est}, we have $|\mathcal{F}_{2}(x')|\le C\delta(x')^{-\frac32}$ and $|\nabla_{x'}^s\mathcal{F}_{2}(x')|\le C\delta(x')^{-\frac{s+3}{2}}$.
Using Proposition \ref{U_esti} with $\gamma=-\frac32$ and \eqref{hatp2tilF12}, we derive, for any $k\ge0$ and $b=1,2$, 
\begin{equation}\label{hatp2esti}
\begin{aligned}
&|\hat{p}_{2}|\le C\delta(x')^{-\frac12},\qquad\quad\,|\nabla_{x'}^{k}\hat{p}
_{2}|\le 
C\delta(x')^{-\frac12-\frac{k}{2}},\\
&|\nabla_{x'}^{k}\tilde{F}_{b2}|\le C\delta(x')^{1-\frac{k}{2}},
\quad|\nabla_{x'}^{k}\tilde{F}_{32}|\le C\delta(x')^{\frac32-\frac{k}{2}}.
\end{aligned}
\end{equation}
Therefore, 
by \eqref{f51polyna}, \eqref{1-SilGilesti}, \eqref{v52form}, \eqref{fz5.51}, and \eqref{5-F1121est}--\eqref{hatp2esti}, we have
\begin{align}\label{f51tilde}
\begin{split}
&({\bf f}_{5}^{1})^{(b)}+\mu \Delta({\bf v}_{5}^{2})^{(b)}-
\partial_{x_b}\hat{p}_{2}=\mu \Delta_{x'}({\bf v}_{5}
^{2})^{(b)}=\sum_{i=0}^{4}\hat{G}_{i1}^{b}(x')x_3^{i},~b=1,2,\\
&({\bf f}_{5}^{1})^{(3)}+\mu \Delta({\bf v}_{5}^{2})^{(3)}-
\partial_{x_3}\hat{p}_{2}=\sum_{i=0}^{3}\bar{S}_{i1}(x')x_3^{i}
+\sum_{i=0}^{5}\bar{G}_{i1}(x')x_3^{i},
\end{split}
\end{align}
where \begin{equation*}
\sum_{i=0}^{3}\bar{S}_{i1}(x')x_3^{i}:=({\bf f}_{5}
^{1})^{(3)}-\partial_{x_3}\hat{p}_{2}+\mu\partial_{x_3x_3}({\bf v}_{5}
^{2})^{(3)},\quad\sum_{i=0}^{5}\bar{G}_{i1}(x_1)x_2^{i}:=\mu 
 \Delta_{x'}({\bf v}_{5}^{2})^{(3)}.
\end{equation*}

For $s\ge 0,$ $b=1,2$, $i=0,1,2,3,4$, $j=0,1,2,3$, and $k=0,1,2,3,4,5$,
\begin{equation*}
|\nabla^{s}_{x'}\hat{G}^b_{i1}|\le C\delta(x')^{-i-\frac{s}{2}},\quad |
\nabla_{x'}^{s}\bar{S}_{j1}|\le C \delta(x')^{-\frac12-j-\frac{s}{2}},\quad 
|\nabla^{s}\bar{G}_{k1}|\le C\delta(x')^{\frac12-k-\frac{s}{2}}.
\end{equation*}
Thus, the leading term in 
${\bf f}_{5}^{1}+\mu \Delta{\bf v}_{5}
^{2}-\nabla\tilde{p}_{2}$ is $\sum_{i=0}^{3}\bar{S}
_{i1}(x')x_3^{i}$, which is of order $\delta(x')^{-\frac12}$. Then we take 
\begin{equation}\label{5tilp2form}
\tilde{p}_{2}=\sum_{i=0}^{3}\frac{1}{i+1}\bar{S}_{i1}(x')x_{3}^{i+1}
\quad\text{in}~\Omega_{2R}.
\end{equation}
Set $\bar{p}_{2}=\tilde{p}_{2}+\hat{p}_{2}$. By using \eqref{f51tilde} and 
\eqref{5tilp2form}, we can write 
\begin{equation*}
({\bf f}_{5}^{2})^{(b)}=\sum_{i=1}^{4}S_{i2}^{b}(x')x_{3}^{i},\,\,b=1,2,
\quad({\bf f}_{5}^{2})^{(3)}=\sum_{i=0}^{5}G_{i2}(x')x_{3}^{i}\quad\text{in}
~\Omega_{2R}.
\end{equation*}
For $k\ge0$, $0\le i\le 4$, and $0\le j\le5$,
\begin{equation}\label{Sj2Gj2esti}
|\nabla_{x'}^{k}S_{i2}|\le C \delta(x')^{-i-\frac{k}{2}},\quad
|\nabla_{x'}^{k}G_{j2}|\le C \delta(x')^{\frac12-j-\frac{k}{2}}.
\end{equation}
Therefore, $|{\bf f}_{5}^{2}|\le C$ in $\Omega_{2R}$.

By using \eqref{v52form}, \eqref{5-F1121est}, \eqref{hatp2esti},
and \eqref{5tilp2form}, a calculation gives for 
$k\ge0$ and $b=1,2,$
\begin{align}\label{v52hatp2tilp2est}
\begin{split}
&|\nabla_{x'}^{k}\partial_{x_3}^{s}({\bf v}_{5}^{2})^{(b)}|
\le\delta(x')^{1-s-\frac{k}{2}}\quad\text{for}~ s\le 4,\quad
\partial_{x_3}^{s}({\bf v}_{5}^{2})^{(b)}=0\quad\text{for}~s\ge5, \\
&|\nabla_{x'}^{k}\partial_{x_3}^{s}({\bf v}_{5}^{2})^{(3)}|
\le\delta(x')^{\frac32-s-\frac{k}{2}}\quad\text{for}~ s\le 5,\quad
\partial_{x_3}^{s}({\bf v}_{5}^{2})^{(3)}=0\quad\text{for}~s\ge 
6, \\
&|\nabla_{x'}^{k}\partial_{x_3}^{s}\tilde{p}_{2}|\le 
C\delta(x')^{\frac12-s-\frac{k}{2}}~\text{for}~0\leq s\le4,~
\partial_{x_3}^{s}\tilde{p}_{2}=0~\text{for}~s\ge 5,\,|
\nabla_{x'}^{k}\hat{p}_{2}|\le C\delta(x')^{-\frac12-\frac{k}{2}}.
\end{split}
\end{align}
By Proposition \ref{keyprop2}, it holds that
\begin{equation}\label{u51-v2est}
\|\nabla^{2}({\bf u}_{1}^{5} -{\bf v}_{5}^{1}-{\bf v}_{5}^{2})\|_{L^{\infty}
(\Omega_{\delta(x')/2}(x'))} +\|\nabla(p_{1}^{5}-\bar{p}_{1}-\bar{p}_{2})\|
_{L^{\infty}(\Omega_{\delta(x')/2}(x'))} \le C.
\end{equation}
By \eqref{v52hatp2tilp2est}, we obtain
\begin{equation}\label{v52est}
\begin{split}
&|\nabla^{m+1}{\bf v}_{5}^{2}|\le C\delta(x')^{-m}
~\text{for}~0\le m\le 2,\quad |\nabla^{m+1}{\bf v}_{5}^{2}|\le C\delta(x')^{-\frac{m+3}{2}}
~\text{for}~m\ge3,\\
&|\nabla^{m}\bar{p}_{2}|\le C\delta(x_1)^{\frac12-m}\quad\text{for}~2\le 
m \le 4,\quad|\nabla^{m}\bar{p}_{2}|\le C \delta(x')^{-\frac{m+3}{2}}\quad\text{for}
~ m\ge5,\\ 
&|\bar{p}_{2}|\leq C\delta(x_1)^{-\frac12},\quad|\nabla\bar{p}_{2}|\leq C\delta(x')^{-1}.
\end{split}
\end{equation}
It follows from \eqref{v51gradesti} and \eqref{v52est} that 
\begin{equation}\label{5gjv2barp2}
|\nabla^{2}{\bf v}^{2}(x)|+|\nabla\bar{p}^{2}(x)|\leq C\delta(x')^{-2}.
\end{equation}
By \eqref{u51-v2est} and \eqref{5gjv2barp2}, we have for $x\in\Omega_{R},$
\begin{equation*}
|\nabla^{2}{\bf u}_{5}^{1}(x)|+|\nabla{p}_{1}^{5}(x)|\le C \delta(x')^{-2} 
\quad\text{in}~\Omega_{R}.
\end{equation*}
For the estimates of $p_1^5(x)$, by the mean value theorem, 
\eqref{v51gradesti}, \eqref{u51-v2est}, and \eqref{v52est}, we have
\begin{align*}
|p_{1}^{5}(x)-p_{1}^{5}(z',0)|\leq&\,| p_{1}^{5}-\bar{p}_1-\bar{p}_2-
(p_{1}^{5}-\bar{p}_1-\bar{p}_2)(z',0)|+|\bar{p}_1+\bar{p}_2|+C\\
\leq&\, C\|\nabla (p_{1}^{5}-\bar{p}_1-\bar{p}_2)\|_{L^\infty(\Omega_{R})}+|
\bar{p}_1+\bar{p}_2|+C
\leq\, C \delta(x')^{-3/2}
\end{align*}
for a fixed point $(z',0)\in\Omega_{R}$ with $|z'|=R/2$.

{\bf Step III. Higher derivative estimates.}
To apply Proposition \ref{keyprop1} to derive higher derivative estimates, we 
successively construct auxiliary functions to make ${\bf f}$ and $g$ in 
\eqref{weqs1} 
as small as possible. In order to cancel out the singularity of $({\bf f}_{5}
^{2})^{(b)},b=1,2,$ we choose ${\bf v}_{5}^{3}(x)$ satisfying
\begin{equation}\label{v53-12}
\begin{aligned}
({\bf v}_{5}^{3}(x))^{(b)}=\Big(\sum_{i=0}^{4}F_{b3}^{i}(x')x_{3}^{i} \Big)
\big(k(x)^2-\frac14\big),~b=1,2,
\end{aligned}
\end{equation}
where 
\begin{equation}\label{u15-F1323form}
F_{b3}^{i}(x')=-\frac{\delta(x')^2S_{i2}^{b}(x')}{\mu(i+1)(i+2)}+(h_1-h_2)
(x')F_{b3}^{i+1}(x')+\frac14(\varepsilon+2h_1(x'))(\varepsilon+2h_2(x'))F_{b3}
^{i+2}(x'),
\end{equation}
such that 
\begin{equation}\label{5cancel1f12}
\mu\partial_{x_3x_3}({\bf v}_{5}^{3})^{(b)} =-({\bf f}_{5}^{2})^{(b)},\quad 
b=1,2.
\end{equation}
Here we use the convention that $F_{b3}^{i}(x')\equiv0$ if $i\notin\{0,\dots,4\}
$ and $F_{33}^{i}=0$ if $i\notin\{0,\dots,5\}$.
Then we take 
\begin{equation}\label{v53-3form}
({\bf v}_{5}^{3})^{(3)}=\sum_{i=0}^{5}F_{33}^{i}(x')x_{3}^{i},
\end{equation}
where 
\begin{equation}
\begin{aligned}
&F_{33}^{i}(x')=(h_1-h_2)(x')F_{33}^{i+1}(x')+\frac14(\varepsilon+2h_1(x'))
(\varepsilon+2h_2(x'))F_{33}^{i+2}(x') \\
&\quad-\frac{\delta(x')^2}{i+2}\sum_{j=1}^{2} \partial_{x_j}\Big(\delta(x')^{-2}
\big(F_{j3}^{i-1}(x')-(h_1-h_2)(x')F_{13}^{i}(x')-\frac14(\varepsilon+2h_1)
(\varepsilon+2h_2)F_{j3}^{i+1}(x')\big)\Big) \label{u15-3-F33form},
\end{aligned}
\end{equation}
such that 
\begin{align}\label{u15-5-2Rx'}
\nabla\cdot{\bf v}_{5}^{3}=&\,-\sum_{i=1}^{2}\partial_{x_i}\Big(\frac{(\varepsilon+2h_1(x'))
(\varepsilon+2h_2(x'))}{4\delta(x')^2}F_{i3}^{0}(x')\Big)\nonumber\\
&\,-\frac{(\varepsilon+2h_1(x'))(\varepsilon+2h_2(x'))}{4\delta(x')^2}
F_{33}^{1}(x')-\frac{h_1(x')-h_2(x')}{\delta(x')^2}F_{33}^{0}(x'):=R(x').
\end{align}
By \eqref{Sj2Gj2esti}, \eqref{u15-F1323form}, \eqref{u15-3-F33form}, and 
\eqref{u15-5-2Rx'}, 
we have for $k\ge 0$, $i=0,1,2,3,4$, and $b=1,2,$
\begin{equation}\label{5F1323esti}
|\nabla_{x'}^kF_{b3}^{i}|\le C\delta(x')^{2-i-\frac{k}{2}},\,\,|
\nabla_{x'}^kF_{33}^{i}|\le C\delta(x')^{\frac52-i-\frac{k}{2}},\,\,|\nabla_{x'}^kF_{33}^{5}|\le C\delta(x')^{-\frac52-\frac{k}{2}},
\end{equation}
and $|\nabla_{x'}^{k}R|\le C \delta(x')^{\frac32-\frac{k}{2}}$. By \eqref{v53-12}, \eqref{5cancel1f12}, and \eqref{v53-3form}, a calculation gives 
\begin{equation}\label{5f12+Delv13}
\begin{aligned}
&{\bf f}_{1}^{2}+\mu\Delta{\bf v}_{5}^{3}=\begin{pmatrix}
\mu \Delta_{x'}({\bf v}_{5}^{3})^{(1)}\\
\mu\Delta_{x'}({\bf v}_{5}^{3})^{(2)}\\
({\bf f}_{5}^{2})^{(3)}+\mu\Delta({\bf v}_{5}^{3})^{(3)}
\end{pmatrix}
=\begin{pmatrix}
\sum_{i=0}^{6}H_{i3}^{1}(x')x_3^{i}\\
\sum_{i=0}^{6}H_{i3}^{2}(x')x_3^{i}\\
\sum_{i=0}^{7}G_{i3}(x')x_3^{i}+\sum_{i=0}^{5}\tilde{S}_{i3}(x')x_3^{i}
\end{pmatrix},
\end{aligned}
\end{equation}
where 
$$\sum_{i=0}^{7}G_{i3}(x')x_3^{i}=\mu\Delta_{x'}({\bf v}_{5}^{3})^{(3)}
\quad\mbox{and}\quad \sum_{i=0}^{5}\tilde{S}_{i2}(x')x_3^{i}= ({\bf f}_{5}
^{2})^{(3)}+\mu\partial_{x_3x_3}({\bf v}_{5}^{3})^{(3)}.$$
By \eqref{5F1323esti}, we obtain for $k\ge0,$
\begin{equation}\label{5HGSi3esti}
|\nabla_{x'}^{k}H_{i3}^{b}|\le C \delta(x')^{1-i-\frac{k}{2}},\,\,b=1,2,\,\,|
\nabla_{x'}^{k}G_{i3}|\le\delta(x')^{\frac32-i-\frac{k}{2}},\,\,\,|
\nabla_{x'}^{k}\tilde{S}_{i3}|\le C\delta(x')^{\frac12-i-\frac{k}{2}}.
\end{equation}
By \eqref{5f12+Delv13} and \eqref{5HGSi3esti}, 
\begin{equation*}
|\Delta_{x'}({\bf v}_{5}^{3})^{(b)}|\leq C\delta(x'),\,\,b=1,2,\, 
|\Delta_{x'}({\bf v}_{5}^{3})^{(3)}|\leq C\delta(x')^{\frac32},\,|({\bf f}_{5}^{2})^{(3)}+\mu\partial_{x_3x_3}({\bf v}_{5}^{3})^{(3)}|
\leq C\delta(x')^{\frac12}.
\end{equation*}
Thus, the leading term in $({\bf f}_{5}^{2})^{(3)}+\mu\partial_{x_3x_3}({\bf v}_{5}
^{3})^{(3)}$ is $\sum_{i=0}^{5}\tilde{S}_{i3}(x')x_3^{i}.$ To cancel this leading term, we take 
\begin{equation}\label{5barp3form}
\bar{p}_{3}=\sum_{i=0}^{5}\frac{\tilde{S}_{i3}(x')}{i+1}x_3^{i+1}
\quad\text{in}~\Omega_{2R},
\end{equation}
such that
$\partial_{x_{3}}\bar{p}_{3}=\sum_{i=0}^{5}\tilde{S}_{i3}(x')x_3^{i}.$

From \eqref{5f12+Delv13} and \eqref{5barp3form},
we can write ${\bf f}_{1}^{3}$ as polynomials in $x_3$:
\begin{equation*}
\begin{split}
({\bf f}_{5}^{3})^{(b)}=&\,\sum_{i=1}^{6}\Big(H_{i3}-
\frac1i\partial_{x_{b}}
\tilde{S}_{(i-1)3}\Big)(x')x_3^{i} +H_{03}(x'):=\sum_{i=0}^{6}S_{i3}^{b}
(x')x_3^{i},~b=1,2,
\\({\bf f}_{5}^{3})^{(3)}=&\,\sum_{i=0}^{7}G_{i3}(x')x_3^{i},
\end{split}
\end{equation*} 
where $S_{i3}^{b}(x')=H_{i3}(x')-\frac1i\partial_{x_{b}}\tilde{S}_{(i-1)3}(x')$, $1\le i\le6$, and $S_{03}^{b}(x')=H_{03}(x')$. 
In view of \eqref{5HGSi3esti}, we have for $k\ge0$,
\begin{equation*}
|\nabla_{x'}^{k}S_{i3}^{b}|\le C\delta(x')^{1-i-\frac{k}{2}},\quad|
\nabla_{x'}^{k}G_{i3}|\le C\delta(x')^{\frac32-i-\frac{k}{2}}.
\end{equation*}
Thus, we have
$|({\bf f}_{5}^{3})^{(b)}|\leq C\delta(x'),\,b=1,2$, $|({\bf f}_{5}
^{3})^{(3)}|\leq C\delta(x')^{\frac32}$ in $\Omega_{2R}$, and then
\begin{equation}\label{f53esti}
|{\bf f}_{5}^{3}|\le |({\bf f}_{5}^{3})^{(1)}|+|({\bf f}_{5}^{3})^{(2)}|+|
({\bf f}_{5}^{3})^{(3)}|\le C\delta(x').
\end{equation}
By using \eqref{v53-12}, \eqref{v53-3form}, \eqref{5F1323esti}, 
\eqref{5HGSi3esti}, and \eqref{5barp3form}, a calculation gives for 
$k\ge0$ and $b=1,2,$
\begin{align}\label{v53barp3est}
\begin{split}
&|\nabla_{x'}^{k}\partial_{x_3}^{s}({\bf v}_{5}^{3})^{(b)}|
\le\delta(x')^{2-s-\frac{k}{2}}\quad\text{for}~ s\le 6,\quad
\partial_{x_3}^{s}({\bf v}_{5}^{3})^{(b)}=0\quad\text{for}~s\ge7, \\
&|\nabla_{x'}^{k}\partial_{x_3}^{s}({\bf v}_{5}^{3})^{(3)}|
\le\delta(x')^{\frac52-s-\frac{k}{2}}\quad\text{for}~ s\le 7,\quad
\partial_{x_3}^{s}({\bf v}_{5}^{3})^{(3)}=0\quad\text{for}~s\ge 
8, \\
&|\nabla_{x'}^{k}\partial_{x_3}^{s}\bar{p}_{3}|\le 
C\delta(x')^{\frac32-s-\frac{k}{2}}\quad\text{for}~0\leq s\le6,\quad
\partial_{x_3}^{s}\bar{p}_{3}=0\quad\text{for}~s\ge 7.
\end{split}
\end{align}
Denote ${\bf v}^{3}={\bf v}_{1}^{1}+{\bf v}_{1}^{2}+{\bf v}_{1}^{3}$ and ${\bar{p}
^3}=\bar{p}_{1}+\bar{p}_{2}+\bar{p}_{3}.$ Recall that $\nabla\cdot{\bf v}
^{3}
=R(x')$ with $|R(x')|\le C\delta(x')^{3/2}.$ By virtue of \eqref{f53esti} 
and Proposition \ref{keyprop1}, we have for any $m\ge1,$
\begin{equation}\label{5u11-v3est}
\begin{aligned}
&\|\nabla^{m+1}({\bf u}_{1}^{1} -{\bf v}^{3})\|_{L^{\infty}(\Omega_{\delta(x')/
2}(x'))} +\|\nabla(p_{1}^{1}-\bar{p}^{3})\|_{L^{\infty}
(\Omega_{\delta(x')/2}(x'))} \\
& \le C(\delta(x')^{2-m} + \delta(x')^{\frac32-m})\le 
C\delta(x')^{\frac{3}{2}-m}.
\end{aligned}
\end{equation}
By \eqref{v53barp3est}, we summarize
\begin{equation}\label{v53est}
\begin{split}
&|\nabla^{m+1}{\bf v}_{5}^{3}|\le C\delta(x')^{1-m}
~\,\text{for}~0\le m\le 4,\quad|\nabla^{m+1}{\bf v}_{5}^{3}|\le C\delta(x')^{-\frac{m+3}{2}}
~\,\text{for}~m\ge5,\\
&|\nabla^{m}\bar{p}_{3}|\le C\delta(x')^{\frac32-m}~\,\text{for}~0\le 
m \le 5,\quad\quad|\nabla^{m}\bar{p}_{3}|\le C \delta(x')^{-\frac{m+3}{2}}~\,\text{for}
~ m\ge6.
\end{split}
\end{equation}
By \eqref{v51gradesti}, \eqref{v52est}, and \eqref{v53est}, we have
\begin{equation}\label{5v3m+1deri}
|\nabla^{m+1}{\bf v}^{3}|+|\nabla^{m}\bar{p}^3|\le C \delta(x')^{-\frac{m+3}
{2}}.
\end{equation}
It follows from \eqref{5u11-v3est} and \eqref{5v3m+1deri} that for $m\ge1$,
\begin{equation*}
|\nabla^{m+1}{\bf u}_{1}^{5}|+|\nabla^{m}p_{1}^{5}|\le C \delta(x')^{-
\frac{m+3}{2}}+\delta(x')^{\frac32-m},
\end{equation*}
which implies that in $\Omega_{R},$
\begin{equation*}
\begin{aligned}
&|\nabla^{m+1}{\bf u}_{1}^{5}|+|\nabla^{m}p_{1}^{5}|\le C \delta(x')^{-
\frac{m+3}{2}}~\text{for}~m\le5,~~
|\nabla^{m+1}{\bf u}_{1}^{5}|+|\nabla^{m}p_{1}^{5}|\le C
\delta(x')^{\frac32-m}~\text{for}~m\ge6.
\end{aligned}
\end{equation*}
This finishes the proof of Proposition \ref{u15-esti}.
\end{proof}

\section{Proofs of the main theorems}\label{sec6}
In this section, we complete the proof of Theorems \ref{main thm1}, \ref{main thm2}, and \ref{main thm3} by using Propositions \ref{u11-general}, \ref{u13-esti}, \ref{u14-esti}, and \ref{u15-esti}.

\subsection{Decomposition of the solution}
Recall that the rigid displacement space in dimension three is given by
$\Psi=\mathrm{span}\Big\{{\boldsymbol\psi}_{1},
{\boldsymbol\psi}_{2},
{\boldsymbol\psi}_{3},
{\boldsymbol\psi}_{4},
{\boldsymbol\psi}_{5},
{\boldsymbol\psi}_{6}
\Big\},$
where ${\boldsymbol\psi}_{\alpha}$ is defined in Section \ref{sec_pri}.
It follows from $e({\bf u})=0$ that
$${\bf u}=\sum_{\alpha=1}^{6}C_{i}^{\alpha}{\boldsymbol\psi}_{\alpha}
\quad\mbox{in}~D_i,\quad i=1,2,$$
where $C_{i}^{\alpha}$ are constants to be determined by $({\bf u},p)$ 
later. We decompose the solution of \eqref{maineqs} as follows:
\begin{align}
{\bf u}(x)=&\,\sum_{i=1}^{2}\sum_{\alpha=1}^{6}C_i^{\alpha}{\bf u}_{i}
^{\alpha}(x)+{\bf u}_{0}(x),\quad 
p(x)=\sum_{i=1}^{2}\sum_{\alpha=1}^{6}C_i^{\alpha}p_{i}^{\alpha}(x)+p_{0}
(x),\quad x\in\,\Omega,\label{ud}
\end{align}
where ${\bf u}_{i}^{\alpha},{\bf u}_{0}\in{C}^{m+1}(\Omega;\mathbb R^3),~p_{i}
^{\alpha}, p_0\in{C}^{m}(\Omega)$, respectively, satisfying
\begin{equation}\label{equ_v12D}
\begin{cases}
\mu\Delta{\bf u}_{i}^\alpha=\nabla p_{i}^{\alpha},\quad\nabla\cdot 
{\bf u}_{i}^{\alpha}=0&\mathrm{in}~\Omega,\\
{\bf u}_{i}^{\alpha}={\boldsymbol\psi}_{\alpha}&\mathrm{on}~\partial{D}
_{i},\\
{\bf u}_{i}^{\alpha}=0&\mathrm{on}~\partial{D_{j}}\cup\partial{D},~j\neq 
i,
\end{cases}\quad i=1,2,
\end{equation}
and
\begin{equation}\label{equ_v32D}
\begin{cases}
\mu \Delta{\bf u}_{0}=\nabla p_0,\quad\nabla\cdot {\bf u}_{0}
=0&\mathrm{in}~\Omega,\\
{\bf u}_{0}=0&\mathrm{on}~\partial{D}_{1}\cup\partial{D_{2}},\\
{\bf u}_{0}={\boldsymbol\varphi}&\mathrm{on}~\partial{D}.
\end{cases}
\end{equation}
We rewrite \eqref{ud} as
\begin{equation*}
{{\bf u}}=\sum_{\alpha=1}^{6}\left(C_{1}^{\alpha}-C_{2}^{\alpha}
\right){\bf u}_{1}^{\alpha}
+ {\bf u}_{b}\quad\mbox{and}\quad p=\sum_{\alpha=1}^{6}\left(C_{1}^{\alpha}-C_{2}^{\alpha}\right)p_{1}
^{\alpha}+p_{b}\quad\mbox{in}~\Omega,
\end{equation*}
where ${\bf u}_{b}:=\sum_{\alpha=1}^{6}C_{2}^{\alpha}({\bf u}_{1}^{\alpha}+{\bf u}
_{2}^{\alpha})+{\bf u}_{0}$ and $ p_{b}:=\sum_{\alpha=1}^{6}C_{2}^{\alpha}
(p_{1}^{\alpha}+p_{2}^{\alpha})+p_{0}.$

\subsection{Upper bounds for the higher derivatives for \texorpdfstring{$h_1(x')\neq h_2(x')$}{}}

Since ${\bf u}_{1}^\alpha+{\bf u}_{2}^\alpha-{\boldsymbol\psi}_\alpha=0$ on $
\partial D_1\cup\partial D_2$, $\alpha=1,2,\ldots,6$, it was proved in 
\cite[Proposition 2.4]{LX1} that, $\nabla\big({\bf u}_{1}^\alpha+{\bf u}_{2}
^\alpha-
{\boldsymbol\psi}_\alpha\big)$ has no singularity in the narrow region. In fact, 
for $m\ge1$, $\nabla^m\big({\bf u}_{1}^\alpha+{\bf u}_{2}^\alpha\big)$ and $
\nabla^{m}\big(p_{1}^{\alpha}+p_{2}^{\alpha})$ are also bounded in the narrow 
region. The proof is the same as that in \cite[Proposition 6.1]{DLTZ1}, thus we omit the details here.
\begin{prop}\label{propu0}
Let ${\bf u}_0,~{\bf u}_{i}^{\alpha}\in{C}^{m+1,\alpha}(\Omega;\mathbb 
R^3),~p_0,~p_{i}^{\alpha}\in{C}^{m,\alpha}(\Omega)$ be the solution to 
\eqref{equ_v12D} and \eqref{equ_v32D}. Then, for any $m\ge0$, we have
\begin{equation*}
\|\nabla^{m+1}({\bf u}_{1}^\alpha+{\bf u}_{2}^\alpha)\|_{L^{\infty}(\Omega_R)}+
\|\nabla^{m}\big(p_{1}^{\alpha}+p_{2}^{\alpha}-(p_{1}^{\alpha}+p_{2}^{\alpha})
(z',0)\big)\|_{L^{\infty}(\Omega_R)}\leq C
\end{equation*}
and
\begin{equation*}
\|\nabla^{m+1}{\bf u}_{0}\|_{L^{\infty}(\Omega_R)}+\|\nabla^{m}\big(p_0-
p_0(z',0)\big)\|_{L^{\infty}(\Omega_R)}\leq C,
\end{equation*}
for some point $z=(z', 0)\in\Omega_{R}$ with $|z'|=R/2$.
\end{prop}
For the coefficients $C_i^\alpha$, $i=1,2$, $\alpha=1,2,\ldots,6$, which depend 
only on the gradient estimates, we have the following result, which can be found in
\cite[Proposition 2.5]{LX2}.
 \begin{prop}\label{lemCialpha}
Let $C_{i}^{\alpha}$ be defined in \eqref{ud}. Then
$$|C_i^{\alpha}|\leq\,C,\quad\,i=1,2,~\alpha=1,2,\ldots,6,$$
and
\begin{equation*}
|C_1^\alpha-C_2^\alpha|\leq \frac{C}{|\ln\varepsilon|},\,\alpha=1,2,5,6,\quad 
|C_1^3-C_2^3|\leq C\varepsilon,\quad|C_1^4-C_2^4|\leq C. 
\end{equation*}
\end{prop}

We now prove the main result, Theorem \ref{main thm1}. 
\begin{proof}[Proof of Theorem \ref{main thm1}]
By virtue of Propositions \ref{u11-general}, \ref{u13-esti}, \ref{u15-esti}, 
\ref{propu0}, and \ref{lemCialpha}, we have, for $x\in\Omega_{R}$,
\begin{equation*}
\begin{aligned}
|\nabla{\bf u}(x)|&\le \sum_{\alpha=1}^{6}|C_{1}^{\alpha}-C_{2}^{\alpha}||
\nabla{\bf u}_{1}^{\alpha}(x)|+C\nonumber\\ 
&\le\frac{C}{|\ln\varepsilon|\delta(x')}+\frac{C\varepsilon} 
{\delta(x')^{\frac{3}{2}}}+\frac{C}{\delta(x')^{\frac12}}+C\le \frac{C}{|\ln\varepsilon|\delta(x')}+\frac{C}{\delta(x')^{\frac12}}\leq \frac{C(1+|\ln\varepsilon||x'|)}{|\ln\varepsilon|
\delta(x')}, 
\end{aligned}
\end{equation*}
and
\begin{equation}\label{p-pz_1esti}
\begin{aligned}
|p(x)-p(z',0)|\leq&\, \sum_{\alpha=1}^{6}|(C_{1}^{\alpha}-C_{2}^{\alpha})
(p_{1}^{\alpha}(x)-p_{1}^{\alpha}(z',0))|+C \\
\leq&\,\frac{C}{|
\ln\varepsilon|\delta(x')^{\frac32}} +\frac{C\varepsilon}{\delta(x')^2}+C 
\leq\, \frac{C}{|\ln\varepsilon|\delta(x')^{\frac32}}+C, 
\end{aligned}
\end{equation}
where $(z',0)\in\Omega_{R}$ with $|z'|=R$ is a fix point.

For higher-order derivative estimates, when $m\ge1$, by using Propositions \ref{u11-general}, \ref{u13-esti}, \ref{u15-esti}, \ref{propu0}, and \ref{lemCialpha} again, we have
\begin{equation*}
 \begin{aligned}
&|\nabla^{m+1}{\bf u}(x)|\leq\, \sum_{\alpha=1}^{6}|C_1^\alpha-C_2^\alpha||
\nabla^{m+1}{\bf u}_{1}^\alpha(x)| +C \nonumber \\
&\leq\frac{C}{|\ln\varepsilon|}(\delta(x')^{-\frac{m+3}{2}}
\chi_{m\le6}+\delta(x')^{\frac32-m}\chi_{m\ge7}) + C\varepsilon\delta(x')^{-
\frac{m+4}{2}}+ C\delta(x')^{-\frac{m+1}{2}} +C \nonumber\\
&\leq 
\begin{cases}
\frac{C}{|\ln\varepsilon|}\delta(x')^{-\frac{m+3}{2}}+C\delta(x')^{-\frac{m+1}{2}} \quad\text{for}~m\le6,\\
\frac{C}{|\ln\varepsilon|}\delta(x')^{\frac32-m}+C\delta(x')^{-
 \frac{m+1}{2}} \quad\text{for}~m\ge7.
\end{cases}
\end{aligned}
\end{equation*}
Similarly,
\begin{equation*}
\begin{aligned}
|\nabla^{m}p(x)|\leq\, \sum_{\alpha=1}^{6}|C_{1}^{\alpha}-C_{2}^{\alpha}||
\nabla^{m}p_{1}^{\alpha}|+C\leq\begin{cases}
\frac{C}{|\ln\varepsilon|}\delta(x')^{-\frac{m+3}{2}}+C\delta(x')^{-
\frac{m+1}{2}} \quad\text{for}~m\le6,\\
\frac{C}{|\ln\varepsilon|}\delta(x')^{\frac32-m}+C\delta(x')^{-
\frac{m+1}{2}} \quad\,\,\,\text{for}~m\ge7.
\end{cases}
\end{aligned}
\end{equation*}
This completes the proof of Theorem \ref{main thm1}.
\end{proof}

\subsection{Upper bounds for the higher derivatives for \texorpdfstring{$h_{1}(x')=h_{2}(x')=\frac{1}{2}h(x')$}{}}

If $h_1(x')= h_2(x')$, we can improve the estimates in Proposition \ref{u11-general}, by constructing a family of auxiliary functions ${\bf v}_{\alpha}^{l}$, \(\alpha=1,2\) in terms of Green function, similar to \cite[Proposition 3.3]{DLTZ1}.
\begin{prop}\label{u11-symmetric}
Under the same assumption as in Theorem \ref{main thm1}, let ${\bf u}_{1}^{\alpha}$ 
and $p_{1}^{\alpha}$ be the solution to \eqref{u,peq1} with $\alpha=1,2$ and 
$h_1(x')=h_2(x')=\frac{1}{2}h(x')$. Then for sufficiently small 
$0<\varepsilon<1/2$, we have
$$\mbox{(i)}~|\nabla{\bf u}_{1}^{\alpha}(x)|\le C(\varepsilon+|x'|^2)^{-1},
\quad |p_1^\alpha(x)-p_1^\alpha(z',0)|\leq C(\varepsilon+|x'|^2)^{-1/2},\quad 
x\in\Omega_{R},$$
for some point $z=(z',0)\in \Omega_{R}$ with $|z'|=R/2$; and
 
\begin{equation*}
\mbox{(ii)}~|\nabla^{m+1}{\bf u}_{1}^{\alpha}(x)|+|\nabla^{m}{p}_{1}^{\alpha}(x)|\le C(\varepsilon+|
x'|^2)^{-\frac{m+2}{2}},\quad\mbox{for}~ m\ge 1,\quad x\in\Omega_{R}.\quad\quad
\end{equation*}
\end{prop}
\begin{proof}
The proof is similar to that in \cite[Proposition 3.3]{DLTZ1}.
We only present the key steps when $\alpha=1$, since the case when $\alpha=2$ is analogous. Since $h_1(x')=h_2(x')=\frac{1}{2}h(x')$, $\delta(x')=\varepsilon+h(x')$, 
and the Keller-type function becomes $k(x)=\frac{x_3}{\delta(x')}$ in $\Omega_{2R}$.
In order to apply Proposition \ref{keyprop1}, we construct a series of 
auxiliary functions ${\bf v}_{1}^{l}$ and $\bar{p}_{l}$, such that \(|\mu\Delta{\bf v}_{1}^{l}-\nabla\bar{p}_{l}|\) as small as possible.

We begin by choosing ${\bf v}_{1}^{1}=(({\bf 
v}_{1}^{1})^{(1)},({\bf v}_{1}^{1})^{(2)},({\bf v}_{1}^{1})^{(3)})$ 
with
\begin{equation}\label{s-v11-1def}
({\bf v}_{1}^{1})^{(1)}(x)=\frac{x_3}{\delta(x')}+\frac{1}{2},\quad ({\bf v}_{1}
^{1})^{(2)}(x)=0\quad\text{for}~x\in\Omega_{2R},
\end{equation}
and 
\begin{equation}\label{s-v11-3def}
({\bf v}_{1}^{1})^{(3)}(x)=-\int_{-\frac{1}{2}\delta(x')}^{x_3}\partial_{x_1}
({\bf v}_{1}^{1})^{(1)}(x',y)\ dy \quad \text{for}~x\in\Omega_{2R},
\end{equation} 
One can see the leading 
term in $\Delta{\bf v}_{1}^{1}$ is $\partial_{x_{3}x_{3}}({\bf v}_{1}
^{1})^{(3)},$ which is of order $\delta(x')^{-\frac32}.$ In order to cancel this leading term, we choose
\begin{equation}\label{s-barp1}
\bar{p}_{1}=\frac{\mu\partial_{x_{1}}h(x')}{\delta(x')^2}x_3\quad\text{in}
~\Omega_{2R}.
\end{equation}
Let ${\bf f}_{1}^{1}=\mu\Delta{\bf v}_{1}^{1}-\nabla\bar{p}_{1},$ by \eqref{s-v11-1def}, \eqref{s-v11-3def}, and \eqref{s-barp1}, a direct calculation shows that
\begin{equation*}
|{\bf f}_{1}^{1}|\le C (|({\bf f}_{1}^{1})^{(1)}|+|({\bf f}_{1}^{1})^{(2)}|+|
({\bf f}_{1}^{1})^{(3)}|)\le C\delta(x')^{-1}\quad \text{in}~\Omega_{2R}.
\end{equation*}

Denote
\begin{equation*}
{\bf f}_{1}^{j}:={\bf f}_{1}^{j-1}(x)+\mu\Delta{\bf v}_{1}^{j}(x)-
\nabla\bar{p}_{j}(x)=\sum_{l=1}^{j}(\mu\Delta{\bf v}_{1}^{l}-\nabla\bar{p}
_{l})(x)\quad\text{for}~ 2\le j\le m+1.
\end{equation*}
Next, we choose auxiliary functions${\bf v}_{1}^{l}(x)$ and $
\bar{p}_{l}(x)$ to further reduce the upper bound of ${\bf f}_1^j$ in $
\Omega_{2R}$. 
For $l\ge 2,$ we choose ${\bf v}_{1}
^{l}=(({\bf v}_{1}^{l})^{(1)},({\bf v}_{1}^{l})^{(2)}, ({\bf v}_{1}
^{l})^{(3)})^{T}$ such that
\begin{equation}\label{s-v1kdef}
\begin{split}
({\bf v}_{1}^{l})^{(i)}(x)=&\,-\frac{1}{\mu}\int_{-\frac{1}{2}\delta(x')}^{
\frac{1}{2}\delta(x')}G(x_3,y)({\bf f}_1^{l-1})^{(i)}(x',y) \ dy,\,\,i=1,2,\\
({\bf v}_{1}^{l})^{(3)}(x)=&\,-\int_{-\frac{1}{2}\delta(x')}^{x_3}
\partial_{x_1}({\bf v}_{1}^{l})^{(1)}(x',y)+\partial_{x_2}({\bf v}_{1}
^{l})^{(2)}(x',y) \ dy
\end{split}\quad \text{in}~\Omega_{2R},
\end{equation}
where
\begin{align*}
G(x_3,y)=\frac{1}{\delta(x')}\left\{ \begin{array}{l}(y+\frac12h(x')+
\frac{\varepsilon}{2})(x_3-\frac12h(x')-\frac{\varepsilon}{2}),\quad-
\frac{\varepsilon}{2}-\frac12h(x')\le y\le x_3,\\\\
(y-\frac12h(x')-\frac{\varepsilon}{2})(x_3+\frac12h(x')+\frac{\varepsilon}{2}),
\quad x_3\le y\le \frac{\varepsilon}{2}+\frac12h(x').
\end{array} \right.
\end{align*} 

Then it is easy to verify that
\begin{equation}\label{s-v1keq}
\begin{split}
\mu\partial_{x_3x_3}({\bf v}_{1}^{l})^{(i)}=&\,-({\bf f}_{1}^{l-1})^{(i)},\,\,i=1,2,\\
\nabla\cdot{\bf v}_{1}^{l}=&\,\partial_{x_1}({\bf v}_{1}^{l})^{(1)}+\partial_{x_2}
({\bf v}_{1}^{l})^{(2)}+\partial_{x_3}
({\bf v}_{1}^{l})^{(3)}=0.
\end{split}
\end{equation}
By \eqref{s-v11-1def} and \eqref{s-barp1}, we see that $({\bf f}_{1}^{1})^{(i)}$, $i=1,2$, are odd functions with respect to $x_3$. By induction, it follows from \eqref{s-v1kdef} that $({\bf v}_{1}^{l})^{(i)}$ and $\bar p_l$ are also odd functions with respect to $x_3$ for any $l\ge2$ and $i=1,2$. Thus, by \eqref{s-v1kdef}, we derive $({\bf v}_{1}^{l})^{(3)}(x', -\frac{1}{2}\delta(x'))=({\bf v}_{1}^{l})^{(3)}(x', \frac{1}{2}\delta(x'))=0$ and ${\bf v}_{1}^{l}=0\,\,\text{on}~\Gamma_{2R}^{+}
\cup\Gamma_{2R}^{-}.$

Furthermore, we take 
\begin{equation}\label{s-barpkdef}
\bar{p}_{l}=\mu\int_{0}^{x_3} \Big(\Delta_{x'}({\bf v}_{1}^{l-1})^{(3)}+
\partial_{x_3x_3}({\bf v}_{1}^{l})^{(3)} \Big) (x',y) \ dy\quad\text{in}
~\Omega_{2R},
\end{equation}
such that
\begin{equation}\label{s-partbarpkeq}
\partial_{x_3}\bar{p}_{l}=\mu\Delta_{x'}({\bf v}_{1}^{l-1})^{(3)}+ \mu 
\partial_{x_3x_3}({\bf v}_{1}^{l})^{(3)}.
\end{equation}
By \eqref{s-v1keq} and \eqref{s-partbarpkeq}, we have
\begin{align}\label{s-flform}
\begin{split}
&({\bf f}_{1}^{l})^{(i)}=({\bf f}_{1}^{l-1})^{(i)} +\mu \Delta({\bf v}_{1}
^{l})^{(i)}-\partial_{x_i}\bar{p}_{l}=\mu\Delta_{x'}({\bf v}_{1}^{l})^{(i)}-
\partial_{x_i}\bar{p}_{l},~i=1,2, \\
&({\bf f}_{1}^{l})^{(3)}=({\bf f}_{1}^{l-1})^{(3)}+\mu \Delta({\bf v}_{1}
^{l})^{(3)}-\partial_{x_3}\bar{p}_l=\mu\Delta_{x'}({\bf v}_{1}^{l})^{(3)}.
\end{split}
\end{align}

We can inductively prove the following estimates, for $j\ge1$ and 
$j=1,2, \dots,m+1$,
 \begin{equation}\label{s-v1jesti}
 |({\bf v}_{1}^{j})^{(1)}|,~|({\bf v}_{1}^{j})^{(2)}|\le C\delta(x')^{j-1},\,\,|
 ({\bf v}_{1}^{j})^{(3)}|\le C\delta(x')^{j-1/2},\,\,|\bar{p}_{j}|\le 
 C\delta(x')^{j-3/2},
 \end{equation}
and for $i=1,2$ and $0\le k\le m+1,$ 
\begin{align}\label{s-v1jhigh-esti}
 \begin{split}
 &|\nabla_{x'}^{k}\partial_{x_3}^{s}({\bf v}_{1}^{j})^{(i)}|
 \le C\delta(x')^{(2j-2s-k-2)/2}~\text{for}~0\le s\le 2j-1,\quad
 \partial_{x_3}^{s}({\bf v}_{1}^{j})^{(i)}=0~\text{for}~s\ge 2j,\\
 &|\nabla_{x'}^{k}\partial_{x_3}^{s}({\bf v}_{1}^{j})^{(3)}|
 \le C\delta(x')^{(2j-2s-k-1)/2}~\text{for}~0\le s\le 2j,\quad
 \partial_{x_3}^{s}({\bf v}_{1}^{j})^{(3)}=0~\text{for}~s\ge 
 2j+1,\\
 &|\nabla_{x'}^{k}\partial_{x_3}^{s}\bar{p}_{j}|\le 
 C\delta(x')^{(2j-2s-k-3)/2}~\text{for}~s\le 2j-1,\quad
 \partial_{x_3}^{s}\bar{p}_{j}=0\quad\text{for}~s\ge 2j.
 \end{split}
\end{align} 

By \eqref{s-flform} and \eqref{s-v1jhigh-esti}, we have
\begin{align*}
\begin{split}
 &|({\bf f}_{1}^{m+1})^{(i)}|\le C \Big(|\Delta_{x'}({\bf v}_{1}^{m+1})^{(i)}|+|
 \partial_{x_i}\bar{p}_{m+1}|\Big)\le C\delta(x')^{m-1},\,\,i=1,2,\\
 &|({\bf f}_{1}^{m+1})^{(3)}|\le C |\Delta_{x'}({\bf v}_{1}^{m+1})^{(3)}|
 \le C \delta(x')^{m-\frac{1}{2}}.
\end{split}
\end{align*}
Thus,
$$|{\bf f}_{1}^{m+1}|\le C \delta(x_1)^{m-1} \quad\mbox{and}~|\nabla^{s}{\bf 
f}_{1}^{m+1}|\le C \delta(x_1)^{m-s-1}\quad\mbox{for}~1\le s\le m.$$

Denote
\begin{equation*}
 {\bf v}^{m+1}(x)=\sum_{l=1}^{m+1}{\bf v}_{1}^{l}(x),\quad \bar{p}^{m+1}
 (x)=\sum_{l=1}^{m+1}\bar{p}_{l}(x).
\end{equation*}
It is easy to verify that ${\bf v}^{m+1}(x)={\bf u}_1^1(x)$ on $\Gamma_{2R}
^{\pm}$, $\nabla\cdot{\bf v}^{m+1}(x)=0$ in $\Omega_{2R}$. Moreover, for $1\le s\le 
 m$,
\begin{equation*}
 |{\bf f}^{m+1}(x)|=|\mu\Delta{\bf v}^{m+1}(x)-\nabla \bar{p}^{m+1}(x)|\leq C 
 \delta(x_1)^{m-1},\quad |\nabla^s{\bf f}^{m+1}(x)|\leq C \delta(x_1)^{m-s-1}.
\end{equation*}
Thus, by virtue of Proposition \ref{keyprop2}, it holds that
\begin{equation}\label{u11-vkestg}
 \|\nabla^{m+1}({\bf u}_{1}^{1} -{\bf v}^{m+1})\|
 _{L^{\infty}(\Omega_{\delta(x')/2}(x'))} +\|\nabla^{m} (p_{1}^{1}-
 \bar{p}^{m+1})\|_{L^{\infty}(\Omega_{\delta(x')/2}(x'))} 
 \le C.
\end{equation}
By
\eqref{s-v1kdef}, \eqref{s-barpkdef}, \eqref{s-v1jesti}, and 
\eqref{s-v1jhigh-esti}, we conclude
\begin{align}\label{v1lestg}
 \begin{split}
 &|\nabla^{m+1}{\bf v}_{1}^{l}|\le C\delta(x')^{-\frac{m+2}{2}}
 \quad\text{for}~l\le 
 \frac{m+1}{2},\\
 &|\nabla^{m+1}{\bf v}_{1}^{l}|\le C\delta(x')^{l-m-2}
 \quad\text{for}~\frac{m+2}{2}\le l\le m+1,\\
 &|\nabla^{m}\bar{p}_{l}|\le C \delta(x')^{-\frac{m+2}{2}}
 \quad\text{for}~ 
 l\le \frac{m}{2},\\
 &|\nabla^{m}\bar{p}_{l}|\le C\delta(x')^{\frac{2l-2m-3}{2}}
 \quad\text{for}~\frac{m+1}{2}\le l\le m+1.
 \end{split}
 \end{align}
It follows from \eqref{v1lestg} that
\begin{equation}\label{zgjdc}
 |\nabla{\bf v}^{m+1}(x)|\leq C\delta(x')^{-1}
\end{equation}
and for $m\ge 1$,
\begin{equation}\label{zgj1dc}
 |\nabla^{m+1}{\bf v}^{m+1}(x)|+|\nabla^{m}\bar{p}^{m+1}(x)|\leq 
 C\delta(x')^{-\frac{m+2}{2}}.
\end{equation}
From \eqref{u11-vkestg}, \eqref{zgjdc}, and \eqref{zgj1dc}, we have
 \begin{equation*}
 |\nabla{\bf u}_{1}^{1}(x)|\leq C\delta(x')^{-1},\quad
 |\nabla^{m+1}{\bf u}_{1}^{1}(x)|+|\nabla^{m}{p}_{1}^{1}(x)|\le C 
 \delta(x')^{-\frac{m+2}{2}} \quad\text{in}~\Omega_{R}.
 \end{equation*}

The estimate of $p_1^\alpha$ is the same as that in \cite[Proposition 3.3]{DLTZ1}, since $\bar{p}_{1}$ 
and $\bar{p}_{2}$ are odd in $x_3$, we thus omit the details here.
This completes the proof of Proposition \ref{u11-symmetric}.
\end{proof}

We are now in a position to prove Theorem \ref{main thm2}.
\begin{proof}[Proof of Theorem \ref{main thm2}]
The proof for the estimates of ${\bf u}$ and $p$ in Theorem \ref{main thm1} 
is similar to that of Theorem \ref{main thm1}. Under the assumptions of Theorem \ref{main thm2}, we have \begin{equation}\label{cxiangdeng}
C_1^\alpha=C_2^\alpha,\quad \alpha=4,5,6.
\end{equation} 
See \cite[Propositions 4.5]{LX2}. Thus, by Propositions \ref{u13-esti}, \ref{lemCialpha}, and \ref{u11-symmetric}, \eqref{p-pz_1esti} becomes
\begin{align*}
|p(x)-p(z',0)|\leq &\,\sum_{\alpha=1}^{3}|(C_{1}^{\alpha}-C_{2}^{\alpha})(p_{1}^{\alpha}(x)-p_{1}
^{\alpha}(z',0))|+C 
\leq\,\frac{C}{\delta(x')^{\frac12}|\ln\varepsilon|}+ \frac{C\varepsilon}
{\delta(x')^2}+C.
\end{align*}
Hence, by Propositions \ref{u13-esti}, \ref{lemCialpha}, and \ref{u11-symmetric}, we have for 
$m\ge 0$,
\begin{align*}
&|\nabla^{m}\sigma[{\bf u},p-p(z',0)]| \\
&\leq\, \sum_{\alpha=1}^{2}|C_1^\alpha-C_2^\alpha|\Big(|\nabla^{m+1}{\bf u}_{1}
^\alpha(x)|+|\nabla^{m}(p_{1}
^{\alpha}-p_{1}^{\alpha}(z',0))|\Big)\\
&\quad +|C_1^3-C_2^3|\Big(|\nabla^{m+1}{\bf u}_{1}
^3(x)|+|\nabla^{m}(p_{1}
^{3}-p_{1}^{3}(z',0))|+C\nonumber\\ 
&\leq \, \frac{C}{|\ln\varepsilon|}\delta(x')^{-\frac{m+2}{2}}
+C\varepsilon\delta(x')^{-\frac{m+4}{2}}+C\leq C\delta(x')^{-\frac{m+2}{2}}+C.
\end{align*}
The proof of Theorem \ref{main thm2} is finished.
\end{proof}

\subsection{Lower bounds for the higher derivatives}

\begin{proof}[Proof of Theorem \ref{main thm3}]
	By \cite[(4.42), Propositions 4.1, p48]{LX2}, we have
	\begin{equation}\label{C1-C2bound}
		|C_{1}^{\alpha}-C_{2}^{\alpha}|\ge\frac{C|\tilde{b}_{1}^{*\alpha}
			[{\boldsymbol{\varphi}}]|}{|\ln\varepsilon|},\quad\alpha=1,2,\quad|C_{1}^{3}-
		C_{2}^{3}|\ge C\varepsilon|\tilde{b}_{1}^{*3}[{\boldsymbol{\varphi}}]|.
	\end{equation}
	By Proposition \ref{lemCialpha} and \eqref{tilv31highesti}, we have
	\begin{equation*}
		|(C_{1}^{3}-C_{2}^{3})\partial_{x_3}({\bf u}_{1}^{3})^{(1)}|\le 
		C\varepsilon(|\partial_{x_3}({\bf v}_{3}^{1})^{(1)}|+C)\leq C \varepsilon^{-
			\frac{1}{2}}.
	\end{equation*}
	It follows from \eqref{s-v11-1def} that for some small $r>0$,
	\begin{equation*}
		|\partial_{x_{3}}({\bf v}_{1}^{1})^{(1)}(r\sqrt{\varepsilon},0,0 )|=|
		(\delta(x')^{-1})(r\sqrt{\varepsilon},0,0 )|\ge C\varepsilon^{-1},
	\end{equation*}
	Hence, 
	\begin{equation}\label{lower2}
		|\partial_{x_3}({\bf u}_{1}^{1})^{(1)}|(r\sqrt{\varepsilon},0,0)
		\ge\, C\varepsilon^{-1}.
	\end{equation}
	Note that $({\bf v}_{2}^{1})^{(1)}=0$, thus
	\begin{equation}\label{lower200}
		|\partial_{x_3}({\bf u}_{1}^{2})^{(1)}|(r\sqrt{\varepsilon},0,0)
		\le\,C.
	\end{equation}
	Combining \eqref{cxiangdeng} and \eqref{C1-C2bound}--\eqref{lower200}, we have
	\begin{align*}
		|\partial_{x_3}{\bf u}^{(1)}|(r\sqrt{\varepsilon},0,0)=&\,
		\Big|\sum_{\alpha=1}^{3}(C_{1}^{\alpha}-C_{2}^{\alpha})\partial_{x_3}({\bf u}
		_{1}^{\alpha})^{(1)}+\partial_{x_3}({\bf u}_{b})^{(1)}\Big|(r\sqrt{\varepsilon},
		0,0)\nonumber\\ \geq&\,|(C_{1}^{1}-C_{2}^{1})\partial_{x_3}({\bf u}_{1}
		^{1})^{(1)}|-C\varepsilon^{-\frac{1}{2}}-C \ge \frac{C|\tilde{b}_{1}^{*1}
			[\boldsymbol{\varphi}]|}{|\ln\varepsilon|}\varepsilon^{-1}.
	\end{align*}
	
    For the higher-order derivative estimates,  by Proposition \ref{lemCialpha} and
	\eqref{s-v1jhigh-esti}, we have, for $m\ge1,$ $\alpha=1,2,$
	\begin{equation}\label{mu11-est}
		|(C_{1}^{\alpha}-C_{2}^{\alpha})\nabla_{x'}^{m-1}\partial_{x_3}^2({\bf u}_{1}
		^{\alpha})^{(1)}|\le C|\ln\varepsilon|^{-1}\Big|\sum_{i=1}^{m+1}\nabla_{x'}
		^{m-1}\partial_{x_3}^2({\bf v}_{\alpha}^{i})^{(1)}+C\Big|\le C|\ln\varepsilon|
		^{-1}\varepsilon^{-\frac{m+1}{2}}.
	\end{equation}
	It follows from \eqref{tildev21form} that for small $r>0$, which may depend 
	on $m$,
	\begin{equation*}
		|\nabla_{x'}^{m-1}\partial_{x_3}^{2}(\tilde{\bf v}_{3}^{1})^{(1)}
		(r\sqrt{\varepsilon},0,0)|\ge C\varepsilon^{-\frac{m+4}{2}},
	\end{equation*}
	and by \eqref{F1121est}, \eqref{tilF11esti}, and \eqref{v3lest},
	$$|\nabla_{x'}^{m-1}\partial_{x_3}^2(\hat{\bf v}_{3}^{1})^{(1)}|+\sum_{i=2}
	^{m+1}|\nabla_{x'}^{m-1}\partial_{x_{3}}^{2}({\bf v}_{3}^{i})^{(1)}|\le 
	C\varepsilon^{-\frac{m+2}{2} }.$$ 
	Hence,
	\begin{equation}\label{mu13esti}
		|\nabla_{x'}^{m-1}\partial_{x_3}^2({\bf u}_{1}^{3})^{(1)}|(r\sqrt{\varepsilon},
		0,0)\ge C\varepsilon^{-\frac{m+4}{2}}.
	\end{equation}
	Combining \eqref{cxiangdeng}, \eqref{C1-C2bound}, \eqref{mu11-est}, and 
	\eqref{mu13esti}, we have 
	\begin{align*}
		|\nabla_{x'}^{m-1}\partial_{x_3}^{2}{\bf u}^{(1)}|(r\sqrt{\varepsilon},0,0)=&\,
		\Big|\sum_{\alpha=1}^{3}(C_{1}^{\alpha}-C_{2}^{\alpha})\nabla_{x'}^{m-1}
		\partial_{x_3}^{2}({\bf u}_{1}^{\alpha})^{(1)}+\nabla_{x'}^{m-1}\partial_{x_3}
		^2({\bf u}_{b})^{(1)}\Big|(r\sqrt{\varepsilon},0,0)\nonumber\\
		\geq&\,|(C_{1}^{3}-C_{2}^{3})\nabla_{x'}^{m-1}\partial_{x_3}^2
		({\bf u}_{1}^{3})^{(1)}|-C|\ln\varepsilon|^{-1}\varepsilon^{-\frac{m+1}{2}}
		\ge C|\tilde{b}_{1}^{*3}[\boldsymbol{\varphi}]|\varepsilon^{-\frac{m+2}{2}}.
	\end{align*}
	
	For the lower bound of the Cauchy stress $\sigma[{\bf u},p-p(z',0)]$, by using 
	\eqref{C1-C2bound}, \eqref{mu11-est}, and \eqref{mu13esti}, we have
	\begin{align*}
		&|\nabla_{x'}^{m-1}\partial_{x_3}\sigma[{\bf u},p-p(z',0)]|(r\sqrt{\varepsilon},
		0,0)\\ 
		&\geq\,\Big|\sum_{\alpha=1}^{3}(C_{1}^{\alpha}-C_{2}^{\alpha})
		\nabla_{x'}^{m-1}\partial_{x_3}\Big(2\mu e({\bf u}_{1}^{\alpha})-(p_{1}
		^{\alpha}-p_{1}^{\alpha}(z',0))
		\mathbb{I}\Big) \Big|(r\sqrt{\varepsilon},0,0)-C\\
		&\geq\, C|(C_{1}^{3}-C_{2}^{3})\nabla_{x'}^{m-1}\partial_{x_3}e_{13}({\bf u}_{1}
		^{3})|(r\sqrt{\varepsilon},0,0)-C |\ln\varepsilon|^{-1}\varepsilon^{-\frac{m+1}
			{2}} \ge C|\tilde{b}_{1}^{*3}[\boldsymbol{\varphi}]|\varepsilon^{-\frac{m+2}
			{2}}.
	\end{align*}
	
	This completes the proof of Theorem \ref{main thm3}.
\end{proof}

\appendix
\section{Proof of Proposition \ref{keyprop1}}

By the regularity of ${\bf u}$ in $\Omega \setminus \Omega_R$ and the definition of ${\bf v}$ stated at the beginning of Section \ref{sec_pri}, we have 
\[\|{\bf w}\|_{C^{m+1}(\Omega \setminus \Omega_R)} + \|q\|_{C^m(\Omega \setminus \Omega_R)} \leq C.\]
To prove Proposition \ref{keyprop1}, we decompose the solution of \eqref{weqs1} into two parts
$${\bf w}:={\bf w}_1+{\bf w}_2,\quad q:=q_1+q_2\quad\mbox{in}~\Omega_{2R},$$
where $({\bf w}_1,q_1)$ satisfies a nonhomogeneous Stokes equation
\begin{align}\label{weqs10}
\begin{cases}
-\mu\Delta{\bf w}_1+\nabla q_1={\bf f}\quad&\mathrm{in}\quad\Omega_{2R},\\
\nabla\cdot {\bf w}_1=0\quad&\mathrm{in}\quad\Omega_{2R},\\
{\bf w}_1=0\quad&\mathrm{on}\quad\partial \Omega_{2R},
\end{cases}
 \end{align}
while $({\bf w}_2,q_2)$ satisfies a homogeneous Stokes equation
\begin{align}\label{weqs2}
\begin{cases}
-\mu\Delta{\bf w}_2+\nabla q_2=0\quad&\mathrm{in}\quad\Omega_{2R},\\
\nabla\cdot {\bf w}_2=g\quad&\mathrm{in}\quad\Omega_{2R},\\
{\bf w}_2=0\quad&\mathrm{on}\quad\Gamma_{2R}^\pm,
\end{cases}
\end{align}
and $\|{\bf w}_{2}\|_{C^{m+1}(\Omega \setminus \Omega_R)} + \|q_{2}\|_{C^m(\Omega \setminus \Omega_R)} \leq C.$

For higher derivative estimates of $({\bf w}_1,q_1)$, we have the following result, which can be found in \cite[Proposition 2.1]{DLTZ1}.
\begin{prop}\label{keyprop2}[\cite{DLTZ1}]
Let ${\bf w}_1\in H^{1}(\Omega_{\frac{3}{2}R})\cap L^{\infty}(\Omega_{2R}
\setminus\Omega_{R})$ and $p_1\in L^{\infty}(\Omega_{2R}\setminus\Omega_{R})$ be the solution to \eqref{weqs10}. For any $m\ge1$, if
\begin{equation*}
|\nabla^k{\bf f}(x)|\leq 
C\delta(x')^{l-k}, \quad\forall\,0\le k\leq m,~l\ge -\frac{3}{2}\quad\mbox{for}
~x=(x',x_d)\in\Omega_{2R},
\end{equation*}
then it holds that
\begin{equation*}
\|\nabla{\bf w}_1\|_{L^{\infty}(\Omega_{\delta(z')/2}(z'))}\leq C\delta(z')^{l+1},
\quad z=(z',z_d) \in \Omega_{R},
\end{equation*}
and
\begin{equation*}
\|\nabla^{m+1}{\bf w}_1\|_{L^{\infty}(\Omega_{\delta(z')/2}(z'))} +\|
\nabla^{m}q_1\|
_{L^{\infty}(\Omega_{\delta(z')/2}(z'))}\le C\delta(z')^{l+1-m},\quad z =(z',z_d) 
\in \Omega_{R}.
\end{equation*}
\end{prop}

For higher derivative estimates of $({\bf w}_2,q_2)$, we have the following result, the proof of which will be given later.
\begin{prop}\label{keyprop3}
Let ${\bf w}_2\in H^{1}(\Omega_{\frac{3}{2}R})\cap L^{\infty}(\Omega_{2R}
\setminus\Omega_{R})$ and $p_2\in L^{\infty}(\Omega_{2R}\setminus\Omega_{R})$ be 
the solution to \eqref{weqs2}. For any $m\ge1$, if
\begin{equation*}
|\nabla^kg(x)|\leq 
C\delta(x')^{\alpha-k}, \quad\forall\,0\le k\leq m,~\alpha\ge -\frac{1}
{2}\quad\mbox{for}~x=(x',x_d)\in\Omega_{2R},
\end{equation*}
then it holds that
\begin{equation*}
\|\nabla{\bf w}_2\|_{L^{\infty}(\Omega_{\delta(z')/2}(z'))}\leq 
C\delta(z')^{\alpha},\quad z=(z',z_d) \in \Omega_{R},
\end{equation*}
and
\begin{equation*}
\|\nabla^{m+1}{\bf w}_2\|_{L^{\infty}(\Omega_{\delta(z')/2}(z'))} +\|
\nabla^{m}q_2\|_{L^{\infty}(\Omega_{\delta(z')/2}(z'))}\le 
C\delta(z')^{\alpha-m},
\quad z =(z',z_d) \in \Omega_{R}.
\end{equation*}
\end{prop}

\begin{proof}[Proof of Proposition \ref{keyprop1}]
Since
${\bf w}:={\bf w}_1+{\bf w}_2, ~q:=q_1+q_2$,
by using Propositions \ref{keyprop2} and \ref{keyprop3}, we have, for any $z=(z',z_d) \in \Omega_{R},$
\begin{equation*}
\|\nabla{\bf w}\|_{L^{\infty}(\Omega_{\delta(z')/2}(z'))}\leq \|\nabla{\bf w}_1\|_{L^{\infty}(\Omega_{\delta(z')/2}(z'))}+\|\nabla{\bf w}_2\|_{L^{\infty}(\Omega_{\delta(z')/2}(z'))}\leq 
C\delta(z')^{l+1}+C\delta(z')^{\alpha},
\end{equation*}
and
\begin{align*}
 &\|\nabla^{m+1}{\bf w}\|_{L^{\infty}(\Omega_{\delta(z')/2}(z'))}+\|
 \nabla^{m}q\|_{L^{\infty}(\Omega_{\delta(z')/2}(z'))}\\
&\leq \|\nabla^{m+1}{\bf w}_1\|_{L^{\infty}(\Omega_{\delta(z')/2}(z'))}+ \|\nabla^{m+1}{\bf w}_2\|_{L^{\infty}(\Omega_{\delta(z')/2}(z'))}+\|
 \nabla^{m}q_1\|_{L^{\infty}(\Omega_{\delta(z')/2}(z'))}\\
&\quad +\|
 \nabla^{m}q_2\|_{L^{\infty}(\Omega_{\delta(z')/2}(z'))}
\le \,
 C\delta(z')^{l+1-m}
 +C\delta(z')^{\alpha-m}.
\end{align*}
This completes the proof of Proposition \ref{keyprop1}.
\end{proof}

We next prove Proposition \ref{keyprop3} for ${\bf w}_{2}$. We first recall the classical result for the divergence equation, the proof can 
be found in \cite[section III.3]{GaldiBook}.
\begin{lemma}\label{div-eqs}
Let $\Omega$ be a bounded Lipschitz domain in $\mathbb{R}^{d}$. Given $f\in 
L^2(\Omega)$ with $\int_{\Omega}f=0.$ Then there exists a function 
${\boldsymbol{\phi}}\in H_{0}^{1}(\Omega;\mathbb{R}^{d})$ such that 
$$\mathrm{div}\,{\boldsymbol{\phi}}=f\quad \mathrm{in}~\Omega\quad\mbox{and}\quad
\|\nabla{\boldsymbol{\phi}}\|_{L^2(\Omega)}\le C\|f\|_{L^2(\Omega)},$$
where $C$ is a constant depending only on $d$. 
\end{lemma}
With Lemma \ref{div-eqs} at hand, we derive the estimates for the pressure term $q.$
\begin{lemma}
Let $({\bf{w}},q)$ be the solution to 
\begin{equation}\label{stokes_eqs}
-\mu\Delta{\bf{w}}+\nabla(q-q_{\Omega}) =0\quad\mathrm{in}~\Omega.
\end{equation}
Then it holds 
\begin{equation}\label{q-q_omegaesti}
\int_{\Omega}|q-q_{\Omega}|^2\ dx \le C \int_{\Omega}|\nabla{\bf w}|^2\ 
dx,
\end{equation}
where $q_{\Omega}:=\frac{1}{|\Omega|}\int_{\Omega}q\ dx$, and $C$ is a constant 
depending only on $d$ and $\mu$. 
\end{lemma}
\begin{proof}
By Lemma \ref{div-eqs}, there exists a function $\boldsymbol{\phi}\in H_{0}^{1}
(\Omega;\mathbb{R}^{d})$ such that 
$\mathrm{div}\,\boldsymbol{\phi}=q-q_{\Omega}$ in $\Omega$, 
with the estimates
\begin{equation}\label{grad-phi-esti}
\|\nabla\boldsymbol{\phi}\|_{L^2(\Omega)}^2 \le C_{0}\|q-q_{\Omega}\|
_{L^2(\Omega)}^2,
\end{equation}
where $C_{0}=C_{0}(d).$ Using this $\boldsymbol{\phi}$ as a test function of 
\eqref{stokes_eqs} and integrating by parts gives
\begin{equation*}
\mu\int_{\Omega}\nabla{\bf w}\cdot\nabla\boldsymbol{\phi}\ dx-
\int_{\Omega}|q-q_{\Omega}|^2\ dx=0.
\end{equation*}
By the Cauchy inequality and \eqref{grad-phi-esti}, we have
\begin{equation*}
\begin{aligned}
\int_{\Omega}|q-q_{\Omega}|^2\ dx&\le\frac{1}{2C_{0}}\int_{\Omega} |
\nabla\boldsymbol{\phi}|^2+C\int_{\Omega}|\nabla{\bf w}|^2 \ dx 
\le\frac12\int_{\Omega}|q-q_{\Omega}|^2\ dx+ C\int_{\Omega}|
\nabla{\bf w}|^2 \ dx,
\end{aligned}
\end{equation*}
which implies \eqref{q-q_omegaesti}. This completes the proof.
\end{proof}

To prove Proposition \ref{keyprop3}, we need the following $L^q$-estimate for 
incompressible Stokes equations in a bounded domain, with partially vanishing 
boundary data, which can be found in \cite[p.278]{GaldiBook}.
\begin{theorem}\label{stokes-eqs-Wkpesti}
Let $\Omega$ be an arbitrary domain in $\mathbb R^{d}$, $d\geq2$, with a 
boundary portion $\sigma$ of class $C^{m+2}$, $m\geq0$. Let $\Omega_0$ be any 
bounded subdomain of $\Omega$ with $\partial\Omega_0\cap\partial\Omega=\sigma$. 
Further, let
\begin{align*}
{\bf u}\in W^{1,q}(\Omega_0), \quad p\in L^q(\Omega_0),\quad 1<q<\infty,
\end{align*}
be such that
\begin{align*}
(\nabla{\bf u},\nabla{\boldsymbol\psi})&=-\langle{\bf f},{\boldsymbol\psi}\rangle+(p,\nabla\cdot{\boldsymbol\psi})\quad\mbox{for~all}~{\boldsymbol\psi}\in C_0^\infty(\Omega_0),\\
(\nabla\cdot{\bf u},\varphi)&=({\bf g}, \varphi)\quad \mbox{for~all}~\varphi\in C_0^\infty(\Omega_0),\quad
{\bf u}=0\quad \mbox{at}~\sigma.
\end{align*}
Then, if ${\bf f}\in W^{m,q}(\Omega_0)$, we have
${\bf u}\in W^{m+2,q}(\Omega'),~ p\in W^{m+1,q}(\Omega'),$
for any $\Omega'$ satisfying 
(1) $\Omega'\subset\Omega$,
(2) $\partial\Omega'\cap\partial\Omega$ is a strictly interior subregion of $\sigma$.
Finally, the following estimate holds
\begin{align*}
 &\|{\bf u}\|_{W^{m+2,q}(\Omega')}+\|p\|_{W^{m+1,q}(\Omega')}
 \leq C\left(\|{\bf 
 f}\|_{W^{m,q}(\Omega_0)}+\|{\bf g}\|_{W^{m+1,q}(\Omega_0)}+\|{\bf u}\|_{W^{1,q}
 (\Omega_0)}+\|p\|_{L^{q}(\Omega_0)}\right),
\end{align*}
where $C=C(d,m,q,\Omega',\Omega_0)$.
\end{theorem}

For convenience, in the sequel we omit the subscripts and denote ${\bf w}={\bf w}_{2}$. Applying Theorem \ref{stokes-eqs-Wkpesti} and a bootstrap argument yields the following proposition.
\begin{prop}\label{w,q-L_infty-esti} 
Let $({\bf w},q)$ be the solution to \eqref{weqs2}. Then for $z=(z',z_d)
\in\Omega_R$, the following estimates hold:
\begin{align*}
&\|\nabla {\bf w}\|_{L^{\infty}(\Omega_{\delta(z')/2}(z'))}
\le \,C\left(\delta(z')^{-\frac d2}\|\nabla {\bf w}\|_{L^{2}
(\Omega_{\delta(z')}(z'))}+\|{\bf g}\|_{L^{\infty}(\Omega_{\delta(z')}(z'))}+
\delta(z')\|{\nabla\bf g}\|_{L^{\infty}(\Omega_{\delta(z')}(z'))}\right),
\end{align*}
and for $m\ge1$,
\begin{align*}
&\|\nabla^{m+1}{\bf w}\|_{L^{\infty}(\Omega_{\delta(z')/2}(z'))}+\|\nabla^{m}q\|
_{L^{\infty}(\Omega_{\delta(z')/2}(z'))}\\
&\leq C\delta(z')^{-m-1}\left(\delta(z')^{1-\frac{d}{2}}\|\nabla {\bf w}\|_{L^{2}
(\Omega_{\delta(z')}(z'))}+\sum_{j=0}^{m+1}\delta(z')^{1+j}\|\nabla^{j}{\bf g}\|
_{L^{\infty}(\Omega_{\delta(z')}(z'))} \right).
\end{align*}
\end{prop}

\begin{proof}[Proof of Proposition \ref{w,q-L_infty-esti}]
For $z=(z', z_{d})\in \Omega_{R}$, we take the change of variables as before
\begin{equation*}
\left\{
\begin{aligned}
&x'-z'=\delta(z') y',\\
&x_d=\delta(z') y_{d},
\end{aligned}\right.
\end{equation*}
which transforms $\Omega_{\delta(z')}(z)$ into a nearly unit size domain $Q_{1}$, 
where 
\begin{align*}
Q_{r}=\left\{y\in\mathbb{R}^{d}:~-\frac{\varepsilon}{2\delta}-\frac{1}{\delta}h_{2}(\delta\,y'+z')<y_{d}
<\frac{\varepsilon}{2\delta}+\frac{1}{\delta}h_{1}(\delta\,y'+z'),~|y'|
<r\right\},~\delta=\delta(z'),
\end{align*}
with top and bottom boundaries denoted by $\hat\Gamma_{1}^{+}$ and $
\hat\Gamma_{1}^{-}$. For $y\in{Q}_{1}$, set
\begin{align*}
\mathcal{W}(y', y_{d})&:={\bf w}(z'+\delta(z')\,y',\delta(z')\,y_{d}), 
\quad~\mathcal{P}(y', y_{d}):=\delta(z') q(z'+\delta(z')\,y',\delta(z')\,y_{d}),
\end{align*}
and 
\begin{equation*}
\mathcal{P}_1(y', y_{d}):=\frac{1}{|Q_1|}\int_{Q_1}\mathcal{P},\quad\mathcal{G}(y', y_{d}):=\delta(z') g(z'+\delta(z')\,y',\delta(z')\,y_{d}),
\end{equation*}
It follows from \eqref{weqs2} that
\begin{align*}
\begin{cases}
-\mu\Delta\mathcal{W}+\nabla(\mathcal{P}-\mathcal{P}_1)=0,
& \mathrm{in}\,Q_1,\\
\nabla \cdot \mathcal{W} =\mathcal{G},&\mathrm{in}\,Q_1,\\
\mathcal{W}=0,&\mathrm{on}\, \hat{\Gamma}^+\cup\hat\Gamma^{-},
\end{cases}
\end{align*}
By applying Theorem \ref{stokes-eqs-Wkpesti} with $m=0$ and $q=4$, we obtain, for any $0<r<R\le 1$,
\begin{equation*}
\|\mathcal{W}\|_{W^{2,4}(Q_{r})}+\|\mathcal{P}-\mathcal{P}_{1}\|_{W^{1,4}(Q_{r})}\leq\,C\left(\|\mathcal{G}\|_{W^{1,4}(Q_R)}+\|\mathcal{W}\|_{W^{1,4}(Q_R)}+\|\mathcal{P}-\mathcal{P}_{1}\|_{L^{4}(Q_R)}\right).
\end{equation*}
Applying the interpolation inequality, for any $\theta_0\in(0,1)$,
\begin{equation*}
\|\mathcal{W}\|_{W^{1,4}(Q_{R})}\leq \theta_0\|\mathcal{W}\|_{W^{2,4}(Q_{R})} +\frac{C}{\theta_0}\|\mathcal{W}\|_{L^{4}(Q_{R})},
\end{equation*}
it follows that
\begin{align*}
&\|\mathcal{W}\|_{W^{2,4}(Q_{r})}+\|\mathcal{P}-\mathcal{P}_{1}\|_{W^{1,4}
(Q_{r})}\nonumber\\
&\leq\, C\theta_{0}\|\mathcal{W}\|_{W^{2,4}(Q_R)}+\frac{C}{\theta_{0}}\|
\mathcal{W}\|_{L^{4}(Q_R)}+C\left(\|\mathcal{G}\|_{W^{1,4}(Q_R)}+\|\mathcal{P}-
\mathcal{P}_{1}\|_{L^{4}(Q_R)}\right).
\end{align*}
 Taking $\theta_0$ small enough such that $C\theta_0\leq\frac{1}{4}$, then by 
Lemma \ref{itLemma}, we get
\begin{align*}
\|\mathcal{W}\|_{W^{2,4}(Q_{r})}+\|\mathcal{P}-\mathcal{P}_{1}\|_{W^{1,4}(Q_{r})}
\leq C\left(\|\mathcal{W}\|_{L^{4}(Q_R)}+\|\mathcal{G}\|_{W^{1,4}(Q_R)}+\|\mathcal{P}-\mathcal{P}_{1}\|_{L^{4}(Q_R)}\right).
\end{align*}
Employing the Sobolev embedding theorem and the Poincar\'e inequality, we have 
\begin{align}\label{sz1}
\begin{split}
&\|\nabla\mathcal{W}\|_{L^{\infty}(Q_{r})}+\|\mathcal{P}-\mathcal{P}_{1}\|
_{L^{\infty}(Q_{r})}\leq\, C\big(\|\mathcal{W}\|_{W^{2,4}
(Q_{r})}+\|\mathcal{P}-\mathcal{P}_{1}\|_{W^{1,4}(Q_{r})}\big)\\
&\leq\,C\big(\|\mathcal{W}\|_{W^{1,2}(Q_R)}+\|\mathcal{G}\|_{W^{1,4}(Q_R)}+
\|\mathcal{P}-\mathcal{P}_{1}\|_{L^{4}(Q_R)} \big)\\
&\leq\, C\left(\|\nabla\mathcal{W}\|_{L^{2}(Q_R)}+\|\mathcal{P}-
\mathcal{P}_{1}\|_{L^{4}(Q_{R})}+\|\mathcal{G}\|_{L^{\infty}(Q_R)}+\|
\mathcal{\nabla G}\|_{L^{\infty}(Q_R)}\right).
\end{split}
\end{align}

In order to control the last term involving $\mathcal{P}$, by means of the interpolation inequality and \eqref{q-q_omegaesti}, for any $\theta_{1}\in (0,1)$, we 
have
\begin{equation}\label{sz2}
\begin{aligned}
\|\mathcal{P}-\mathcal{P}_{1}\|_{L^{4}(Q_R)}&\leq \theta_{1} \|\mathcal{P}-
\mathcal{P}_{1}\|_{L^{\infty}(Q_R)}+C(\theta_{1})\|\mathcal{P}-\mathcal{P}
_{1}\|_{L^{2}(Q_R)} \\
&\leq \theta_{1} \|\mathcal{P}-\mathcal{P}_{1}\|_{L^{\infty}(Q_R)}
+C(\theta_1)\|\nabla\mathcal{W}\|_{L^{2}(Q_R)}.
\end{aligned}
\end{equation}
By \eqref{sz1} and \ref{sz2}, 
\begin{align*}
&\|\nabla\mathcal{W}\|_{L^{\infty}(Q_{r})}+\|\mathcal{P}-\mathcal{P}_{1}\|
_{L^{\infty}(Q_{r})}\\
& \leq C \theta_{1} \|\mathcal{P}-\mathcal{P}_{1}\|_{L^{\infty}(Q_R)}
+C\left(\|\nabla\mathcal{W}\|_{L^{2}(Q_R)}+\|\mathcal{G}\|_{L^{\infty}(Q_R)}+\|
\mathcal{\nabla G}\|_{L^{\infty}(Q_R)}\right).
\end{align*}
Chosing $\theta_1$ small enough such that $C\theta_1\leq\frac{1}{2}$ and using Lemma \ref{itLemma} again, we have
\begin{equation*}
\|\nabla\mathcal{W}\|_{L^{\infty}(Q_{r})}+\|\mathcal{P}-\mathcal{P}_{1}\|
_{L^{\infty}(Q_{r})}\le C\left(\|\nabla\mathcal{W}\|_{L^{2}(Q_R)}+\|
\mathcal{G}\|_{L^{\infty}(Q_R)}+\|\mathcal{\nabla G}\|_{L^{\infty}(Q_R)}
\right).
\end{equation*}
Taking $r=\frac{1}{2},$ $R=1$, and rescaling back to $({\bf w},q)$, we obtain
\begin{align*}
&\|\nabla {\bf w}\|_{L^{\infty}(\Omega_{\delta(z')/2}(z'))}\\
&\le
\,C\left(\delta(z')^{-\frac d2}\|\nabla {\bf w}\|_{L^{2}(\Omega_{\delta(z')}(z'))}+\|{\bf g}\|_{L^{\infty}(\Omega_{\delta(z')}(z'))}+\delta(z')\|{\nabla\bf g}\|_{L^{\infty}(\Omega_{\delta(z')}(z'))}\right).
\end{align*}
By applying Theorem \ref{stokes-eqs-Wkpesti} with $m=1,2,\dots$ and repeating the process above, we conclude that Proposition \ref{w,q-L_infty-esti} holds.
\end{proof}
From Proposition \ref{w,q-L_infty-esti}, the $W^{m+1,\infty}$ estimates of ${\bf 
w}$ depends on the local energy estimates $\|\nabla {\bf w}\|_{L^{2}
(\Omega_{\delta(z')}(z'))}$ and the pointwise estimates of $\nabla^{j}g$, 
$j=0,1,\dots,m+1$, in the narrow region. Indeed, we are able to demonstrate that 
the local energy estimates $\|\nabla {\bf w}\|_{L^{2}(\Omega_{\delta(z')}(z'))}$ are also closely related to the pointwise upper bound of $g(x)$ itself in the 
narrow region. To this end, we first prove that the global energy of ${\bf w}$ 
is bounded in the domain $\Omega_{\frac{3}{2}R}$.

\begin{lemma}\label{lemmaenergy}
Let $({\bf w},q)$ be the solution to \eqref{weqs2}. If 
\begin{equation}\label{cond_f}
|g(x)|\leq C\delta(x')^{-1/2}\quad\mbox{in}~ \Omega_{2R}
\quad\mbox{and}~ \|{\bf w}\|_{L^\infty(\Omega_{2R}\setminus\Omega_{R})}\leq C,
\end{equation}
then
$\int_{\Omega_{\frac{3}{2}R}}|\nabla {\bf w}|^{2}\leq C.$
\end{lemma}

\begin{proof}

For $1\le t<s\le 2$, let $\eta$ be a smooth function satisfying $\eta(x')=1$ if 
$|x'|<tR$, $\eta(x')=0$ if $|x'|>sR$, $0\leqslant\eta(x')\leqslant1$ if $tR\le |
x'|\le sR$, and $|\nabla_{x'}\eta(x')|\leq \frac{C}{s-t}$. Suppose $({\bf 
w},q)$ 
is the solution to \eqref{weqs2}, then it also verifies 
\begin{equation}
\label{w12}-\mu\Delta{\bf w}+\nabla(q-q_{_{sR}})=0,
\end{equation}
where $q_{_{sR}}:=\frac{1}{|\Omega_{sR}|}\int_{\Omega_{sR}}q$. Multiplying the 
first equation in \eqref{w12} by $\eta^2{\bf w}$ and integrating by parts, we 
get
\begin{align}\label{zy1}
\mu\int_{\Omega_{sR}}\eta^2|\nabla{\bf w}|^2\leq 2\mu\int_{\Omega_{sR}}|
\nabla_{x'}\eta||\eta\nabla{\bf w}||{\bf w}|+\int_{\Omega_{sR}}\big|
(q-q_{sR})\nabla\cdot({\bf w}\eta^2)\big|.
\end{align}
By virtue of the Cauchy inequality and \eqref{cond_f}, we have
\begin{equation}\label{zy2}
\begin{aligned}
\int_{\Omega_{sR}}|\nabla_{x'}\eta||\eta\nabla{\bf w}||{\bf w}|&\leq \frac{1}{8}
\int_{\Omega_{sR}}\eta^2|\nabla{\bf w}|^2+\frac{C}{(s-t)^2}\int_{\Omega_{2R}
\setminus\Omega_{R}}|{\bf w}|^2\\
&\leq \frac{1}{8}\int_{\Omega_{sR}}\eta^2|
\nabla{\bf w}|^2+\frac{C}{(s-t)^2}.
\end{aligned}
\end{equation}
It follows from \eqref{cond_f} that
\begin{equation}\label{g_L2_esti}
\int_{\Omega_{sR}}g(x)^2\leq C\int_{\Omega_{sR}}\frac{1}{\delta(x')}\leq C.
\end{equation}
For the last term on the right-hand side of \eqref{zy1}, by 
\eqref{q-q_omegaesti}, \eqref{cond_f}, \eqref{g_L2_esti}, and the Cauchy 
inequality, we derive
\begin{equation}\label{zy3}
\begin{aligned}
 \int_{\Omega_{sR}}\big|(q-q_{sR})\nabla\cdot({\bf w}\eta^2)\big|&\leq 
 \int_{\Omega_{sR}}|q-q_{sR}||g|+|q-q_{sR}||{\bf w}||\nabla \eta|\\
 &\le \frac{\mu}{4}\int_{\Omega_{sR}}|\nabla {\bf w}|^2+\frac{C}{(s-t)^2}.
\end{aligned}
\end{equation}
Substituting \eqref{zy2} and \eqref{zy3} into \eqref{zy1}, we obtain
\begin{equation*}
\int_{\Omega_{tR}}|\nabla {\bf w}|^2\leq\, \frac13\int_{\Omega_{sR}}
|\nabla{\bf w}|^2+\frac{C}{(s-t)^2}.
\end{equation*}
Hence, by Lemma \ref{itLemma} with $\rho=\frac32$ and $R=2$, we complete the 
proof.
\end{proof}

The following Caccioppoli-type inequality is a key step to build an adapted version of iteration formula, used in \cite{DLTZ1}
\begin{lemma}\label{cpll} (Caccioppoli-type Inequality) Let $({\bf w},q)$ be the 
solution to \eqref{weqs2}. Then for $z=(z',z_d)\in\Omega_{R}$ and for $0<t<s\leq 
C\sqrt{\delta(z')}$, there holds
\begin{align}\label{iterating1}
\int_{\Omega_{t}(z')}|\nabla {\bf w}|^{2}\leq\,&
\,\Big(\frac{1}{4}
+\frac{C\delta(z')^2}{(s-t)^{2}}\Big)\int_{\Omega_{s}(z')}|\nabla {\bf w}|^2+C\int_{\Omega_{s}(z')}|g|^{2}.
\end{align}
\end{lemma}

\begin{proof}
Since ${\bf w}=0$ on $\partial D_1\cup\partial D_2$, employing Poincar\'e 
inequality, it is not difficult to deduce (see \cite[(3.36) and (3.39)]{BLL2})
\begin{equation}\label{estwDwl2}
\int_{\Omega_{s}(z')}|{\bf w}|^2\leq
C\delta(z')^2\int_{\Omega_{s}(z')}|\nabla {\bf w}|^2.
\end{equation}
For $0<t<s\leq R$, let $\eta$ be a smooth function satisfying $\eta(x')=1$ if $|
x'-z'|<t$, $\eta(x')=0$ if $|x'-z'|>s$, $0\leqslant\eta(x')\leqslant1$ if 
$t\leqslant|x'-z'|<s$,
and
$|\nabla_{x'}\eta(x')|\leq \frac{2}{s-t}$. Let $q_{s;z'}=\frac{1}{|\Omega_{s}
(z')|}\int_{\Omega_{s}(z')}q$. Multiplying the equation 
$$-\mu\Delta{\bf w}+\nabla(q-q_{s;z'})=0\quad\mbox{in}~\Omega_{s}(z')$$ 
by $\eta^{2}{\bf w}$ and integrating by parts, we obtain
\begin{align}\label{energynarr}
\mu\int_{\Omega_{s}(z')}\eta^{2}|\nabla {\bf w}|^{2}=
-\mu\int_{\Omega_{s}(z')}({\bf w}\nabla {\bf w})\cdot\nabla\eta^{2}+\int_{\Omega_{s}(z')}\nabla\cdot(\eta^{2}{\bf w})\left(q-q_{s;z'}\right).
\end{align}
Next, we shall estimate the terms on the right-hand side of \eqref{energynarr} 
one by one. 

For the first term in \eqref{energynarr}, Cauchy's inequality gives
\begin{align}\label{w1Dw1eta}
\int_{\Omega_{s}(z')}|{\bf w}\nabla {\bf w}||\nabla\eta^{2}|
\leq \frac{1}{8}\int_{\Omega_{s}(z')}|\nabla {\bf w}|^{2}
+\frac{C}{(s-t)^2}\int_{\Omega_{s}(z')}|{\bf w}|^{2}.
\end{align}
For the last term, by using the second line of \eqref{weqs2}, 
\eqref{q-q_omegaesti}, and Young's inequality, we have
\begin{align}\label{etaw1q1}
 \left|\int_{\Omega_{s}(z')}\left(q-q_{s;z'}\right)\nabla\cdot(\eta^{2}{\bf w}) 
 \right|\leq \frac{\mu}{8}\int_{\Omega_{s}(z')}|\nabla {\bf w}|
 ^2+C\int_{\Omega_{s}(z')}|{\bf w}|^{2}|\nabla\eta|^{2} +C\int_{\Omega_{s}(z')}|
 g|^{2},
\end{align}
where $C$ depends only on $d$ and $\mu$.
Combining \eqref{estwDwl2} and \eqref{w1Dw1eta}--\eqref{etaw1q1} with
\eqref{energynarr} yields \eqref{iterating1}.
\end{proof}

By virtue of Lemmas \ref{lemmaenergy}, \ref{cpll}, and the iteration 
technique, we can deduce the following lemma, which illustrates the 
relationship between the local estimates of 
\(\| \nabla {\bf w}\|_{L^{2}
(\Omega_{\delta(z')})(z')}\) and the pointwise estimates of \(g(x)\) in 
\(\Omega_{2R}\).

\begin{lemma}\label{w_L2_energy}
Let ${\bf w}\in H^{1}(\Omega_{\frac{3}{2}R})$ be the solution to 
\eqref{weqs2}. If
\begin{equation}\label{jstj}
|g(x)|\leq C\delta(x')^{\alpha}, \quad\alpha\ge -\frac{1}{2}, \quad x\in 
\Omega_{2R},
\end{equation}
then $\int_{\Omega_{\delta(z')}(z')}\left|\nabla{\bf w}\right|^2 \le 
C\delta(z')^{d+2\alpha}$, for $z \in \Omega_{R}.$
\end{lemma} 
\begin{proof}
Denote $E(t):=\int_{\Omega_{t}(z')}\left|\nabla{\bf w}\right|^2$, by 
\eqref{jstj}, we have
\begin{align}\label{jbnl1}
\int_{\Omega_s(z')}|g|^2\leq Cs^{d-1}\delta(z')^{2\alpha+1},~0<s\le C\sqrt{\delta(z')}.
\end{align}
By using \eqref{iterating1} and \eqref{jbnl1}, we derive 
\begin{equation}\label{itera1}
E(t)\le\Big(\frac{1}{4}+\Big(\frac{c_0\delta(z')}{s-t}\Big)^2\Big)E(s)+Cs^{d-1}\delta(z')^{2\alpha+1},
\end{equation}
where $c_{0}$ is a constant and we fix it now. Let $k_{0}=\left[\frac{1}{8c_{0}
\sqrt{\delta(z')}}\right]$ and $t_{i}=\delta(z')+2c_{0}i\delta(z'),$ $i=0,1,2,
\dots,k_{0}$. Then taking $s=t_{i+1}$ and $t=t_i$ in \eqref{itera1}, we have%
\begin{align*}
E(t_{i})\leq \frac{1}{2}E(t_{i+1})+C(i+1)^{d-1}\delta(z')^{d+2\alpha}.
\end{align*}
After $k_{0}$ iterations, and using Lemma \ref{lemmaenergy}, we obtain
\begin{align*}
E(t_0)&\leq\,\Big(\frac{1}{2}\Big)^{k_0}E(t_{k_{0}})+C\delta(z')^{d+2a}
\sum_{i=0}^{k_0-1}\Big(\frac{1}{2}\Big)^i(i+1)^{d-1}\\
&\leq\,\Big(\frac{1}{2}\Big)^{k_0}\int_{\Omega_{\frac{3}{2}R}}\left|
\nabla{\bf w}\right|^2+C\delta(z')^{d+2\alpha} \le C\delta(z')^{d+2\alpha},
\end{align*}
for sufficiently small $\epsilon$. This completes the proof.
\end{proof}

\begin{proof}[Proof of Proposition \ref{keyprop3}]
Combining Lemma \ref{w_L2_energy} and Proposition \ref{w,q-L_infty-esti}, we 
finish the proof of Proposition \ref{keyprop3}.
\end{proof}

\section*{Declarations}

\noindent{\bf Data availability.} Data sharing not applicable to this article as no data sets were generated or analyzed during the current study.

\noindent{\bf Financial interests.} The authors have no relevant financial or non-financial interests to disclose.

\noindent{\bf Conflict of interest.} The authors declare that they have no conflict of interest.

\end{document}